\documentclass[11pt, twoside]{amsart}

\usepackage[T1]{fontenc}
\usepackage{amssymb,amsmath,amstext}

\newtheorem{theor}{Theorem}[section]
\newtheorem{lem}[theor]{Lemma}
\newtheorem{defin}[theor]{Definition}

\newtheorem{prop}[theor]{Proposition} 

\newtheorem{notation}[theor]{Notation}
\newtheorem{exam}[theor]{Example}

\newtheorem{cor}[theor]{Corollary}
\newtheorem{rem}[theor]{Remark}

\newtheorem{assump}[theor]{Assumption}

\newtheorem{terminologyandnotation}[theor]{Terminology and notation}

\numberwithin{equation}{section}

\newcommand{\es}{\emptyset}

\newcommand{\nts}{\negthickspace}
\newcommand{\uhrc}{\nts \upharpoonright \nts}

\newcommand{\mcA}{\mathcal{A}}
\newcommand{\mcB}{\mathcal{B}}
\newcommand{\mcC}{\mathcal{C}}

\newcommand{\mcM}{\mathcal{M}}
\newcommand{\mcN}{\mathcal{N}}
\newcommand{\mcO}{\mathcal{O}}

\newcommand{\mcS}{\mathcal{S}}

\newcommand{\mbC}{\mathbf{C}}

\newcommand{\mbK}{\mathbf{K}}
\newcommand{\mbS}{\mathbf{S}}
\newcommand{\mbT}{\mathbf{T}}
\newcommand{\mbU}{\mathbf{U}}

\newcommand{\mbX}{\mathbf{X}}
\newcommand{\mbY}{\mathbf{Y}}

\newcommand{\mbbN}{\mathbb{N}}
\newcommand{\mbbR}{\mathbb{R}}

\newcommand{\mbbQ}{\mathbb{Q}}

\newcommand{\Aut}{\mathrm{Aut}}

\newcommand{\Spt}{\mathrm{Spt}}
\newcommand{\spt}{\mathrm{spt}}
\newcommand{\Orb}{\mathrm{Orb}}
\newcommand{\orb}{\mathrm{orb}}

\newcommand{\liso}{\big|\big[}
\newcommand{\riso}{\big]\big|}

\title[Limit laws of random nonrigid structures]
{Limit laws and automorphism groups \\ of random nonrigid structures}

\author{Ove Ahlman and Vera Koponen}

\address{Ove Ahlman, Department of Mathematics, Uppsala University, Box 480,
75106 Uppsala, Sweden.}

\email{ove@math.uu.se}

\address{Vera Koponen, Department of Mathematics, Uppsala University, Box 480,
75106 Uppsala, Sweden.}

\email{vera@math.uu.se}

\begin{document}

\begin{abstract} 
A systematic study is made, for an arbitrary finite relational language with at least one symbol of
arity at least 2, of classes of nonrigid finite structures. 
The well known results that almost all finite structures are rigid and that
the class of finite structures has a zero-one law are, in the present context, 
the first layer in a hierarchy of classes of finite structures with increasingly more complex 
automorphism groups.
Such a hierarchy can be defined in more than one way. 
For example, the $k$th level of the hierarchy can consist of all structures 
having at least $k$ elements which are moved by some automorphism.
Or we can consider, for any finite group $G$, all finite structures $\mcM$
such that $G$ is a subgroup of the group of autmorphisms of $\mcM$;
in this case the ``hierarchy'' is a partial order.
In both cases, as well as variants of them, each ``level'' satisfies a logical limit law,
but not a zero-one law (unless $k = 0$ or $G$ is trivial).
Moreover, the number of (labelled or unlabelled) $n$-element structures in one place of the hierarchy divided
by the number of $n$-element structures in another place always converges to a rational number or to $\infty$
as $n \to \infty$.
All instances of the respective result are proved by an essentially uniform argument.

\end{abstract}

\maketitle

\begin{center}{\large Contents}\end{center}
{\small
\contentsline{section}{\ref{introduction}. Introduction}{\pageref{introduction}}
\contentsline{section}{\ref{Upper bounds of the support of automorphisms}. 
Upper bounds of the support of automorphisms}{\pageref{Upper bounds of the support of automorphisms}}
\contentsline{section}{\ref{Upper bounds of the support of structures}. 
Upper bounds of the support of structures}{\pageref{Upper bounds of the support of structures}}
\contentsline{section}{\ref{asymptotic estimates}. 
Asymptotic estimates of the number of structures with bounded support}{\pageref{asymptotic estimates}}
\contentsline{section}{\ref{Comparing different groups}.
Comparing different groups}{\pageref{Comparing different groups}}
\contentsline{section}{\ref{Logical limit laws}.
Logical limit laws}{\pageref{Logical limit laws}}
\contentsline{section}{\ref{Unlabelled structures}.
Unlabelled structures}{\pageref{Unlabelled structures}}
\contentsline{section}{References}{\pageref{References}}
}

\section{Introduction}\label{introduction}

\noindent
In a sequence of articles
\cite{ER, Fag77, FU, Har, Ober} it has been shown that for any finite relational vocabulary
(also called signature), the proportion of labelled (as well as unlabelled) $n$-element 
structures which are {\em rigid}, i.e. has no nontrivial automorphism, approaches 1 as 
$n$ approaches infinity.
By the work of Glebskii et. al. \cite{Gleb} and Fagin \cite{Fag76}, 
for any sentence $\varphi$, the proportion of $n$-element structures 
(labelled or unlabelled) in which $\varphi$ is true approaches either 0 or 1 as
$n$ tends to infinity. In other words, the class of finite structures satisfy a
(labelled and unlabelled) zero-one law.

However, the asymptotic behaviour of nonrigid $n$-element structures appears to have been neglected,
besides work of Cameron \cite{Cam80, Cam05} in the case of unlabelled undirected graphs. 
Possibly because the class of nonrigid finite structures
make up to only a ``measure zero'' subclass of the class of all finite structures.
Nevertheless, for any integer $k$, the number of (nonisomorphic) $n$-element structures with at least $k$ elements
which are moved by some automorphism grows exponentially with $n$;
and the same holds for the number of $n$-element structures whose automorphism group contains some specified group.
(This follows from the proofs in Section~\ref{Upper bounds of the support of automorphisms}.)
But more interestingly, consideration of finite structures whose automorphism group has
a certain (minimum) complexity gives rise to an infinitude of natural classes of finite (nonrigid)
structures with logical limit laws (Theorem~\ref{main theorem about limit laws for first-order logic}).
Each such class has the property that there are more than one but only finitely many ``convergence points'',
all of which are rational; that is, there is a finite set $A$ of rational numbers such that $|A| > 1$
and, for every sentence $\varphi$, the proportion of $n$-element structures in the class which satisfy $\varphi$
converges to a number in $A$. Moreover, in a sense that can be made precise, there are only
finitely many (but more than one) ``limit theories'' of any such class, 
all of which are $\aleph_0$-categorical and simple with
SU-rank one.\footnote{
For any finite relational language with at least one symbol of arity at least 2 and integer $l \geq 2$,
the class of all finite structures and the class of all (strongly) $l$-colourable finite structures 
\cite{KPR, Kop12b} have a zero-one law with a ``limit/almost sure'' theory
which is $\aleph_0$-categorical and simple with SU-rank 1.
The class of all finite partial orders has a zero-one law \cite{Com} with a limit theory which
is probably $\aleph_0$-categorical (because the ``height'' of a finite partial order is almost always 3
\cite{KR}), although we have not checked this.}
It appears like the classes of nonrigid structures considered here are the first 
nontrivial and ``naturally occuring'' classes of finite structures with such limit law behaviour.\footnote{
A trivial example can constructed by adding a new unary relation symbol $R$ to a vocabulary with some
relation symbol of arity at least 2 and letting the interpretation of $R$ be
a singleton set in half of all $n$-element structures in the initial vocabulary and the empty set in the other half.}

Furthermore, for any two classes $\mbC$ and $\mbK$ of finite structures that are
associated with some (minimum) complexity
of the automorphism group, the number of (labelled or unlabelled) $n$-element structures
which belong to $\mbC$ divided by the number of $n$-element structures which belong to $\mbK$
converges to a rational number or to $\infty$ as $n \to \infty$ 
(Theorem~\ref{main theorem about comparissons between different groups}
and Remark~\ref{remark on limit for other complexity measures than subgroups}). 

In general, this study gives fairly complete answers, for any finite relational
vocabulary with at least one relation symbol with arity at least 2 and for labelled as
well as unlabelled structures, to questions initiated by Cameron long ago
(in particular Theorems~1 and~2 in \cite{Cam80}),
but also to other natural variations of his questions and to the problem of whether
logical limit laws hold for classes of structures defined in terms of the complexity of their automorphism group.

A more detailed study, for any $m \in \mbbN$, of the typical autmorphism groups of finite 
structures such that at least $m$ elements are moved by some automorphism is carried out
in \cite{Kop13c}. Roughly speaking, \cite{Kop13c} shows that almost all finite structures
with some minimum complexity of their automorphism group have as simple automorphism group as
the minimum complexity allows.

Before stating the main results we introduce some basic terminology, 
notation and assumptions that will be used throughout.
We fix a finite {\em vocabulary}, also called {\em signature}, 
$\{R_1, \ldots, R_\rho\}$ of (only) {\em relation symbols} where $R_i$ has arity $r_i$.
Let $r = \max\{ r_1, \ldots, r_\rho \}$
and we call $r$ the {\em maximal arity}. 
{\em We always assume that $r \geq 2$,}
although this assumption is sometimes repeated.
By a {\em structure}, we mean a structure for the above vocabulary, that is, a tuple 
$\mcM = (M, R_1^\mcM, \ldots, R_\rho^\mcM)$ where
$M$ is a set, called the {\em universe} of $\mcM$, and, 
for each $i = 1, \ldots, \rho$, $R_i^\mcM \subseteq M^{r_i}$.
The relation $R_i^{\mcM}$ is called the {\em interpretation} of $R_i$ in $\mcM$.
For every positive integer $n$ let $[n] = \{1, \ldots, n\}$ and let 
$\mbS_n$ be the set of all structures with universe $[n]$ and let
$\mbS = \bigcup_{n=1}^\infty \mbS_n$. 
For every structure $\mcM$, let $\Aut(\mcM)$ denote the group of automorphisms of $\mcM$.
(For basic model theory, see \cite{Mar, Roth} for example.)

For groups $G$ and $H$, $G \cong H$ means that they are isomorphic (as abstract groups) and
$G \leq H$ means that $G$ is isomorphic to a subgroup of $H$. 
For structures $\mcM$ and $\mcN$, $\mcM \cong \mcN$ means that they are isomorphic.
Let $\mbbN$, $\mbbN^+$, $\mbbQ$ and $\mbbR$ denote the sets of nonnegative integers, 
positive integers, rational and real numbers, respectively.

\begin{theor}\label{main theorem about comparissons between different groups}
For any two finite groups $G$ and $H$, each one of the following limits exists in $\mbbQ \cup \{\infty\}$:
\begin{align*}
&\lim_{n\to\infty} 
\frac{\big| \{\mcM \in \mbS_n : H \leq \Aut(\mcM)\} \big|}{\big| \{\mcM \in \mbS_n : G \leq \Aut(\mcM)\} \big|}, \quad
\lim_{n\to\infty} 
\frac{\big| \{\mcM \in \mbS_n : H \cong \Aut(\mcM)\} \big|}{\big| \{\mcM \in \mbS_n : G \cong \Aut(\mcM)\} \big|} 
\quad \text{ and} \\
&\lim_{n\to\infty} 
\frac{\big| \{\mcM \in \mbS_n : G \cong \Aut(\mcM)\} \big|}{\big| \{\mcM \in \mbS_n : G \leq \Aut(\mcM)\} \big|}.
\end{align*}
\end{theor}

\noindent
Before stating the remaining main results, we introduce some more notation which will be used
throughout the article. For a set $A$, $|A|$ denotes its cardinality and $Sym(A)$ denotes
the group of all permutations of $A$. 
If $f_1, \ldots, f_k \in Sym(A)$ then $\langle f_1, \ldots, f_k \rangle$ denotes the subgroup
of $Sym(A)$ generated by $f_1, \ldots, f_k$,
\[ 
\Spt(f_1, \ldots, f_k) \ = \ \{a \in A : g(a) \neq a \text{ for some } g \in \langle f_1, \ldots, f_k \rangle \}
\]
and $\spt(f_1, \ldots, f_k) = |\Spt(f_1, \ldots, f_k)|$.
We call $\Spt(f_1, \ldots, f_k)$ the {\em support} of the sequence $f_1, \ldots, f_k$.
For a finite structure $\mcM$ we let
\begin{align*}
\spt(\mcM) \ &= \ \max\{\spt(f) : f \in \Aut(\mcM)\},\\
\Spt^*(\mcM) \ &= \ \{a \in M : a \in \Spt(f) \text{ for some } f \in \Aut(\mcM)\}, \ \text{ and}\\
\spt^*(\mcM) \ &= \ \big|\Spt^*(\mcM)\big|.
\end{align*}
The set $\Spt^*(\mcM)$ is called the {\em support} of $\mcM$. 
Note that we always have $\spt(\mcM) \leq \spt^*(\mcM)$.
Throughout, we use the following notation for $p, p' \in \mbbN$:
\begin{align*}
\mbS_n(\spt = p) \ &= \ \{\mcM \in \mbS_n : \spt(\mcM) = p\},\\
\mbS_n(\spt \geq p) \ &= \ \{\mcM \in \mbS_n : \spt(\mcM) \geq p\},\\
\mbS_n(\spt \leq p) \ &= \ \{\mcM \in \mbS_n : \spt(\mcM) \leq p\},\\
\mbS_n(\spt^* = p) \ &= \ \{\mcM \in \mbS_n : \spt^*(\mcM) = p \},\\
\mbS_n(\spt^* \geq p) \ &= \ \{\mcM \in \mbS_n : \spt^*(\mcM) \geq p \},\\
\mbS_n(\spt^* \leq p) \ &= \ \{\mcM \in \mbS_n : \spt^*(\mcM) \leq p \},\\
\mbS_n(p \leq \spt^* \leq p') \ &= \ \{\mcM \in \mbS_n : p \leq \spt^*(\mcM) \leq p' \}.
\end{align*}
Whenever $\mbS'_n \subseteq \mbS_n$ is defined for $n \in \mbbN^+$ we let 
$\mbS' = \bigcup_{n=1}^{\infty} \mbS'_n$.
The expression {\em almost all $\mcM \in \mbS'$ has property $P$} means that
\[
\lim_{n\to\infty} \ \frac{\big| \{ \mcM \in \mbS'_n : \text{ $\mcM$ has $P$} \}\big|}{\big| \mbS'_n \big|}
\ = \ 1.
\]

\noindent
Suppose that $\mbS'_n \subseteq \mbS_n$ for all $n \in \mbbN^+$.
We say that $\mbS' = \bigcup_{n \in \mbbN^+} \mbS'_n$ 
has a {\em limit law} if for every first-order sentence $\varphi$
over the vocabulary, the proportion of $\mcM \in \mbS'_n$ which satisfy $\varphi$ converges as $n \to \infty$.
If the limit converges to 0 or 1 for every first-order sentence $\varphi$, then we say that $\mbS'$
has a {\em zero-one law}.

\begin{theor}\label{main theorem about limit laws for first-order logic}
(i) For every finite group $G$, 
$\{\mcM \in \mbS : G \cong \Aut(\mcM)\}$ and \\
$\{\mcM \in \mbS : G \leq \Aut(\mcM)\}$ have a limit law.\\
(ii) For every integer $m \geq 2$, 
$\mbS(\spt^* = m)$, $\mbS(\spt \geq m)$ and $\mbS(\spt^* \geq m)$ have a limit law.\\
(iii) In each case of the previous parts there is a finite set $A \subseteq \mbbQ$ such that, for 
every first-order sentence $\varphi$, the proportion of structures of the kind considered
converges to some $a \in A$ as $n \to \infty$.
\end{theor}

\noindent
However, in each case of 
Theorem~\ref{main theorem about limit laws for first-order logic}
we do not we have a zero-one law, if $G$ is nontrivial,
as explained in Remark~\ref{remark on limit laws}.

\begin{theor}\label{main theorem about unlabelled structures}
Theorems~\ref{main theorem about comparissons between different groups}
and~\ref{main theorem about limit laws for first-order logic}.
also hold in the unlabelled case, that is, if we only count structures up to isomorphism.
\end{theor}

\begin{rem}\label{remark about asymptotic formulas}{\rm ({\bf Asymptotic estimates})
The results, in particular Propositions~\ref{estimate of cardinality of S-n(A)}
and~\ref{general quotient proposition}
and Lemmas~\ref{spt-star equals a union of S(A,H)},
\ref{different sets equals a union of S(A,H)},
\ref{consequence of H and H' being equivalent over A}
and~\ref{moving out the summation sign} give, in principle, a
method of finding, for any finite group $G$, an asymptotic formula of
the number of $\mcM \in \mbS_n$ such that $G \leq \Aut(\mcM)$.
The same is true if `$\leq$' is replaced by `$\cong$' or if we instead consider, for some arbitrary fixed
integer $m \geq 2$,
$\big| \mbS_n(\spt \geq m) \big|$, $\big| \mbS_n(\spt^* \geq m) \big|$ or 
$\big| \mbS_n(\spt^* = m) \big|$ as $n \to \infty$.
}\end{rem}

\begin{rem}\label{remark about irreflexive and/or symmetric relations}{\rm
({\bf Irreflexive and symmetric relations})
(i) Suppose that every relation symbol is always interpreted as an {\em irreflexive}
relation, that is, if $\mcM \models R_i(a_1, \ldots, a_{r_i})$ then $a_j \neq a_{j'}$ whenever $j \neq j'$.
Then 
Theorems~\ref{main theorem about comparissons between different groups}~--
~\ref{main theorem about unlabelled structures}
remain true, but some modifications have to be made in some proofs
and in some technical results of the article. 

(ii) Suppose that every relation symbol is always interpreted as an irreflexive and 
{\em symmetric} relation, where the later means that
if $\mcM \models R_i(a_1, \ldots, a_{r_i})$ then 
$\mcM \models R_i(a_{\pi(1)}, \ldots, a_{\pi(r_i)})$ for every permutation $\pi$ of $[r_i]$.
Again Theorems~\ref{main theorem about comparissons between different groups}~--
~\ref{main theorem about unlabelled structures}
remain true, with minor modifications in some proofs
and technical results. 
}\end{rem}

\noindent
Here follows an outline of the article.
We deal with labelled structures until the last section, 
where we show why the main results also hold for unlabelled structures.
In Section~\ref{Upper bounds of the support of automorphisms} 
we show that for every $m \in \mbbN$ there is a number $t$, 
depending only on $m$ and the vocabulary, 
such that almost all $\mcM \in \mbS(\spt \geq m)$ have {\em no} automorphism the support of 
which contains more than $t$ elements.
In Section~\ref{Upper bounds of the support of structures} 
we show, by a Ramsey type argument, that if $\mcM$ is finite
and for every $f \in \Aut(\mcM)$, $\spt(f) \leq t$, then there are at most $t^{t+2}$ elements 
$a \in M$ such that $g(a) \neq a$ for some $g \in \Aut(\mcM)$.
More briefly, with the notation after Theorem~\ref{main theorem about comparissons between different groups}:
if $\spt(\mcM) \leq t$ then $\spt^*(\mcM) \leq t^{t+2}$.
A consequence of these results is that for every $m \in \mbbN$ there is $T \in \mbbN$ 
such that almost all $\mcM \in \mbS(\spt \geq m)$ have the property that 
at most $T$ elements are moved by some automorphism.
(In the case of unlabelled undirected graphs this was proved, in a different way, by Cameron \cite{Cam80}.)

Section~\ref{asymptotic estimates} considers asymptotic estimates that are needed later.
In this section, a structure $\mcA \in \mbS$ and subgroup $H$ of $\Aut(\mcA)$ is given
and an asymptotic estimate is proved for the number of $\mcM \in \mbS_n$ such that
$\spt^*(\mcM) = |A|$ and there is an embedding $f : \mcA \to \mcM$  such that
$H_f = \{f \sigma f^{-1} : \sigma \in H\}$ is a subgroup of the group
$\{g \uhrc \Spt^*(\mcM) : g \in \Aut(\mcM)\}$. The set of such structures is denoted $\mbS_n(\mcA, H)$.
Sets of this sort are the ``building blocks'' of other sets of structures considered here, 
in the sense that almost all
structures of any set of structures in the main theorems belong to a finite union of sets of the
form $\bigcup_{n=1}^{\infty}\mbS_n(\mcA, H)$.
In Section~\ref{Comparing different groups} we use the results from previous sections and in particular
the asymptotic estimate of $\mbS_n(\mcA, H)$ to prove 
Theorem~\ref{main theorem about comparissons between different groups}, 
in the form of Propositions~\ref{comparisson of two subgroups},
~\ref{proportion of structures with subgroup G that have automorphism group G}
and~\ref{comparisson of two groups as automorphism groups}.

Theorem~\ref{main theorem about limit laws for first-order logic}, 
about logical limit laws, is proved in Section~\ref{Logical limit laws}.
Again, the set $\mbS_n(\mcA, H)$ plays a central role.
In fact, the main task is to prove that $\mbS(\mcA, H)$ has a zero-one
law. This and Proposition~\ref{general quotient proposition}
implies Theorem~\ref{main theorem about limit laws for first-order logic}.
The final Section~\ref{Unlabelled structures} shows why all main results also hold
for unlabelled structures. This is summarised in 
Theorem~\ref{unlabelled limit laws} which implies
Theorem~\ref{main theorem about unlabelled structures}.

\begin{terminologyandnotation}{\rm
We use the calligraphic letters $\mcA, \mcB, \mcC, \mcM, \mcN$ to denote structures and the corresponding
noncalligraphic letters $A, B, C, M, N$ to denote their universes. 
Usually the universe will
be $[n] = \{1, \ldots, n\}$ for some $n \in \mbbN^+$. 
We sometimes write $\bar{a}$ to denote a finite tuple $(a_1, \ldots, a_n)$,
and if $\bar{a} = (a_1, \ldots, a_n)$ and $\bar{b} = (b_1, \ldots, b_m)$, 
then we let $\bar{a}\bar{b} = (a_1, \ldots, a_n, b_1, \ldots, b_m)$.
If $\mcM$ is a structure and $A \subseteq M$, then $\mcM \uhrc A$ denotes the substructure of $\mcM$ with universe $A$.

Let $H$ and $H'$ be permutation groups on sets $\Omega$ and $\Omega'$, respectively.
A bijection $f : \Omega \to \Omega'$ is called an {\em isomorphism from $H$ to $H'$ as 
permutation groups} if
$H' = \{ fhf^{-1} : h \in H\}$.
We say that $H$ and $H'$ are {\em isomorphic as permutation groups} if such $f$ exists, and this clearly
implies that they are isomorphic as abstract groups.
We let $H \cong_P H'$ mean that $H$ and $H'$ are isomorphic as permutation groups.
If $f : A \to B$ is a function and $X \subseteq A$, then $f\uhrc X$ denotes the restriction of $f$ to $X$.
If $H$ is a permutation group on $\Omega$ and $X \subseteq \Omega$ is the union of some of the orbits
of $H$ on $\Omega$, then we define $H \uhrc X = \{h \uhrc X : h \in H\}$ which is a permutation group on $X$,
and we call $H \uhrc X$ the {\em restriction of $H$ to $X$}.

If $f$ is a permutation of $\Omega$ then $a \in \Omega$ is called a {\em fixed point of $f$}
if $f(a) = a$. If $H$ is a group of permutations of $\Omega$ then $a \in \Omega$ is called
a {\em fixed point of $H$} if $a$ is a fixed point of {\em every} $h \in H$.
For a structure $\mcA$, $a \in A$ is called a {\em fixed point of $\mcA$} if $a$ is a
fixed point of $\Aut(\mcA)$.
For any nonempty set $\Omega$, $Sym(\Omega)$ denotes the symmetric group of $\Omega$,
i.e. the group of all permuations of $\Omega$, and $Sym_n = Sym([n])$.

If $G$ is a group and $g_1, \ldots, g_n \in G$ then $\langle g_1, \ldots, g_n \rangle$
denotes the subgroup of $G$ generated by $g_1, \ldots, g_n$.
For a permutation group $G$ on a set $\Omega$ we let $\Orb(G)$ be the set of orbits of $G$ on $\Omega$
and $\orb(G) = \big| \Orb(G) \big|$. Such $G$ also acts on $\Omega^m$, the set of ordered $m$-tuples
of elements from $\Omega$, by 
the action $g(a_1, \ldots, a_m) = \big(g(a_1), \ldots, g(a_m)\big)$ for every 
$g \in G$ and $(a_1, \ldots, a_m) \in \Omega^m$. 
When refering to ``the orbits of $G$ on $\Omega^m$'' we mean the orbits with respect to this action,
unless something else is said.
We let $\Orb^m(G)$ be the set of orbits of $G$ on $\Omega^m$ and
$\orb^m(G) = \big|\Orb^m(G)\big|$. 
For $\pi_1, \ldots, \pi_k \in Sym_n$ we let 
$\Orb(\pi_1, \ldots, \pi_k) = \Orb(\langle \pi_1, \ldots, \pi_k \rangle)$,
$\orb(\pi_1, \ldots, \pi_k) = \orb(\langle \pi_1, \ldots, \pi_k \rangle)$,
$\Orb^m(\pi_1, \ldots, \pi_k) = \Orb^m(\langle \pi_1, \ldots, \pi_k \rangle)$ and
$\orb^m(\pi_1, \ldots, \pi_k) = \orb^m(\langle \pi_1, \ldots, \pi_k \rangle)$.
For unexplained notions such as `action, orbit' etc., see for example \cite{DM}.

We will also use the terminology and notation that was introduced between 
Theorems~\ref{main theorem about comparissons between different groups} 
--~\ref{main theorem about unlabelled structures} as well as the following notation:
if $f_1, \ldots, f_k$ are permutations of $[n]$, then 
\[
\mbS_n(f_1, \ldots, f_k) \ = \ \{\mcM \in \mbS_n : f_1, \ldots, f_k \in \Aut(\mcM)\}.
\]
By $f(n) \sim g(n)$ (as $n \to \infty$) we mean that $f(n)/g(n) \to 1$ as $n \to \infty$.
It will be convenient to use the notation $\exp_2(x) = 2^x$.
}\end{terminologyandnotation}

\section{Upper bounds of the support of automorphisms}\label{Upper bounds of the support of automorphisms}

\noindent
The main result of this section, Proposition~\ref{automorphisms with large support are unusual}, 
is that for any $m \in \mbbN$ there is $t \in \mbbN$ such that the 
proportion of $\mcM \in \mbS_n(\spt \geq m)$ such that $\spt(\mcM) \leq t$ approaches 1 as $n \to \infty$.
We also derive a couple of corollaries of this which are important for the rest of the article.
The following elementary result, often called Burnside's lemma or theorem\footnote{
But was actually proved earlier by Cauchy and Frobenius, according to~\cite{DM}}, will be used. Proofs are found in \cite{Burn, DM}, for example.

\begin{prop}\label{orbitformula}
If $G$ is a group of permutations of a finite set $M$ then
\[\orb(G) \ = \ \frac{1}{|G|} \sum_{g\in G} \big|\{a\in M : g(a) = a\}\big| \]
\end{prop}

\noindent
Recall that $[n] = \{1, \ldots,n\}$ and by $[n]^r$ we denote the set of {\em ordered} 
$r$-tuples of elements from $[n]$. 

\begin{lem}\label{orbitlemma}
Suppose that $d, n \in \mbbN^+$, $\pi_1, \ldots, \pi_s \in Sym_n$ 
and $\spt(\pi_1, \ldots, \pi_s) = p$.
Then
\[ \frac{n^d \ + \ (p! - 1)(n-p)^d}{p!} \ \leq \ 
\orb^d\big(\pi_1, \ldots, \pi_s\big) \ \leq \ n^d \ - \ \frac{pn^{d-1}}{2}.\] 
\end{lem}

\noindent
{\bf Proof.}
Suppose that $d,n\in \mbbN^+$, $\pi_1, \ldots, \pi_s \in Sym_n$ and $\spt(\pi_1, \ldots, \pi_s) = p$.
For each $\pi \in Sym_n$ let $\widetilde \pi \in Sym([n]^d)$ be defined by
$\widetilde \pi(x_1, \ldots, x_d) = (\pi(x_1), \ldots, \pi(x_d))$.
We consider the subgroup $G = \langle \pi_1, \ldots, \pi_s \rangle$ of $Sym_n$ and the 
subgroup $\widetilde G = \langle \widetilde\pi_1, \ldots, \widetilde\pi_s \rangle$ of $Sym([n]^d)$.
Note that $\orb(\widetilde G) = \orb^d(G) = \orb^d\big(\pi_1, \ldots, \pi_s\big)$. 
Also observe that, by the assumption
that $\spt(G) = \spt(\pi_1, \ldots, \pi_s) = p$, every $g \in G$ has at least $n - p$ fixed points.
Therefore every $g \in \widetilde G$ has at least $(n-p)^d$ fixed points.
In particular, the identity permutation has $n^d$ fixed points.
Therefore we get, by also using Proposition~\ref{orbitformula},
\begin{align*}
\orb^d(\pi_1, \ldots, \pi_s) \ &= \ \orb(\widetilde G)
\ = \ 
\frac{1}{\big|\widetilde G\big|}
\sum_{\widetilde g \in \widetilde G} \big|\{\bar a \in [n]^d : \widetilde g(\bar a) = \bar a\}\big| \\
&\geq \ 
\frac{(\big|\widetilde G\big| - 1)(n-p)^d \ + \ n^d}{\big|\widetilde G\big|} 
= \ 
(n-p)^d \ + \ \frac{n^d - (n-p)^d}{\big|\widetilde G\big|}\\
&\geq \
(n-p)^d \ + \ \frac{n^d - (n-p)^d}{p!}
\ = \ 
\frac{n^d \ + \ (p! - 1)(n-p)^d}{p!}
\end{align*}

\noindent
On the other hand we also have that
\[\orb(G) \ \leq \ (n - p)  +  \frac{p}{2} \ = \ n - \frac{p}{2}\]
which implies that
\[\orb^d\big(\pi_1, \ldots, \pi_s\big) \ = \ \orb(\widetilde G) \ \leq \ 
\orb(G) \cdot n^{d-1} \ \leq \ n^d - \frac{pn^{d-1}}{2},\]
because if $(a_1, \ldots, a_d)$ and $(b_1, \ldots, b_d)$ belong to the same
orbit of $\widetilde G$, then $a_1$ and $b_1$ belong to the same orbit of $G$.
\hfill $\square$
\\

\noindent
Recall that $r \geq 2$ is the maximal arity among relation symbols in the vocabulary.

\begin{prop}\label{automorphisms with large support are unusual}
Suppose that $m, t \in \mbbN$, $f_1, \ldots, f_s \in Sym_n$
and $\spt(f_1, \ldots, f_s) = m$.
For all sufficiently large $n$ the following holds, where $k$ is the number
of $r$-ary relation symbols and the bound $\mcO( \ )$ depends only on $m, t$ and the
vocabulary:
\begin{equation*}
\frac{\big|\mbS_n(\spt \geq t)\big|}{\big|\mbS_n(f_1, \ldots, f_s)\big|} \ \leq \ 
\exp_2\Bigg(
k\frac{\big(2(m! - 1)rm \ - \ (t-1)m!\big)n^{r-1}}{2(m!)} \ \pm \ \mcO\big(n^{r-2}\big) \Bigg).
\end{equation*}
Hence, if $t > 2r(m! - 1)m/m! + 1$ then the quotient approaches 0 as $n \to \infty$.
\end{prop}

\noindent
{\bf Proof.} 
For each $i = 1, \ldots, r$, let $k_i$ be the number of $i$-ary relation symbols.
Suppose that $m, t \in \mbbN$, $f_1, \ldots, f_s \in Sym_n$
and $\spt(f_1, \ldots, f_s) = m$. 
Observe that for every $i$ and every $i$-ary relation symbol $R$
we have: if $\bar a, \bar b \in [n]^i$ belong to the same orbit of 
$\langle f_1, \ldots, f_s\rangle$ and 
$\mcM \in \mbS_n(f_1, \ldots, f_s)$, then $\mcM \models R(\bar a)$
if and only if $\mcM \models R(\bar b)$.
Since this is the only restriction on members of $\mbS_n(f_1, \ldots,f_s)$ we get

\begin{equation}\label{equality for automorphisms}
\big|\mbS_n(f_1, \ldots, f_s)\big| \ = \ 
\exp_2\Bigg(\sum_{i=1}^r k_i
\orb^i\big(f_1, \ldots, f_s\big)\Bigg).
\end{equation}
For every $\pi \in Sym_n$ with $\spt(\pi) \geq t$ we have 
$\mbS_n(\pi) \subseteq \mbS_n(\spt \geq t)$ and therefore
\begin{equation}\label{upper bound for spt less or equal to m}
\big|\mbS_n(\spt \geq t)\big| \ \leq \ 
\sum_{\substack{\pi \in Sym_n \\ \spt(\pi) \geq t}}  |\mbS_n(\pi)| \ = \ 
\sum_{\substack{\pi \in Sym_n \\ \spt(\pi) \geq t}}
\exp_2\Bigg( \sum_{i=1}^r k_i \orb^i(\pi)\Bigg).  
\end{equation}
By first applying Lemma~\ref{orbitlemma} 
on $f_1, \ldots, f_s$ and then on an arbitrary $\pi \in Sym_n$ we get, for each $i = 1, \ldots, r$,

\begin{equation}\label{first lower bound on number of orbits}
\frac{n^i+(m! - 1)(n-m)^i}{m!} \ \leq \ \orb^i\big(f_1, \ldots, f_s\big), \quad \text{ and}
\end{equation}
\begin{equation}\label{first upper bound on number of orbits}
\orb^i(\pi) \ \leq \ n^i - \frac{n^{i-1}\spt(\pi)}{2} \quad \text{ for every $\pi \in Sym_n$}.
\end{equation}
A straightforward computation\footnote{
Set $a = \exp_2\Big( \log_2 n \ - \ \sum_{i=1}^r k_i \frac{n^{i-1}}{2} \Big)$
and we have $\sum_{j=t}^n a^j \leq a^t/(1-a) \leq a^{t-1}$ if $n$ is large enough.}
shows that for all sufficiently large $n$
\begin{align}\label{a convergent series}
&\sum_{j=t}^n \exp_2\Bigg( j \log_2 n \ - \ j \sum_{i=1}^r k_i \frac{n^{i-1}}{2} \Bigg) \\ 
\leq \ 
&\exp_2\Bigg( -k_r\frac{(t-1)n^{r-1}}{2} \ \pm \ \mcO\big(n^{r-2} + \log_2 n\big) \Bigg),
\nonumber
\end{align}
where the bound $\mcO( \ )$ depends only on the vocabulary.
Notice that the number of $\pi\in Sym_n$ with $\spt(\pi) = j$ is 
$\displaystyle \binom{n}{j} j! \leq n^j$. 
By also using~(\ref{equality for automorphisms})--(\ref{a convergent series}) we now get

\begin{align*}
&\frac{\big|\mbS_n(\spt \geq t)\big|}{\big|\mbS_n(f_1, \ldots, f_s)\big|} 
\ \leq \  \underset{\spt(\pi) \geq t}{\sum_{\pi\in Sym_n}}  
\exp_2\Bigg( \sum_{i=1}^r k_i \orb^i(\pi) \ - \ \sum_{i=1}^r k_i \orb^i(f_1, \ldots, f_s)\Bigg) \\
&\leq \ 
\underset{\spt(\pi) \geq t}{\sum_{\pi\in Sym_n}} 
\exp_2\Bigg( \sum_{i=1}^r 
k_i \bigg[n^i \ - \ \frac{n^{i-1} \spt(\pi)}{2}\bigg]  
\ - \ 
\sum_{i=1}^r k_i\frac{n^i + (m! - 1)(n - m)^i}{m!}\Bigg) \\
&\leq \ 
\sum_{j=t}^n n^{j}
\exp_2\Bigg(
\sum_{i=1}^r k_i\bigg[n^i \ - \ \frac{jn^{i-1}}{2}\bigg] 
\ - \
\sum_{i=1}^r k_i\frac{n^i + (m! - 1)(n - m)^i}{m!} \Bigg) \\
&= \ 
\exp_2\Bigg( 
\sum_{i=1}^r k_i \bigg[ n^i \ - \ \frac{n^i + (m! - 1)(n - m)^i}{m!}\bigg] \Bigg) 
\sum_{j=t}^n \exp_2\Bigg(j \log_2 n - j\sum_{i=1}^r k_i \frac{n^{i-1}}{2}\Bigg) \\
&\leq \ 
\exp_2\Bigg( k_1\frac{(m! - 1)m}{m!} \ + \ 
\sum_{i=2}^r k_i\frac{(m! - 1)imn^{i-1} \pm \mcO\big(n^{i-2}\big) }{m!} \Bigg) \cdot \\
& \hspace{50mm} 
\exp_2\Bigg( -k_r\frac{(t-1)n^{r-1}}{2} \ \pm \ \mcO\big(n^{r-2} + \log_2 n\big) \Bigg) \\
&= \ 
\exp_2\Bigg(
k_r\frac{\big(2(m! - 1)rm \ - \ (t-1)m!\big)n^{r-1}}{2(m!)} \ \pm \ 
\mcO\big(n^{r-2} + \log_2 n\big) \Bigg).
\end{align*}
\hfill $\square$

\begin{rem}{\rm 
Suppose that we require that a relation symbol $R_i$ of arity $r_i \geq 2$ is always interpreted as
an irreflexive and symmetric relation. Then we need to use a modification of 
Lemma~\ref{orbitlemma} where, for $\pi_1, \ldots, \pi_s \in Sym_n$, we consider
the orbits of $G = \langle \pi_1, \ldots, \pi_s \rangle$ on 
the set of $r_i$-subsets of $[n]$ by the action $g(\{a_1, \ldots, a_{r_i}\}) = \{g(a_1), \ldots, g(a_s)\}$ for
every $g \in G$ and $r_i$-subset $\{a_1, \ldots, a_{r_i}\} \subseteq [n]$.
By slightly modifying the proof of Lemma~\ref{orbitlemma} one gets that if $q$ is the
number of orbits of $G$ by its action on the set of $r_i$-subsets of $[n]$, then
\[
\frac{\binom{n}{d} \ - \ (p! - 1) \binom{n-p}{d}}{p!} \ \leq \ q \ \leq \ 
n \binom{n}{d} \ - \ \frac{p}{2}\binom{n}{d-1}.
\]
By using this when estimating (the appropriate analogues of) $\orb^i(\pi)$  and
$\orb^i(f_1, \ldots, f_s)$ in the proof of 
Proposition~\ref{automorphisms with large support are unusual} for each $i$-ary
relation symbol (where $i \geq 2$) that is always interpreted as an irreflexive and
symmetric relation, one gets a similar upper bound,
by a bit more involved computations. 
Similar adaptations work if we require that some relation symbols are always interpreted
as irreflexive, but not necessarily symmetric, relations.
}\end{rem}

\begin{cor}\label{corollary to automorphisms with large support are unusual}
Let $m \in \mbbN$. If $t > 2r(m! - 1)m/m! + 1$ then
\[
\lim_{n\to\infty} \ \frac{\big|\mbS_n(\spt \geq t)\big|}{\big|\mbS_n(\spt \geq m)\big|} \ = \ 
\lim_{n\to\infty} \ \frac{\big|\mbS_n(\spt \geq t)\big|}{\big|\mbS_n(\spt^* \geq m)\big|}  \ = \ 0.
\]
\end{cor}

\noindent
{\bf Proof.}
This follows immediately from Proposition~\ref{automorphisms with large support are unusual}, 
because if $f \in Sym_n$ and $\spt(f) = m$,
then $\mbS_n(f) \subseteq \mbS_n(\spt \geq m)
\subseteq \mbS_n(\spt^* \geq m)$.
\hfill $\square$

\begin{cor}\label{having G as a subgroup almost allways implies bounded support of automorphisms}
Suppose that $G$ is a finite group which is isomorphic to a group of permutations
of $[m]$. If $t = 2r(m! - 1)m/m! + 1$ then 
\[
\lim_{n \to \infty} 
\frac{\big| \{\mcM \in \mbS_n : G \leq \Aut(\mcM) \text{ and } \spt(\mcM) \leq t) \} \big|}
{\big| \{ \mcM \in \mbS_n : G \leq \Aut(\mcM) \} \big|} \ = \ 1.
\]
\end{cor}

\noindent
{\bf Proof.}
Let $H = \{h_1, \ldots, h_s\}$ be a permutation group on $[m]$ such that $H \cong G$. 
Let $t = 2r(m! - 1)m/m! + 1$.
Extend each $h_i$ to a permutation $h'_i$ of $[n]$
by letting $h'_i(j) = j$ for every $j > m$ and $h'_i(j) = h_i(j)$ for every $j \leq m$. 
Observe that for every $\mcM \in \mbS_n(h'_1, \ldots, h'_s)$, $G \leq \Aut(\mcM)$.
From Proposition~\ref{automorphisms with large support are unusual}
we get
\[
\frac{\big| \big\{ \mcM \in \mbS_n : 
G \leq \Aut(\mcM) \text{ and } \spt(\mcM) > t \big\} \big|}
{\big| \big\{ \mcM \in \mbS_n : G \leq \Aut(\mcM) \big\} \big|} 
\ \leq \ 
\frac{\big| \mbS_n(\spt > t) \big|}{\big| \mbS_n(h'_1, \ldots, h'_s) \big|} \ \to \ 0,
\]
as $n \to \infty$.
\hfill $\square$

\section{Upper bounds of the support of structures}\label{Upper bounds of the support of structures}

\noindent
In this section we prove that for every $t \in \mbbN$ there is $T \in \mbbN$, depending only on $t$,
such that for every finite structure $\mcM$, if $\spt(\mcM) \leq t$ then $\spt^*(\mcM) \leq T$. 
In other words, if no automorphism of $\mcM$ moves more than $t$ elements, 
then not more than $T$ elements of $\mcM$ are moved by some automorphism.
This is stated by Proposition~\ref{support of structure is bounded in terms of support of automorphisms}.
Corollaries~\ref{almost surely bounded support of automorphisms implies bounded support of structure}
and~\ref{having G as a subgroup almost always implies bounded support of the structure}
will be used in later sections.

\begin{defin}{\rm
Let $\mcM \in \mbS$ and $X\subseteq M$. \\
(i) For $f \in \Aut(\mcM)$ let $d(f,X) = |\Spt(f)- X|$.\\
(ii) We call $f \in \Aut(\mcM)$ {\em maximal} if 
for all $g\in \Aut(\mcM)$, if $\Spt(f)\subseteq \Spt(g)$ then $\Spt(f)=\Spt(g)$. \\
(iii) Let $\Aut^*(\mcM) = \{f\in\Aut(\mcM) :f \text{ is maximal} \}$.\\
(iv) For $\mcM\in \mbS$, a sequence $f_0, \ldots, f_n \subseteq\Aut^*(\mcM)$ is called a 
{\em special sequence of automorphisms} of $\mcM$ if it satisfies the two following condition:
\begin{itemize}

\item[] For each $k = 0, \ldots, n-1$, 
\[d\big(f_{k+1}, \Spt(f_0, \ldots, f_k)\big) \ = \ 
\max_{g\in\Aut^*(\mcM)} d\big(g, \Spt(f_0, \ldots, f_k)\big).\]
\end{itemize}
}\end{defin}

\begin{notation}{\rm
Whenever a special sequence of automorphisms $f_0, \ldots, f_n \in \Aut^*(\mcM)$, $k \leq n$
and $g \in \Aut(\mcM)$ are given, then we may use the abbreviation
\[d_k(g) \ = \ d\big(g, \Spt(f_0, \ldots, f_k)\big).\]
}\end{notation}

\noindent
The following lemma states some basic facts about special sequences of automorphisms.

\begin{lem}\label{dkfacts} 
Let $\mcM\in\mbS$ and let $f_0, \ldots, f_n \in \Aut^*(\mcM)$ be
a special sequence of automorphisms. Then
\begin{enumerate}
\item for all $0 \leq k \leq n$ and all $g\in\Aut(\mcM)$, $d_k(g)\geq d_{k+1}(g)$,
\item~\label{secoundfact} if $k+1 \leq p \leq n$ then $d_k(f_{k+1}) \geq d_k(f_p)$ and
\item if $0 \leq k < n$ and $d_k(f_{k+1}) = 0$ then for all $g\in\Aut^*(\mcM)$, 
$\Spt(g)\subseteq \Spt(f_1,\ldots,f_k)$.

\end{enumerate}
\end{lem}

\noindent
{\bf Proof.}
Let $\mcM\in\mbS$ and let $f_0, \ldots, f_n \in \Aut^*(\mcM)$ be
a special sequence of automorphisms.

(1) Suppose that $g \in \Aut(\mcM)$. 
As $\Spt(f_0, \ldots, f_k) \ \subseteq \ \Spt(f_0, \ldots, f_{k+1})$ we get
\[ \big| \Spt(g) \ \setminus \ \Spt(f_0, \ldots, f_k) \big| \ \geq \ 
\big| \Spt(g) \ \setminus \ \Spt(f_0, \ldots, f_{k+1}) \big|,\]
that is, $d_k(g) \geq d_{k+1}(g)$.

(2) Suppose that $k+1 \leq p \leq n$.
Since 
\[d_k(f_{k+1}) \ = \ \max_{g \in \Aut^*(\mcM)} d_k(g)\]
we get $d_k(f_{k+1}) \geq d_k(f_p)$.

(3) If $0 \leq k < n$ and $d_k(f_{k+1}) = 0$, then
$\max_{g\in \Aut^*(\mcM)} d_k(g) = 0$, so $\Spt(g) \subseteq \Spt(f_1, \ldots, f_k)$
for every $g \in \Aut^*(\mcM)$.
\hfill $\square$
\\

\noindent
Now to a less obvious claim:

\begin{lem}\label{distancelemma}
Let $\mcM \in \mbS$. Suppose that $f_0, \ldots, f_n \in \Aut^*(\mcM)$ is
a special sequence and $1 \leq k < p$.
If $d_k(f_{p})>0$ then there is 
$x \in \Spt(f_k) \setminus \Spt(f_0,\ldots,f_{k-1})$ such that $x\notin\Spt(f_p)$
\end{lem}

\noindent
{\bf Proof.}
Let $\mcM \in \mbS$, let $f_0, \ldots, f_n \in \Aut^*(\mcM)$ be
a special sequence and suppose that $1 \leq k < p$ and $d_k(f_p) > 0$.
We use the abbreviations $\Spt(k) = \Spt(f_k)$ and $\Spt(0, \ldots, k) = \Spt (f_0, \ldots, f_k)$.
Let 
\[X \ = \ \Spt(k) \ \setminus \ \Spt(0,\ldots,k-1).\]
For a contradiction, we assume that $X \subseteq \Spt(p)$.
Since $d_k(f_{p})>0$ we know that there is an element $a \in \Spt(p)$ such that $a\notin \Spt(0,\ldots,k)$. 
Then Lemma~\ref{dkfacts}~(1) together with $d_k(f_{p})>0$ gives us that $d_{k-1}(f_{p})>0$.
By Lemma~\ref{dkfacts}~(2) we get $d_{k-1}(f_k)>0$, which implies that $X \neq \es$. 
Also notice that $a\notin X$, by the choice of $a$.
From the definition of $X$ and the assumption that $X \subseteq \Spt(p)$
it follows that $X \subseteq \Spt(p) \setminus \Spt(0, \ldots, k-1)$.
By the choice of $a$ we have $a \in \Spt(p) \setminus \Spt(0, \ldots, k-1)$,
so we get 
\[X \cup \{a\} \ \subseteq \ \Spt(p) \ \setminus \ \Spt(0, \ldots, k-1),\]
and recall that $a \notin X$. Hence we get
\begin{align*}
d_{k-1}(k) \ &= \ \big|\Spt(k) - \Spt(0,\ldots,k-1)\big| \ = \ |X| \\ 
&< \ 
\big|X\cup\{a\}\big| \ \leq \ \big|\Spt(p) \ \setminus \ \Spt(0,\ldots,k-1)\big| \ = \ d_{k-1}(p),
\end{align*} 
i.e. $d_{k-1}(k) < d_{k-1}(p)$ which contradicts Lemma~\ref{dkfacts}~(\ref{secoundfact}).
\hfill $\square$
\\

\noindent
The next proposition tells that, for each $k \geq 2$,
$\mbS(\spt \leq k) \subseteq \mbS\big(\spt^* \leq k^{k+2}\big)$.

\begin{prop}\label{support of structure is bounded in terms of support of automorphisms}
For every integer $k \geq 2$ and every $\mcM \in \mbS(\spt \leq k)$
we have $\spt^*(\mcM) \leq k^{k+2}$. 
\end{prop}

\noindent
{\bf Proof.}
Fix any integer $k \geq 2$.
For $i = 0, \ldots, k$, let $l_i = k^{k-i+1}$.
Note that $l_0 = k^{k+1}$ and 
$l_i =  k l_{i+1}$ for each $i$.
Suppose that $\mcM \in \mbS(\spt \leq k)$ and, for a contradiction, that
$\spt^*(\mcM) > k^{k+2}$.

By definition, any $f_0 \in \Aut^*(\mcM)$ is a special sequence of length 1.
Now let $f_0, \ldots, f_n \in \Aut^*(\mcM)$ be any special sequence and
suppose that $n < l_0$. By the assumption that $\mcM \in \mbS(\spt \leq k)$
we have $\big|\Spt(f_0, \ldots, f_n)\big| \leq kl_0 = k^{k+2}$.
From the assumption that $\spt^*(\mcM) > k^{k+2}$ it now follows that there is 
$g \in \Aut(\mcM)$ such that $d_n(g) = \big|\Spt(g) \setminus \Spt(f_0, \ldots,f_n)\big| > 0$. 
Hence there is also a maximal $f \in \Aut^*(\mcM)$ such that
$d_n(f) > 0$. If we choose $f_{n+1} \in \Aut^*(\mcM)$
so that $d_n(f_{n+1}) = \max_{g \in \Aut^*(\mcM)} d_n(g)$,  
then $f_0, \ldots, f_{n+1}$ is a special sequence.
This proves that there is a special sequence
$f_0, \ldots, f_{l_0} \in \Aut^*(\mcM)$ such that $d_p(f_{p+1}) > 0$ for every
$p = 0, \ldots, l_0 - 1$. 
We fix this special sequence for the rest of the proof and use
the abbreviations $\Spt(p) = \Spt(f_p)$ and $\Spt(0, \ldots, p) = \Spt(f_0, \ldots, f_p)$.

We will prove that there are a subsequence (of distinct numbers) 
$t_1, \ldots, t_{k+1}$ of the sequence $0, \ldots, l_0$ and elements
$b_i \in \Spt(f_{t_i})$, for $i = 1, \ldots, k+1$, such that 
$b_i \notin \Spt(f_{t_j})$ if $j \neq i$; so $i \neq j$ implies $b_i \neq b_j$. 
Then $b_1, \ldots, b_{k+1} \in  \Spt(f_{t_1}\circ...\circ f_{t_{k+1}})$,
where of course the composition $f_{t_1}\circ...\circ f_{t_{k+1}}$ belongs to $\Aut(\mcM)$.
This contradicts the assumption that $\mcM \in \mbS(\spt \leq k)$.

We will inductively define sequences $t^i_0, \ldots, t^i_{l_i}$, for $i = 0, \ldots, k+1$,
of indices from which we can extract a sequence $t_1, \ldots, t_{k+1}$ as above.
Let $t^0_j = j$ for $j = 0, \ldots, l_0 = k^{k+1}$.
For each $p = 2, \ldots, l_0$, there is, by Lemma~\ref{distancelemma},
$a_p \in \Spt(1) \setminus \Spt(0)$ such that $a_p \notin \Spt(p)$.
As $|\Spt(1)| \leq k$ there are $b_1 \in \Spt(1) \setminus \Spt(0)$
and a subsequence of distinct numbers $t^1_1, \ldots, t^1_{l_1}$ of 
the sequence $2, \ldots, l_0$ such that, for all $p = t^1_1, \ldots, t^1_{l_1}$, $a_p = b_1$.
Let $t_1 = t^1_0 = 1$.

Now suppose that $m \leq k$ and that, for $i = 1, \ldots, m$, $t^i_0, \ldots, t^i_{l_i}$
is a subsequence (of distinct numbers) of $t^{i-1}_0, \ldots, t^{i-1}_{l_{i-1}}$,
$b_i \in \Spt(t^i_0) \setminus \Spt(0, \ldots, t^i_0 - 1)$ and
$b_i \notin \Spt(p)$ for all $p = t^i_1, \ldots, t^i_{l_i}$.
By Lemma~\ref{distancelemma}, there is for each $p = t^m_2, \ldots, t^m_{l_m}$
an element $a_p \in \Spt(t^m_1) \setminus \Spt(0, \ldots, t^m_1 - 1)$  such that
$a_p \notin \Spt(p)$. 
Since $\big|\Spt(t^m_1)\big| \leq k$ there are $b_{m+1} \in \Spt(t^m_1)$
and a subsequence $t^{m+1}_1, \ldots, t^{m+1}_{l_{m+1}}$ of 
$t^m_2, \ldots, t^m_{l_m}$ such that, for all $p = t^{m+1}_1, \ldots, t^{m+1}_{l_{m+1}}$,
$a_p = b_{m+1}$. 
Let $t_{m+1} = t^{m+1}_0 = t^m_1$.
When $t^i_0, \ldots, t^i_{l_i}$ are 
defined for every $i = 0, \ldots, k+1$
and $b_i$ for every $i = 1, \ldots, k+1$, then, 
as already indicated, we take $t_i = t^i_0$ for $i = 1, \ldots, k+1$.
\hfill $\square$

\begin{rem}{\rm
Notice that the proofs up to now of this section do not need the assumption
that we have considered a structure $\mcM$ and its automorphisms. We could, more generally, 
have considered a set $M$ and a group of permutations $H$ of $M$.
If we do this, we get the following version of 
Proposition~\ref{support of structure is bounded in terms of support of automorphisms}:
{\em If $k \geq 2$ is an integer and $H$ is a group of permutations of  a set $M$
such that $\spt(h) \leq k$ for every $h \in H$, then
\[ 
\big| \{ a \in M : h(a) \neq a \text{ for some } h \in H \} \big| \ \leq \ k^{k+2}.
\]
}}\end{rem}

\begin{cor}\label{almost surely bounded support of automorphisms implies bounded support of structure}
Let $m \in \mbbN$. If $k = 2r(m! - 1)m/m! + 1$ and $T = k^{k+2}$ then
\begin{align*}
&\lim_{n\to\infty} \frac{\big| \mbS_n(\spt \geq m) \ \cap \ \mbS_n(\spt^* \leq T) \big|}
{\big| \mbS_n(\spt \geq m) \big|} \ = \\ 
&\lim_{n\to\infty} \frac{\big| \mbS_n(\spt^* \geq m) \ \cap \ \mbS_n(\spt^* \leq T) \big|}
{\big| \mbS_n(\spt^* \geq m) \big|}
\ = \ 1.
\end{align*}
\end{cor}

\noindent
{\bf Proof.} 
Let $k = 2r(m! - 1)m/m! + 1$ and $T = k^{k+2}$.
By Corollary~\ref{corollary to automorphisms with large support are unusual},
\[\big| \mbS_n(\spt \geq m) \big| \ = \ \big(1 + o(1)\big)\big| \mbS_n(m \leq \spt \leq k)\big|\]
and by Proposition~\ref{support of structure is bounded in terms of support of automorphisms},
\[\big| \mbS_n(m \leq \spt \leq k) \big| \ = \ 
\big| \mbS_n(m \leq \spt \leq k) \ \cap \ \mbS_n(\spt^* \leq T) \big|,\]
so we get 
$\big| \mbS_n(\spt \geq m) \big| \ = \ \big(1 + o(1)\big)
\big| \mbS_n(m \leq \spt \leq k) \ \cap \ \mbS_n(\spt^* \leq T) \big|$.
The other limit is proved in the same way.
\hfill $\square$

\begin{cor}\label{having G as a subgroup almost always implies bounded support of the structure}
Suppose that $G$ is a finite group which is isomorphic to a group of permutations
of $[m]$ where $m \in \mbbN^+$.
 Then there is $T \in \mbbN$, depending only on $G$ and the vocabulary, such that 
\[
\lim_{n\to\infty} 
\frac{\big| \{ \mcM \in \mbS_n : G \leq \Aut(\mcM) \text{ and } \spt^*(\mcM) \leq T \} \big|}
{\big| \{\mcM \in \mbS_n : G \leq \Aut(\mcM) \} \big|} \ = \ 1. 
\]
\end{cor}

\noindent
{\bf Proof.}
By Corollary~\ref{having G as a subgroup almost allways implies bounded support of automorphisms} 
we know that if $k = 2r(m! - 1)m/m! + 1$ then
\[
\lim_{n \to \infty} 
\frac{\big| \{\mcM \in \mbS_n : G \leq \Aut(\mcM) \text{ and } \spt(\mcM) \leq k) \} \big|}
{\big| \{ \mcM \in \mbS_n : G \leq \Aut(\mcM) \} \big|} \ = \ 1.
\]
Let $T = k^{k+2}$.
As Proposition~\ref{support of structure is bounded in terms of support of automorphisms}
says that $\mbS_n(\spt \leq k) \subseteq \mbS_n(\spt^* \leq T)$ we are done.
\hfill $\square$

\section{Asymptotic estimates of the number of structures with bounded support}
\label{asymptotic estimates}

\noindent
By Corollary~\ref{almost surely bounded support of automorphisms implies bounded support of structure},
for arbitrary fixed $m \in \mbbN$ and all large enough $n$, 
an overwhelming part of the members of 
$\mbS_n(\spt \geq m)$ belong $\mbS_n(\spt^* \leq T)$ for some $T$ depending only on $m$ and the vocabulary.
We will show that an overwhelming part of the members of, for example,
$\mbS_n(\spt \geq m)$ for large enough $n$, belong to a
finite union of sets of the form $\mbS_n(\mcA, H)$, defined below, where the structure
$\mcA$ and permutation group $H$ depend only on the vocabulary and $m$.
In order to understand the asymptotic behaviour of 
$\mbS_n(\spt \geq m)$ we will therefore, in this section, find asymptotic estimates
of $\big| \mbS_n(\mcA, H) \big|$ as $n \to \infty$.
As will become clear in the sequel, the sets of the form $\mbS_n(\mcA, H)$
are the ``atomic'' pieces of our analysis, and questions about, for example, $\mbS_n(\spt \geq m)$
or $\{\mcM \in \mbS_n : G \leq \Aut(\mcM)\}$, for a fixed $G$,
will be reduced to analysing quotients of the form 
$\big|\mbS_n(\mcA', H')\big| \Big/ \big|\mbS_n(\mcA, H)\big|$ as $n \to \infty$.

Recall that if $H$ is a group of permutations of $\Omega$ and $X \subseteq \Omega$ 
is the union of some of the orbits of $H$ on $\Omega$, then $H \uhrc X = \{h \uhrc X : h \in H\}$ which is a permutation group on $X$. For every structure $\mcM$, $\Spt^*(\mcM)$ is the union of all nonsingleton 
orbits of $\Aut(\mcM)$ on $M$, 
so it always makes sense to speak about $\Aut(\mcM) \uhrc \Spt^*(\mcM)$ and we always have
$\Aut(\mcM) \uhrc \Spt^*(\mcM) \cong \Aut(\mcM)$.

\begin{defin}\label{definition of S-n(A)}{\rm
Let $\mcA \in \mbS$ be such that $\Aut(\mcA)$ has no fixed point.
Suppose that $H$ is a subgroup of $\Aut(\mcA)$ such that $H$ has no fixed point.
For each integer $n > 0$, $\mbS_n(\mcA, H)$ is the set of $\mcM \in \mbS_n$ 
such that there is an embedding $f : \mcA \to \mcM$ such that $\Spt^*(\mcM)$ is the
image of $f$ and 
$H_f = \{f\sigma f^{-1} : \sigma \in H\}$ is a subgroup of $\Aut(\mcM) \uhrc \Spt^*(\mcM)$.
Note that $H_f \cong_P H$. Let $\mbS(\mcA, H) = \bigcup_{n \in \mbbN^+} \mbS_n(\mcA, H)$.
}\end{defin}

\begin{lem}\label{spt-star equals a union of S(A,H)}
Let $m \geq 2$ be an integer. 
There are $\mcA_1, \ldots, \mcA_l \in \mbS_m$ without any fixed point and, for each 
$i = 1, \ldots, l$, subgroups $H_{i,1}, \ldots, H_{i,l_i} \subseteq \Aut(\mcA_i)$ 
without any fixed point such that
\[
\mbS(\spt^* = m) \ = \ \bigcup_{i=1}^l \bigcup_{j=1}^{l_i} \mbS(\mcA_i, H_{i,j}).
\]
\end{lem}

\noindent
{\bf Proof.}
Let $\mcA_1, \ldots, \mcA_l$ enumerate all structures of $\mbS_m$ that do not have any fixed point.
Suppose that $\mcM \in \mbS(\spt^* = m)$. 
Then $\mcM \uhrc \Spt^*(\mcM) \cong \mcA_i$ for some $i$.
If $K = \Aut(\mcM) \uhrc \Spt^*(\mcM)$, 
$f : \mcA_i \to \mcM \uhrc \Spt^*(\mcM)$ is an isomorphism
and $H = \{f^{-1}\sigma f : \sigma \in K\}$,
then $H$ is a subgroup of $\Aut(\mcA_i)$ without any fixed point.
From the definition of $\mbS(\mcA_i, H)$ it follows that $\mcM \in \mbS(\mcA_i, H)$.
Hence every $\mcM \in \mbS(\spt^* = m)$ belongs to $\mbS(\mcA_i, H)$ for some $i$ 
and some subgroup $H \subseteq \Aut(\mcA_i)$.
Conversely, for every $i = 1, \ldots, l$ and every subgroup
$H \subseteq \Aut(\mcA_i)$ we have $\mbS(\mcA_i, H) \subseteq \mbS(\spt^* = m)$,
since $\spt^*(\mcM) = m$ for every $\mcM \in \mbS(\mcA_i, H)$.
\hfill $\square$

\begin{lem}\label{different sets equals a union of S(A,H)}
(i) Let $m \geq 2$ be an integer.
There are finitely many $\mcA_1, \ldots, \mcA_l \in \mbS$ without any fixed point and, 
for each $i = 1, \ldots, l$,
subgroups $H_{i,1}, \ldots, H_{i,l_i} \subseteq \Aut(\mcA_i)$ without any fixed point such that
\[
\big| \mbS_n(\spt^* \geq m) \big| \ \sim \ \Bigg| \bigcup_{i=1}^l \bigcup_{j=1}^{l_i} \mbS_n(\mcA_i, H_{i,j}) \Bigg|
\quad \text{ as } n \to \infty.
\]
(ii) Part~(ii) holds if `$\spt^* \geq m$' is replaced by `$\spt \geq m$'.\\
(iii) Let $G$ be a nontrivial finite group.
There are finitely many $\mcA_1, \ldots, \mcA_l \in \mbS$ without any fixed point and, 
for each $i = 1, \ldots, l$,
subgroups $H_{i,1}, \ldots, H_{i,l_i} \subseteq \Aut(\mcA_i)$ without any fixed point such that
$G \leq H_{i,j}$ for all $i$ and $j$ and
\[
\big| \{ \mcM \in \mbS_n : G \leq \Aut(\mcM) \} \big| \ \sim \ 
\Bigg| \bigcup_{i=1}^l \bigcup_{j=1}^{l_i} \mbS_n(\mcA_i, H_{i,j}) \Bigg|
\quad \text{ as } n \to \infty.
\]
\end{lem}

\noindent
{\bf Proof.}
(i) By Corollary~\ref{almost surely bounded support of automorphisms implies bounded support of structure},
there is an integer $T$ such that 
\[
\big| \mbS_n(\spt^* \geq m) \big| \ \sim \ \big| \mbS_n(m \leq \spt^* \leq T) \big|
\quad \text{ as } n \to \infty.
\]
Since $\mbS_n(m \leq \spt^* \leq T) = \bigcup_{i=m}^T \mbS_n(\spt^* = m)$, part~(i)
follows from Lemma~\ref{spt-star equals a union of S(A,H)}.

(ii) By Corollary~\ref{almost surely bounded support of automorphisms implies bounded support of structure},
there is $T$ such that 
\[
\big| \mbS_n(\spt \geq m) \big| \ \sim \ \big| \mbS_n(\spt \geq m) \ \cap \ \mbS_n(\spt^* \leq T) \big|
\quad \text{ as } n \to \infty.
\]
As every $\mcM \in \mbS_n(\spt \geq m) \ \cap \ \mbS_n(\spt^* \leq T)$ belongs to 
$\mbS_n(\spt^* = p)$ for some $m \leq p \leq T$, we get~(ii) from 
Lemma~\ref{spt-star equals a union of S(A,H)}.

(iii) By Corollary~\ref{having G as a subgroup almost always implies bounded support of the structure},
there is an integer $T$ such that 
\[
\big| \{\mcM \in \mbS_n : G \leq \Aut(\mcM)\} \big| \ \sim \ 
\big| \{\mcM \in \mbS_n : G \leq \Aut(\mcM) \text{ and } \spt^*(\mcM) \leq T\} \big|
\]
as $n \to \infty$.
Since every $\mcM \in \{\mcM \in \mbS_n : G \leq \Aut(\mcM) \text{ and } \spt^*(\mcM) \leq T\}$
belongs to $\mbS_n(\spt^* = p)$ for some $p \leq T$, we also get part~(iii) from
Lemma~\ref{spt-star equals a union of S(A,H)} and its proof, which shows that
we only need to consider $\mcA_i$ and $H_{i,j}$ such that $G \leq H_{i,j}$.
\hfill $\square$
\\

\noindent
As suggested by the previous lemma, an essential step towards the main results is to asymptotically estimate
$\big| \mbS_n(\mcA, H) \big|$ for any $\mcA \in \mbS$ without a fixed point and any subgroup 
$H \subseteq \Aut(\mcA)$ without a fixed point.

\begin{prop}\label{estimate of cardinality of S-n(A)}
Suppose that $\mcA \in \mbS$ has no fixed point.
Let $H$ be a subgroup of $\Aut(\mcA)$ such that $H$ has no fixed point.
Let $p = |A|$, for every $i = 1, \ldots, r-1$
let $q_i$ be the number of orbits of $H$ on $A^i$ and,
for every $i = 1, \ldots r$, let $k_i$ be the number of relation symbols with arity $i$.
There is an integer $c(\mcA,H) > 0$, depending only on $\mcA$, $H$ 
and the vocabulary, such that
\[\big| \mbS_n(\mcA, H) \big| \ \sim \ 
c(\mcA,H) \binom{n}{p} \exp_2\Bigg( \sum_{i=1}^r k_i(n-p)^i \ + \ 
\sum_{j = 2}^r \sum_{i=1}^{j-1} \ k_j \binom{j}{i}q_i(n-p)^{j - i} \Bigg).\]
\end{prop}

\noindent
As will be explained below, Proposition~\ref{estimate of cardinality of S-n(A)}
is a consequence of Lemma~\ref{estimate of S-n(A-X,H)} which in turn
follows from Lemmas~\ref{estimate for Pi-1,...,Pi-r-1}--\ref{a negligible intersection}.

\begin{assump}\label{assumptions for section about asymptotic estimates}{\rm
For the rest of this section we assume the following, although the assumptions may be restated:
\begin{itemize}{\em
\item[] Suppose that $\mcA \in \mbS$  has no fixed point.
\item[] Let $H$ be a subgroup of $\Aut(\mcA)$ such that $H$ has no fixed point.}
\end{itemize}
Also let {\em
\begin{itemize}
\item[] $p = |A|$, 
\item[] for every $i = 1, \ldots, r-1$
let $q_i$ be the number of orbits of $H$ on $A^i$ and,
\item[] for every $i = 1, \ldots r$, let $k_i$ be the number of relation symbols with arity $i$.
\end{itemize}
}
}\end{assump}

\noindent
We consider the number of ways in which the relation symbols can be interpreted on $[n]$ so that
the resulting structure belongs to $\mbS_n(\mcA, H)$.
Let $c_\mcA$ be the number of structures in $\mbS_p$ that are isomorphic to $\mcA$.
First, it is clear that we can choose the set $X \subseteq [n]$ which is going
to be the support of the structure in $\binom{n}{p}$ ways, since we want that $|X| = p = |A|$.
Then we can choose interpretations of the relation symbols on $X$ in $c_\mcA$ ways so that
the resulting substructure with universe $X$, call it $\mcA_X$, is isomorphic to $\mcA$.
Now suppose that $X \subseteq [n]$ of cardinality $p$ and $\mcA_X \cong \mcA$ with universe $X$ are fixed.
Let 
\begin{equation}\label{definition of S-n-A-X-H}
\mbS_n(\mcA_X, H) \ = \ \big\{ \mcM \in \mbS_n(\mcA, H) : \mcM \uhrc \Spt^*(\mcM) = \mcA_X \big\}.
\end{equation}
Note that the condition $\mcM \uhrc \Spt^*(\mcM) = \mcA_X$ means that $\mcM \uhrc \Spt^*(\mcM)$
is {\em identical} with $\mcA_X$.
Also observe that if $X, X' \subseteq [n]$ and $X \neq X'$ then 
$\mbS_n(\mcA_X, H)$ and $\mbS_n(\mcA_{X'}, H)$ are disjoint.
Moverover, if both $\mcA'_X$ and $\mcA_X$ have universe $X$ and are isomorphic with $\mcA$,
but $\mcA'_X \neq \mcA_X$, then $\mbS_n(\mcA_X, H)$ and $\mbS_n(\mcA'_X, H)$ are disjoint.
Therefore Proposition~\ref{estimate of cardinality of S-n(A)} follows from the following:

\begin{lem}\label{estimate of S-n(A-X,H)}
Suppose that $X \subseteq [n]$ and $|X| = |A| = p$.
There is an integer $d(\mcA, H) > 0$, 
depending only on $\mcA$, $H$ and the vocabulary, such that
\[ \big| \mbS_n(\mcA_X, H) \big| \ \sim \ 
d(\mcA, H) \exp_2\Bigg( \sum_{i=1}^r k_i(n-p)^i \ + \ 
\sum_{j = 2}^r \sum_{i=1}^{j-1} \ k_j \binom{j}{i}q_i(n-p)^{j - i} \Bigg).\]
\end{lem}

\noindent
Lemma~\ref{estimate of S-n(A-X,H)} follows from
Lemmas~\ref{estimate for Pi-1,...,Pi-r-1}--\ref{a negligible intersection}, 
as we will show after proving them.
We begin with some preparatory work.
{\em Until Lemma~\ref{estimate of S-n(A-X,H)} has been proved we fix $X \subseteq [n]$
such that $|X| = |A| = p$ and $\mcA_X \cong \mcA$ with universe $X$.}
For every isomorphism $f : \mcA \to \mcA_X$, let
\[H_f = \{f\sigma f^{-1} : \sigma \in H \},\]
so $H_f$ is a subgroup of $\Aut(\mcA_X)$ and  $H_f \cong_P H$.

Suppose that $\mcM \in \mbS_n(\mcA_X, H)$.
By the definition of $\mbS_n(\mcA_X, H)$, $\mcM \uhrc \Spt^*(\mcM) = \mcA_X$ and there is an isomorphism
$f : \mcA \to \mcA_X$ such that $H_f$
is a subgroup of $\Aut(\mcM) \uhrc \Spt^*(\mcM)$.
For each $t = 1, \ldots, r-1$, the orbits of $H_f$ on $X^t$ forms a partition, denoted $\Pi_t$, of $X^t$.
Since $\Spt^*(\mcM) = X$ the following holds for $\mcM$:
\begin{itemize}
\item[(a)] Whenever $2 \leq j \leq r$, $R \in \{R_1, \ldots, R_\rho\}$ is a $j$-ary relation symbol,
$1 \leq i < j$, $(a_1, \ldots, a_i) \in X^i$ and $(a'_1, \ldots, a'_i) \in X^i$ belong to the
same part of $\Pi_i$ and 
$(a_{i+1}, \ldots, a_j) = (a'_{i+1}, \ldots, a'_j) \in \big([n] \setminus X\big)^{j - i}$,
then, for every $\pi \in Sym_j$, 
either both of 
\[(a_{\pi(1)}, \ldots, a_{\pi(j)}) \ \text{ and } \ 
(a'_{\pi(1)}, \ldots, a'_{\pi(j)}),\]
or none of them, belong to the interpretation of $R$.
\end{itemize}

\begin{defin}\label{definition of a structure respecting a pi-sequence}{\rm
If, for every $t = 1, \ldots, r-1$, $\Pi_t$ is a partition of $X^t$ such that~(a) holds (for $\mcM$),
then we say that $\mcM$ {\em respects} (the sequence of partitions) $\Pi_1, \ldots, \Pi_{r-1}$.
}\end{defin}

\noindent
In other words, if $\mcM \in \mbS_n(\mcA_X, H)$ then $\mcM \uhrc \Spt^*(\mcM) = \mcA_X$ and, 
for each $t = 1, \ldots, r-1$,
there is a partition $\Pi_t$ of $X^t$ such that 
$\mcM$ respects $(\Pi_1, \ldots, \Pi_{r-1})$ and for some isomorphism
$f : \mcA \to \mcA_X$, $\Pi_t$ is the set of orbits of $H_f$ on $X^t$ for $t = 1, \ldots, r-1$.
Conversely, if $\mcM \in \mbS_n$ is such that $\mcM \uhrc \Spt^*(\mcM) = \mcA_X$ and,
for each $t = 1, \ldots, r-1$, there is a partition $\Pi_t$ of $X^t$ such that
$\mcM$ respects $(\Pi_1, \ldots, \Pi_{r-1})$
and for some isomorphism $f : \mcA \to \mcA_X$,
$\Pi_t$ is the set of orbits of $H_f$ on $X^t$ for $t = 1, \ldots, r-1$,
then $H_f$ is a subgroup of $\mcM \uhrc \Spt^*(\mcM)$ and therefore $\mcM \in \mbS_n(\mcA_X, H)$.

\begin{defin}\label{definition of sequence of (A,H)-partitions}{\rm
A sequence $\Pi_1, \ldots, \Pi_{r-1}$ 
is called a {\em sequence of $(\mcA_X, H)$-partitions} if the following holds:
\begin{itemize}
\item[(b)] there is an isomorphism $f : \mcA \to \mcA_X$ such that, 
for each $t = 1, \ldots, r-1$, $\Pi_t$ is
the set of orbits of $H_f$ on $X^t$.
\end{itemize}
}\end{defin}

\noindent
For every sequence of $(\mcA_X, H)$-partitions $\Pi_1, \ldots, \Pi_{r-1}$
we define
\begin{align}\label{definition of S-n-A-X-Pi}
&\mbS_n(\mcA_X, \Pi_1, \ldots, \Pi_{r-1}) \ = \\ 
&\big\{ \mcM \in \mbS_n : \mcM \uhrc X = \mcA_X, \ \Spt^*(\mcM) = X 
\text{ and $\mcM$ respects } \Pi_1, \ldots, \Pi_{r-1} \}
\nonumber
\end{align}
and
\begin{align}\label{definition of T-n-A-X-Pi}
&\mbT_n(\mcA_X, \Pi_1, \ldots, \Pi_{r-1}) \ = \\
&\big\{ \mcM \in \mbS_n : \ \mcM \uhrc X = \mcA_X \text{ and $\mcM$ respects }
\Pi_1, \ldots, \Pi_{r-1} \}.
\nonumber
\end{align}

\noindent
If $\Pi_1, \ldots, \Pi_{r-1}$ is a sequence of $(\mcA_X, H)$-partitions and
$\mcM \in \mbT_n(\mcA_X, \Pi_1, \ldots, \Pi_{r-1})$ then 
there is an isomorphism $f : \mcA \to \mcA_X$ such that, for $i = 1, \ldots, r-1$,
$\Pi_i$ is the set of orbits of $H_f$ on $X^i$ and,
by Assumption~\ref{assumptions for section about asymptotic estimates},
every orbit of $H_f$ on $X$ has at least two members; hence $X \subseteq \Spt^*(\mcM)$.
Consequently, 
\[
\mbS_n(\mcA_X, \Pi_1, \ldots, \Pi_{r-1}) \ \subseteq \ \mbT_n(\mcA_X, \Pi_1, \ldots, \Pi_{r-1})
\]
From the argument before the definition of $\mbS_n(\mcA_X, \Pi_1, \ldots, \Pi_{r-1})$
it follows that 
\begin{equation}\label{union of all S-n(A-X, Pi)}
\mbS_n(\mcA_X, H) \ = \ \bigcup_{\Pi_1, \ldots, \Pi_{r-1}} \mbS_n(\mcA_X, \Pi_1, \ldots, \Pi_{r-1}),
\end{equation}
where the union ranges over all sequences $\Pi_1, \ldots, \Pi_{r-1}$
of $(\mcA_X, H)$-partitions.
The next step in the proof of Lemma~\ref{estimate of S-n(A-X,H)}
is to estimate $\big| \mbS_n(\mcA_X, \Pi_1, \ldots, \Pi_{r-1}) \big|$.
Then we deal with the slightly problematic issue that even if 
$\Pi_1, \ldots, \Pi_{r-1}$ and $\Pi'_1, \ldots, \Pi'_{r-1}$ are different
sequences of $(\mcA_X, H)$-partitions it may be the case that
$\mbS(\mcA_X, \Pi_1, \ldots, \Pi_{r-1})$ and $\mbS(\mcA_X, \Pi'_1, \ldots, \Pi'_{r-1})$
have nonempty intersection. However, as we will show, their intersection
will always be negligibly small, which implies that we can add the asymptotic estimates of the
cardinalities of all $\mbS_n(\mcA_X, \Pi_1, \ldots, \Pi_{r-1})$ to get an asymptotic estimate of the
cardinality of $\mbS_n(\mcA_X, H)$.
Recall that for $i = 1, \ldots, r$, $k_i$ is the number of $i$-ary relation symbols.
Also, $p = |A| = |X|$ and, for $i = 1, \ldots, r-1$, $q_i$ is the number of orbits
of $H$ on $A^i$. 

\begin{lem}\label{estimate for Pi-1,...,Pi-r-1}
If $\Pi_1, \ldots, \Pi_{r-1}$ is a sequence of $(\mcA_X, H)$-partitions, then
\[
\big| \mbT_n(\mcA_X, \Pi_1, \ldots, \Pi_{r-1}) \big| \ = \ 
\exp_2\Bigg( \sum_{i=1}^r k_i(n-p)^i \ + \ 
\sum_{j = 2}^r \sum_{i=1}^{j-1} \ k_j \binom{j}{i}q_i(n-p)^{j - i} \Bigg).
\]
Moreover, there is $\varepsilon : \mbbN \to \mbbR$, depending only on $\mcA$, $H$ and the vocabulary,
such that $\lim_{n\to\infty} \varepsilon(n) = 0$ and for all large enough $n$ 
the proportion of $\mcM \in \mbT_n(\mcA_X, \Pi_1, \ldots, \Pi_{r-1})$ 
such that $\mcM \notin \mcS_n(\mcA_X, \Pi_1, \ldots, \Pi_{r-1})$ is at most $\varepsilon(n)$.
\end{lem}

\noindent
{\bf Proof.}
Suppose that $\Pi_1, \ldots, \Pi_{r-1}$ is a sequence of $(\mcA_X, H)$-partitions,
so there is an isomorphism $f : \mcA \to \mcA_X$ such that, 
for each $t = 1, \ldots, r-1$, $\Pi_t$ is
the set of orbits of $H_f$ on $X^t$.
Since $H_f \cong_P H$ it follows that 
$\Pi_t$ partitions $X^t$ into $q_t$ parts, for every $t = 1, \ldots, r-1$.
Let
\[
\gamma(n) \ = \ 
\exp_2\Bigg( \sum_{i=1}^r k_i(n-p)^i \ + \ 
\sum_{j = 2}^r \sum_{i=1}^{j-1} \ k_j \binom{j}{i}q_i(n-p)^{j - i} \Bigg).
\]
First we will prove that $\big| \mbT_n(\mcA_X, \Pi_1, \ldots, \Pi_{r-1}) \big| = \gamma(n)$.
As observed before Lemma~\ref{estimate for Pi-1,...,Pi-r-1},
\[
\mbS_n(\mcA_X, \Pi_1, \ldots, \Pi_{r-1}) \ \subseteq \ \mbT_n(\mcA_X, \Pi_1, \ldots, \Pi_{r-1})
\]
and $X \subseteq \Spt^*(\mcM)$ for every $\mcM \in \mbT_n(\mcA_X, \Pi_1, \ldots, \Pi_{r-1})$.
Then we show that the proportion of $\mcM \in \mbT_n(\mcA_X, \Pi_1, \ldots, \Pi_{r-1})$ 
such that $X$ is a proper subset
of $\Spt^*(\mcM)$ approaches 0 as $n \to \infty$. 
Moreover, we will get a bound $\varepsilon(n)$ as in the lemma. 
For the rest of the proof of this lemma we use the abbreviation
\[\mbT_n \ = \ \mbT_n(\mcA_X, \Pi_1, \ldots, \Pi_{r-1}).\]
To determine $\big| \mbT_n \big|$ we consider the number of ways 
in which the relation symbols can be interpreted on $[n]$
so that the resulting structure $\mcM$ has the properties that
$\mcM \uhrc X = \mcA_X$ and $\mcM$ respects $\Pi_1, \ldots, \Pi_{r-1}$, that is,~(a)
holds for $\mcM$.
Since the substructure on $X$ must be $\mcA_X$, there is only one choice for the
interpretations on tuples all of which coordinates belong to $X$.

Now we consider in how many ways the relation symbols can be interpreted on
tuples that intersect both $X$ and $[n] \setminus X$ so that 
resulting structure respects $\Pi_1, \ldots, \Pi_{r-1}$; 
so in this stage we only consider
relation symbols of arity at least 2.
Let $R \in \{R_1, \ldots, R_\rho\}$ 
be a relation symbol of arity $j \geq 2$ and let $1 \leq i \leq j - 1$.
We consider the number of ways in which $R$ can be interpreted on $j$-tuples $\bar{a} \in [n]^j$
with exactly $i$ coordinates of $\bar{a}$ from $X$ in such a way that the resulting
structure respects $\Pi_1, \ldots, \Pi_{r-1}$.

Suppose that
\[a_1, \ldots, a_i, a'_1, \ldots, a'_i \in X \ \text{ and } \ b_{i+1}, \ldots, b_j \in [n] \setminus X\]
and that the $i$-tuples $(a_1, \ldots, a_i)$ and $(a'_1, \ldots, a'_i)$ belong to the same part of $\Pi_i$.
Since we want~(a) to be satisfied we have the choice of
letting both $j$-tuples
\[(a_1, \ldots, a_i, b_{i+1}, \ldots, b_j) \ \text{ and } \ (a'_1, \ldots, a'_i , b_{i+1}, \ldots, b_j),\] 
or none of them, belong to the interpretation of $R$ (and this independently of other choices).
We considered the case when $a_1, \ldots, a_i$ and $a'_1, \ldots, a'_i$ 
occured in the first $i$ positions of the respective
$j$-tuple, but the same is clearly true if  $a_1, \ldots, a_i$ and $a'_1, \ldots, a'_i$
take other positions in the respective $j$-tuples, but still so that $a_l$ 
preceeds $a_{l'}$ if $l < l'$ and $a_l$ takes position $t$ if and only if $a'_l$ takes position $t$.\footnote{
We consider only the given order of $a_1, \ldots, a_i$ and $a'_1, \ldots, a'_i$ because, in general,
an $i$-tuple obtained by reordering $a_1, \ldots, a_i$ need not belong to the same part of $\Pi_i$
as $(a_1, \ldots, a_i)$.}
There are $\binom{j}{i}$ ways in which $i$ positions in an $j$-tuple can be chosen.
Therefore the number of ways to choose the interpretation of $R$ on $j$-tuples with
exactly $i$ coordinates in $X$ in such a way that~(a) is satisfied is 
\[\exp_2\bigg( \binom{j}{i}q_i(n-p)^{j - i} \bigg),\]
where we recall that $q_i$ is the number of parts of the partition $\Pi_i$ of $X^i$.\footnote{
If we assume that $R$ is always interpreted as an irreflexive and symmetric relation,
then the corresponding number is $\exp_2(q'_i\binom{n-p}{j-i})$ where $q'_i$ is the 
number of orbits of the action of $H$ on $\{B \subseteq A : |B| = i\}$ given by
$h(\{b_1, \ldots, b_i\}) = \{h(b_1), \ldots, h(b_i)\}$ for every $h \in H$ and 
$i$-subset $\{b_1, \ldots, b_i\}$ of $A$.}
If $i' \neq i$ and $1 \leq i' \leq j - 1$
then the corresponding number of choices for $j$-tuples with exactly $i'$
coordinates in $X$ is independent from the previously made choices.
Therefore the number of ways in which $R$ can be interpreted on tuples that intersect both
$X$ and $[n] \setminus X$ is
\[\exp_2\bigg( \sum_{i=1}^{j-1} \binom{j}{i}q_i(n-p)^{j - i} \bigg).\]
The same argument can be carried out for every relation symbol $R$ of arity at least 2.
The number of choices for each such $R$ is independent of previously made choices.
Therefore the number of ways in which all relation symbols with arity at least 2 can
be interpreted on tuples that intersect both $X$ and $[n] \setminus X$ in such a way that~(a)
is satisfied is 
\begin{equation}\label{number of relationships between X and n-X in general case}
\exp_2\bigg( \sum_{j = 2}^r \sum_{i=1}^{j-1} \ k_j \binom{j}{i}q_i(n-p)^{j - i} \bigg).
\end{equation}

Finally we consider interpretations on tuples $\bar{a}$ such that none of the coordinates
of $\bar{a}$ belongs to $X$.
If $R$ has arity $i$, then there are $2^{(n-p)^i}$ ways in which to
interpret $R$ on tuples $\bar{a} \in ([n] \setminus X)^i$, independently of other choices.
As there are $k_i$ relation symbols of arity $i$, the number of ways to interpret all relation
symbols on $[n] \setminus X$ is 
\begin{equation}\label{number of interpretations on X and on n-X in general case}
\exp_2\bigg( \sum_{i=1}^r k_i(n-p)^i \bigg),
\end{equation}

Suppose that a structure $\mcM$ has been constructed by making the choices described above.
Then, by construction, $\mcM \uhrc X = \mcA_X$ and $\mcM$ respects $\Pi_1, \ldots, \Pi_{r-1}$.
By assumption, $H$ has no fixed point which implies that every part of the
partition $\Pi_1$ of $X$ has at least two members. 
Since $\mcM$ respects $\Pi_1, \ldots, \Pi_{r-1}$ and $\Pi_1, \ldots, \Pi_{r-1}$ 
is a sequence of $(\mcA_X, H)$-partitions
it follows that $X \subseteq \Spt^*(\mcM)$.
It is also clear that every member of $\mbT_n$ can be obtained in exactly one way by 
making choices as described by the construction.
Hence, by multiplying~(\ref{number of interpretations on X and on n-X in general case})
and~(\ref{number of relationships between X and n-X in general case}),
we see that $\big| \mbT_n \big| = \gamma(n)$. 

It remains to prove that for all large enough $n$,
\begin{equation}\label{structures in T-n with bigger support than X}
\frac{\big| \big\{ \mcM \in \mbT_n : \Spt^*(\mcM) \neq X \big\} \big|}{\big| \mbT_n \big|} 
\ \leq \ \varepsilon(n),
\end{equation}
where $\lim_{n\to\infty}\varepsilon(n) = 0$ and
$\varepsilon$ depends only on $\mcA, H$ and the vocabulary.
After defining $\mbT_n = \mbT_n(\mcA_X, \Pi_1, \ldots, \Pi_{r-1})$, 
see~(\ref{definition of T-n-A-X-Pi}), 
we observed that if $\mcM \in \mbT_n$ then $X \subseteq \Spt^*(\mcM)$.
Since $\Pi_1, \ldots, \Pi_{r-1}$
is a sequence of $(\mcA_X, H)$-partitions, there is an isomorphism
$f : \mcA \to \mcA_X$ such that, for each $t = 1, \ldots, r-1$,
$\Pi_t$ is the set of orbits of $H_f = \{f\sigma f^{-1} : \sigma \in H\}$ on $X^t$.
Let $H_f = \{h_1, \ldots, h_s\}$ and extend every $h_i \in H_f$ to $h'_i \in Sym_n$
by $h'_i(x) = h_i(x)$ if $x \in X$ and $h'_i(x) = x$ if $x \in [n] \setminus X$.
Then $\Spt(h'_1, \ldots, h'_s) = X$ and hence $\spt(h'_1, \ldots, h'_s) = |X| = |A| = p$.

If $\mcM \in \mbS_n$, then $\mcM$ belongs to $\mbS_n(h'_1, \ldots, h'_s)$ if and only if
the following condition holds:
for every $t = 1, \ldots, r$ and every $t$-ary relation symbol $R$, if
$\bar{a}$ and $\bar{b}$ are two $t$-tuples from the same orbit
of $\langle h'_1, \ldots, h'_s \rangle$ on $[n]^t$ (which here denotes the set of ordered $t$-tuples 
of elements from $[n]$), 
then either $\mcM \models R(\bar{a}) \wedge R(\bar{b})$ or $\mcM \models \neg R(\bar{a}) \wedge \neg R(\bar{b})$.
As $X$ is a union of orbits of $\langle h'_1, \ldots, h'_s \rangle$ it follows that if we define
\[
\mbS_n(h'_1, \ldots, h'_s, \mcA_X) \ = \ 
\big\{ \mcM \in \mbS_n(h'_1, \ldots, h'_s) : \mcM \uhrc X = \mcA_X \},
\]
then there is a constant $0 < c \leq 1$, depending only on $\mcA$, $H$ and the vocabulary, 
such that 
\begin{equation}\label{relationship between S-n and S-n with A-x}
\big| \mbS_n(h'_1, \ldots, h'_s, \mcA_X) \big| \ = \ c \big| \mbS_n(h'_1, \ldots, h'_s) \big|.
\end{equation}
From the definition of $h'_1, \ldots, h'_s$ it follows that
\begin{equation}\label{S-n with A-X is inluded in T-X}
\mbS_n(h'_1, \ldots, h'_s, \mcA_X) \ \subseteq \ \mbT_n.
\end{equation}
By~(\ref{relationship between S-n and S-n with A-x}), (\ref{S-n with A-X is inluded in T-X}) 
and Propositions~\ref{automorphisms with large support are unusual}
and~\ref{support of structure is bounded in terms of support of automorphisms},
there are $\lambda, p_0 > 0$, depending only on $\mcA$, $H$ and the vocabulary, such that
for all sufficiently large $n$,
\begin{align*}
\frac{\big| \mbS_n(\spt^* > p_0) \big|}{\big| \mbT_n \big|} \ \leq \ 
\frac{\big| \mbS_n(\spt^* > p_0) \big|}{c \big| \mbS_n(h'_1, \ldots, h'_s) \big|} \ \leq \ 
2^{-\lambda n^{r-1}}.
\end{align*}
Hence, for all large enough $n$, the proportion of $\mcM \in \mbT_n$ such that $\spt^*(\mcM) \leq p_0$
is at least $1 - 2^{-\lambda n^{r-1}}$.

Fix any $a \in [n] \setminus X$ and $a' \in [n]$ such that $a \neq a'$.
From the definition of $\mbT_n$ it is clear that for every sequence of distinct $(r-1)$-tuples 
$\bar{b}_1, \ldots, \bar{b}_\kappa \in \big([n] \setminus (X \cup \{a, a'\})\big)^{r-1}$,
the proportion of $\mcM \in \mbT_n$ that satisfies the following is $2^{-\kappa}$:
\begin{equation}\label{tuples outside support satisfying the same relations}
\text{for every } i = 1, \ldots, \kappa, \ 
\mcM \models R(a,\bar{b}_i) \ \Longleftrightarrow \mcM \models R(a',\bar{b}_i).
\end{equation}
Observe that if $\mcM \in \mbT_n$, $\spt^*(\mcM) \leq p_0$ and $g(a) = a'$ for some
$g \in \Aut(\mcM)$,
then there is a sequence of distinct $(r-1)$-tuples
$\bar{b}_1, \ldots, \bar{b}_\kappa \in \big([n] \setminus (X \cup \{a, a'\})\big)^{r-1}$
such that $\kappa = 2^{(n - p_0 - 2)^{r-1}}$ 
and~(\ref{tuples outside support satisfying the same relations}) is satisfied.
Hence the proportion of $\mcM \in \mbT_n$ such that $\spt^*(\mcM) \leq p_0$
and $g(a) = a'$ for some $g \in \Aut(\mcM)$ is at most $2^{-(n - p_0 - 2)^{r-1}}$.
As the proportion of $\mcM \in \mbT_n$ such that $\spt^*(\mcM) \leq p_0$
is at least $1 - 2^{-\lambda n^{r-1}}$, it follows that that the proportion of
$\mcM \in \mbT_n$ with an automorphism $g$ such that $g(a) = a'$ is at most
$2^{-(n - p_0 - 2)^{r-1}} + 2^{-\lambda n^{r-1}}$.
It follows that the proportion of $\mcM \in \mbT_n$ which have 
distinct elements $a \in [n] \setminus X$ and $a' \in [n]$ and an automorphism
$g$ such that $g(a) = a'$ is at most $n^2 \Big( 2^{-(n - p_0 - 2)^{r-1}} + 2^{-\lambda n^{r-1}} \Big)$.
This immediately implies~(\ref{structures in T-n with bigger support than X}),
so the proof of Lemma~\ref{estimate for Pi-1,...,Pi-r-1} is finished.
\hfill $\square$

\begin{rem}\label{remark about irreflexive and symmetric relations for the asymptotic formula}{\rm
If we assume that all relation symbols are always interpreted as irreflexive and symmetric
relations then we get
\[
\big| \mbT_n(\mcA_X, \Pi_1, \ldots, \Pi_{r-1}) \big| \ = \ 
\exp_2\Bigg( \sum_{i=1}^r k_i \binom{n-p}{i} \ + \ 
\sum_{j = 2}^r \sum_{i=1}^{j-1} \ k_j q'_i \binom{n-p}{j - i} \Bigg),
\]
where $q'_i$ is the number of orbits of the action of $H$ on $\{B \subseteq A : |B| = i\}$ given by
$h(\{b_1, \ldots, b_i\}) = \{h(b_1), \ldots, h(b_i)\}$ for every $h \in H$ and 
$i$-subset $\{b_1, \ldots, b_i\}$ of $A$.
Under the same assumptions we still have 
$\big| \mbS_n(\mcA_X, \Pi_1, \ldots, \Pi_{r-1}) \big| \sim 
\big| \mbT_n(\mcA_X, \Pi_1, \ldots, \Pi_{r-1}) \big|$ by the same argument as above
(and a modification of Proposition~\ref{automorphisms with large support are unusual}).
}\end{rem}

\begin{lem}\label{an unusual property}
Suppose that $\Pi_1, \ldots, \Pi_{r-1}$ is a sequence of $(\mcA_X, H)$-partitions.
For each $1 \leq i < r$, 
the proportion of $\mcM \in \mbS_n(\mcA_X,\Pi_1, \ldots, \Pi_{r-1})$ with the following property
is at most $\varepsilon(n)$ where $\varepsilon(n) \to 0$ as $n \to 0$ and the function $\varepsilon$ 
depends only on $\mcA$, $H$
and the vocabulary:
\begin{itemize}
\item[$(\dagger)$] There are an $r$-ary relation symbol $R$,
different parts $P, P' \in \Pi_i$, $\bar{a} = (a_1, \ldots, a_i) \in P$
and $\bar{a}' = (a'_1, \ldots, a'_i) \in P'$ such that
for every $\bar{b} = (b_{i+1}, \ldots, b_r) \in ([n] \setminus X)^{r-i}$,
\[\mcM \models R(\bar{a}, \bar{b}) \ \Longleftrightarrow \ \mcM \models R(\bar{a}', \bar{b}).\] 
\end{itemize}
\end{lem}

\noindent
{\bf Proof.}
By Lemma~\ref{estimate for Pi-1,...,Pi-r-1} it suffices to prove that 
the proportion of $\mcM \in \mbT_n(\mcA_X,\Pi_1, \ldots, \Pi_{r-1})$ with property~($\dagger$)
is at most $\varepsilon(n)$ where $\varepsilon(n) \to 0$ as $n \to 0$ and $\varepsilon$ 
depends only on $\mcA$, $H$ and the vocabulary.
Suppose that $\mcM \in \mbT_n(\mcA_X,\Pi_1, \ldots, \Pi_t)$ and~$(\dagger)$ holds, so
there are different parts $P, P' \in \Pi_i$,
$\bar{a} = (a_1, \ldots, a_i) \in P$ and
$\bar{a}' = (a'_1, \ldots, a'_i) \in P'$ such that for 
{\em every} $\bar{b} = (b_{i+1}, \ldots, b_r) \in ([n] \setminus X)^{r-i}$,
\[\mcM \models R(\bar{a}, \bar{b}) \ \Longleftrightarrow \ \mcM \models R(\bar{a}', \bar{b}).\] 
Fix these tuples $\bar{a}$ and $\bar{a}'$.
The number of ways in which we can interpret $R$, in such $\mcM$,
on tuples of the form $\bar{a}\bar{b}$ and $\bar{a}'\bar{b}$ where $\bar{b} \in ([n] \setminus X)^{r-i}$
is $\displaystyle 2^{(n-p)^{r-i}}$, independently of how $R$ is interpreted on other tuples and independently
of how other relation symbols are interpreted.

On the other hand, for $\mcM \in \mbT_n(\mcA_X,\Pi_1, \ldots, \Pi_{r-1})$ without
property~$(\dagger)$, the number of ways in which
$R$ can be interpreted on tuples of the form $\bar{a}\bar{b}$ and $\bar{a}'\bar{b}$ where 
$\bar{b} \in ([n] \setminus X)^{r-i}$ is 
$\displaystyle 4^{(n-p)^{r-i}}$, independently of how $R$ is interpreted on other tuples and independently
of how other relation symbols are interpreted.
Therefore the proportion of $\mcM \in \mbT_n(\mcA_X,\Pi_1, \ldots, \Pi_{r-1})$ with property~$(\dagger)$
is at most $\displaystyle 2^{(n-p)^{r-i}} \Big/ 4^{(n-p)^{r-i}} \leq 2^{-(n-p)}$.
\hfill $\square$

\begin{lem}\label{a negligible intersection}
If  $\Pi_1, \ldots, \Pi_{r-1}$ and 
$\Pi'_1, \ldots, \Pi'_{r-1}$ are two different sequences 
of $(\mcA_X, H)$-partitions,
then 
\[
\frac{\big| \mbS_n(\mcA_X, \Pi_1, \ldots, \Pi_{r-1}) \ \cap \ \mbS_n(\mcA_X, \Pi'_1, \ldots, \Pi'_{r-1}) \big|}
{\big| \mbS_n(\mcA_X, \Pi_1, \ldots, \Pi_{r-1}) \ \cup \ \mbS_n(\mcA_X, \Pi'_1, \ldots, \Pi'_{r-1}) \big|}
\ \leq \ \varepsilon(n).\]
where $\varepsilon(n) \to 0$ as $n \to 0$ and the function $\varepsilon$ 
depends only on $\mcA$, $H$ and the vocabulary.
\end{lem}

\noindent
{\bf Proof.}
Suppose that $\Pi_1, \ldots, \Pi_{r-1}$ and 
$\Pi'_1, \ldots, \Pi'_{r-1}$ are different sequences of $(\mcA_X, H)$-partitions
and that 
\[\mcM \ \in \  \mbS_n(\mcA_X, \Pi_1, \ldots, \Pi_{r-1}) \ \cap \ \mbS_n(\mcA_X, \Pi'_1, \ldots, \Pi'_{r-1}).\]
Then for some $1 \leq i < r$, there are $\bar{a}, \bar{a}' \in X^i$ such that 
$\bar{a}$ and $\bar{a}'$ are in the same part of the partition $\Pi'_i$ but in different
parts of the partition $\Pi_i$, or vice versa.
In the first case, $\mcM$ has property~$(\dagger)$ 
from Lemma~\ref{an unusual property} (for every $r$-ary relation symbol $R$)
when seen as a member of $\mbS_n(X, \Pi_1, \ldots, \Pi_{r-1})$.
In the second case, $\mcM$ has property~$(\dagger)$ 
when seen as a member of $\mbS_n(X, \Pi'_1, \ldots, \Pi'_{r-1})$.
Therefore, using Lemma~\ref{an unusual property}, 
the quotient of the lemma is at most $2 \varepsilon(n)$ where $\varepsilon(n) \to 0$ 
as $n \to 0$ and the function $\varepsilon$ 
depends only on $\mcA$, $H$ and the vocabulary..
\hfill $\square$
\\

\noindent
{\bf Proof of Lemma~\ref{estimate of S-n(A-X,H)}.}
Let
\[\gamma(n) \ = \ 
\exp_2\Bigg( \sum_{i=1}^r k_i(n-p)^i \ + \ 
\sum_{j = 2}^r \sum_{i=1}^{j-1} \ k_j \binom{j}{i}q_i(n-p)^{j - i} \Bigg)\]
and let $d(\mcA, H)$ be the number of different sequences $\Pi_1, \ldots, \Pi_{r-1}$
of $(\mcA_X, H)$-partitions.
Hence, $d(\mcA, H)$ is finite and depends only on $\mcA$, $H$ and the vocabulary.
We prove that $\big| \mbS_n(\mcA_X, H) \big| \sim d(\mcA, H)\gamma(n)$.
From~(\ref{union of all S-n(A-X, Pi)}) it follows that
\begin{equation}\label{upper bound of S-n-X}
\big| \mbS_n(\mcA_X, H) \big| \ \leq \ d(\mcA, H) \gamma(n).
\end{equation}
Let $\mbU_n$ be the union of all intersections
\[ \mbS_n(\mcA_X, \Pi_1, \ldots, \Pi_{r-1}) \ \cap \ \mbS_n(\mcA_X, \Pi'_1, \ldots, \Pi'_{r-1}) \]
where $\Pi_1, \ldots, \Pi_{r-1}$ and $\Pi'_1, \ldots, \Pi'_{r-1}$ range over all unordered
pairs of different sequences of $(\mcA_X, H)$-partitions.
If the sums below ranges over such unordered pairs, then,
by Lemma~\ref{a negligible intersection}, we have 
\begin{align*}
\big| \mbU_n \big| \ 
&\leq \ 
\sum \Big| \mbS_n(\mcA_X, \Pi_1, \ldots, \Pi_{r-1}) \ \cap \ \mbS_n(\mcA_X, \Pi'_1, \ldots, \Pi'_{r-1}) \Big|  \\
&\leq \ \varepsilon(n) \sum 
\Big( \big| \mbS_n(\mcA_X, \Pi_1, \ldots, \Pi_{r-1}) \big| \ + \ 
\big| \mbS_n(\mcA_X, \Pi'_1, \ldots, \Pi'_{r-1}) \big| \Big) \\
&\sim \ \varepsilon(n) \cdot \binom{d(\mcA, H)}{2} \cdot 2\gamma(n),
\end{align*}
where $\varepsilon(n) \to 0$ as $n \to \infty$
By Lemma~\ref{estimate for Pi-1,...,Pi-r-1}, 
$\mbS_n(\mcA_X, \Pi_1, \ldots, \Pi_{r-1}) \sim \gamma(n)$ for
every sequence $\Pi_1, \ldots, \Pi_{r-1}$ of $(\mcA_X, H)$-partitions. 
It follows that, for every such sequence,
\[
\big| \mbS_n(\mcA_X, \Pi_1, \ldots, \Pi_{r-1}) \big| \ - \ \big| \mbU_n \big| \ \sim \ \gamma(n).
\]
Since $\mbS_n(\mcA_X, \Pi_1, \ldots, \Pi_{r-1}) \ \setminus \ \mbU_n$
and $\mbS_n(\mcA_X, \Pi'_1, \ldots, \Pi'_{r-1}) \ \setminus \ \mbU_n$ are disjoint
if $\Pi_1, \ldots, \Pi_{r-1}$ and $\Pi'_1, \ldots, \Pi'_{r-1}$ are different sequences, it follows that
\begin{align*}
\big| \mbS_n(\mcA_X, H) \big| \ 
&\geq \ 
\sum \big| \mbS_n(\mcA_X, \Pi_1, \ldots, \Pi_{r-1}) \ \setminus \ \mbU_n \big| \\
&\geq \ 
\sum \Big( \big| \mbS_n(\mcA_X, \Pi_1, \ldots, \Pi_{r-1}) \big| \ - \ \big| \mbU_n \big| \Big)
\ \sim \ d(\mcA, H)\gamma(n),
\end{align*}
where the sums range over all sequences $\Pi_1, \ldots, \Pi_{r-1}$ 
of $(\mcA_X, H)$-partitions.
This together with~(\ref{upper bound of S-n-X}) implies that 
$\big| \mbS_n(\mcA_X, H) \big| \sim d(\mcA, H)\gamma(n)$,
so Lemma~\ref{estimate of S-n(A-X,H)} is proved
\hfill $\square$
\\

\noindent
As explained in the paragraph after the statement of 
Proposition~\ref{estimate of cardinality of S-n(A)}, it follows from
Lemma~\ref{estimate of S-n(A-X,H)}, so now we have also proved
Proposition~\ref{estimate of cardinality of S-n(A)}.
We can now derive two corollaries of this proposition.
These corollaries, as well as the proposition itself will be used in the next section.
It will be convenient to use the following notation:

\begin{defin}\label{definition of p(H) etc}{\rm
Suppose that $H$ is a group of permutations of the set $\Omega$.
Then $p(H) = |\Omega|$, 
$q(H)$ is the number of orbits of $H$ on $\Omega$ and
$s(H)$ is the number of orbits of $H$ on $\Omega^2$.
}\end{defin}

\begin{cor}\label{estimate of the number of structures for given p,q when r=2}
Suppose that $r = 2$, that $\mcA \in \mbS$  has no fixed point
and let $H$ be a subgroup of $\Aut(\mcA)$ without any fixed point.
Let $p = p(H) = |A|$, let $q = q(H)$ and for $i = 1,2$ let $k_i$ 
be the number of relation symbols of arity $i$.
Then there is an integer $c(\mcA,H) > 0$, depending only on $\mcA$, $H$ and the vocabulary, such that
\[
\big| \mbS_n(\mcA, H) \big| \ \sim \ 
c(\mcA,H) \binom{n}{p} \exp_2\Big( k_2 n^2 \ - \ 2k_2 (p-q)n \ + \ k_1 n \ + \ k_2p^2 \ - \ k_1 p \Big).
\]
\end{cor}

\noindent
{\bf Proof.}
By Proposition~\ref{estimate of cardinality of S-n(A)}
with $r = 2$ and $q = q_1$, there is an integer $c(\mcA,H) > 0$,
depending only on $\mcA$, $H$ and the vocabulary, such that
\begin{align*}
\big| \mbS_n(\mcA, H) \big| \ &\sim \ 
c(\mcA,H) \binom{n}{p} \exp_2\Bigg( \sum_{i=1}^2 k_i(n-p)^i \ + \ 2k_2 q(n-p) \Bigg) \\
&= \ 
c(\mcA,H) \binom{n}{p} \exp_2\Big( k_2 n^2 \ - \ 2k_2 (p-q)n \ + \ k_1 n \ + \ k_2p^2 \ - \ k_1 p \Big).
\end{align*}
\hfill $\square$

\begin{cor}\label{estimate of the number of structures for given p,q,s when r>2}
Suppose that $r > 2$, that $\mcA \in \mbS$  has no fixed point
and that $H$ is a subgroup of $\Aut(\mcA)$ without any fixed point.
Let $p = p(H) = |A|$, let $q = q(H)$ and $s = s(H)$.
Moreover, let $k$ be the number of $r$-ary relation symbols, let $l$ be the
number of $(r-1)$-ary relation symbols, let $m$ be the number of $(r-2)$-ary relation symbols 
and define 
\[
\beta(x,y,z) \ = \ k\binom{r}{2}x^2 \ - \ kr(r-1)xy \ - \ l(r-1)x \ + \ l(r-1)y \ + \ k\binom{r}{2}z.
\]
Then there is an integer $c(\mcA,H)$, depending only on $\mcA$, $H$ and the vocabulary,
such that 
\begin{align*}
&\big|\mbS_n(\mcA, H)\big| \ \sim \ c(\mcA,H) \binom{n}{p} \cdot \\
&\exp_2\bigg( k n^r \ - \ 
\big(kr(p-q) \ - \ l\big)n^{r-1} \ + \ 
\big(\beta(p,q,s) \ + \ m\big)n^{r-2} \ + \
\mcO\big(n^{r-3}\big) \bigg).
\end{align*}
\end{cor}

\noindent
{\bf Proof.}
For every $i = 1, \ldots r-1$, let $q_i$ be the number of 
orbits of $H$ on $A^i$. 
For every $j = 1, \ldots, r$, let $k_j$ be the number of relation symbols
of arity $j$. 
So  we have 
$q_1 = q$, $q_2 = s$, $k_r = k$, $k_{r-1} = l$ and $k_{r-2} = m$.
By Proposition~\ref{estimate of cardinality of S-n(A)}, there is an integer $c(\mcA,H) > 0$,
depending only on $\mcA$, $H$ and the vocabulary, such that 
\[
\big| \mbS_n(\mcA, H) \big| \ \sim \ 
c(\mcA,H) \binom{n}{p} \exp_2\big(\lambda(n)\big),
\]
where
\begin{align*}
\lambda(n) \ &= \ 
\sum_{i=1}^r k_i(n-p)^i \ + \ 
\sum_{j = 2}^r \sum_{i=1}^{j-1} \ k_j \binom{j}{i}q_i(n-p)^{j - i} \\
&= \ 
\Bigg( \sum_{i=1}^{r-3} k_i (n-p)^i \Bigg) \ + \ 
m(n-p)^{r-2} \ + \ l(n-p)^{r-1} \ + \ k(n-p)^r \\ 
&+ \
\Bigg( \sum_{j = 2}^{r-2} \sum_{i=1}^{j - 1} k_j \binom{j}{i} q_i (n-p)^{j - i} \Bigg)
\ + \ \Bigg( \sum_{i=1}^{r-2} l \binom{r-1}{i} q_i (n-p)^{r-1-i} \Bigg) \\
&+ \ 
\sum_{i=1}^{r-1} k \binom{r}{i} (n-p)^{r-i} \\
&= \ m(n-p)^{r-2} \ + \ l(n-p)^{r-1} \ + \ k(n-p)^r \\
&+ \ l(r-1)q(n-p)^{r-2} \ + \ krq(n-p)^{r-1} \ + \ k\binom{r}{2}s(n-p)^{r-2} \ + \ \mcO\big(n^{r-3}\big) \\
&= \ k n^r \ - \ 
\big(kr(p-q) \ - \ l\big)n^{r-1} \ + \ 
\big(\beta(p,q,s) \ + \ m\big)n^{r-2} \ + \
\mcO\big(n^{r-3}\big). 
\end{align*}
\hfill $\square$

\section{Comparing different groups}\label{Comparing different groups}

\noindent
In this section we use the analysis from Section~\ref{asymptotic estimates}
to prove Theorem~\ref{main theorem about comparissons between different groups},
which collects the statements of Propositions~\ref{comparisson of two subgroups},
~\ref{proportion of structures with subgroup G that have automorphism group G}
and~\ref{comparisson of two groups as automorphism groups}.
The main technical result of the section is Proposition~\ref{general quotient proposition}
which helps to break down more complex problems to problems about 
quotients of the form $\mbS_n(\mcA, H) \big/ \mbS_n(\mcA', H')$, where the meaning of
$\mbS_n(\mcA, H)$ was given by Definition~\ref{definition of S-n(A)}.
Also recall Definition~\ref{definition of p(H) etc} of $p(H)$, $q(H)$ and $s(H)$ for a permutation group $H$.
As usual, {\em $r$ denotes the maximal arity and in this section $k$ denotes the number of
$r$-ary relation symbols and $l$ denotes the number of $(r-1)$-ary relation symbols.} 
Parts~(ii) and~(iii) of the next result will not be used in this article, but
they are used in~\cite{Kop13c}.

\begin{prop}\label{comparisson of structures with different support-substructures}
Suppose that $\mcA, \mcA' \in \mbS$ are such that neither $\Aut(\mcA)$ nor $\Aut(\mcA')$ has a
fixed point. Moreover, suppose that $H$ is a subgroup of $\Aut(\mcA)$ without fixed any point and 
that $H'$ is a subgroup of $\Aut(\mcA')$ without any fixed point. 
Let $p = p(H)$, $q = q(H)$, $s = s(H)$,
$p' = p(H')$, $q' = q(H')$ and $s' = s(H')$.\\
(i) The following limit exists in $\mbbQ \cup \{\infty\}$:
\[
\lim_{n\to\infty} \frac{\big|\mbS_n(\mcA', H')\big|}{\big|\mbS_n(\mcA, H)\big|}. 
\]
(ii) Suppose that $r = 2$.
\begin{itemize}
\item[(a)] If $p - q < p' - q'$ or if $p - q = p' - q'$ and $p > p'$, then
\[
\lim_{n\to\infty} \frac{\big|\mbS_n(\mcA', H')\big|}{\big|\mbS_n(\mcA, H)\big|} \ = \ 0.
\]
\item[(b)] If $p - q = p - q'$ and $p = p'$ then there is a rational number $a > 0$,
depending only on $\mcA$, $\mcA'$, $H$, $H'$ and the vocabulary, such that
\[
\lim_{n\to\infty} \frac{\big|\mbS_n(\mcA', H')\big|}{\big|\mbS_n(\mcA, H)\big|} \ = \ a.
\]
\end{itemize}
(iii) Suppose that $r > 2$ and let $\beta(x,y,z)$ be as in 
Corollary~\ref{estimate of the number of structures for given p,q,s when r>2}.
If any one of the two conditions
\begin{itemize}
\item[] $p - q < p' - q'$, or
\item[] $p - q = p' - q'$ and $\beta(p,q,s) > \beta(p',q',s')$
\end{itemize}
hold, then
\[
\lim_{n\to\infty} \frac{\big|\mbS_n(\mcA', H')\big|}{\big|\mbS_n(\mcA, H)\big|} \ = \ 0.
\]
\end{prop}

\noindent
{\bf Proof.}
(i) From Proposition~\ref{estimate of cardinality of S-n(A)}
it follows that there are integers $C, C' > 0$ and
polynomials $\lambda(x)$, $\lambda'(x)$ with integer coefficients, 
depending only on $\mcA$, $\mcA'$, $H$, $H'$ and the vocabulary, 
such that
\[
\frac{\big| \mbS_n(\mcA', H') \big|}{\big| \mbS_n(\mcA, H) \big|} \ \sim \ 
\frac{C' \binom{n}{p'}}{C \binom{n}{p}} \exp_2\big( \lambda'(n) - \lambda(n) \big).
\]
Depending on whether the leading term in the polynomial $\lambda'(n) - \lambda(n)$ has
positive degree and negative coefficient, positive degree and positive coefficient, or is constant, 
\[\exp_2\big( \lambda'(n) - \lambda(n) \big)\]
approaches 0, $\infty$, or a positive real as $n \to \infty$, respectively.
In the first case 
\[\big| \mbS_n(\mcA', H') \big| \Big/ \big| \mbS_n(\mcA, H) \big|\] 
approaches 0.
In the second case it approaches $\infty$.
In the third case, when $\lambda'(n) - \lambda(n)$ is constant,
we get the conclusion by considering whether $p > p'$, $p=p'$ or $p < p'$.

(ii) Suppose that $r = 2$.
Then Corollary~\ref{estimate of the number of structures for given p,q when r=2} says
that for some positive integers $C$ and $C'$, depending only on
$\mcA$, $\mcA'$, $H$, $H'$ and the vocabulary, we have
\begin{align*}
\frac{\big| \mbS_n(\mcA', H') \big|}{\big| \mbS_n(\mcA, H) \big|} \ \sim \ 
\frac{C' \binom{n}{p'}}{C \binom{n}{p}} 
\exp_2\Big( 2k\big[ (p - q) - (p' - q') \big]n \ + \ k\big[(p')^2 - (p)^2\big] \ + \ l[p - p'] \Big).
\end{align*}
From this we immediately get claims~(a) and~(b).

(iii) Suppose that $r > 2$.
Then Corollary~\ref{estimate of the number of structures for given p,q,s when r>2} implies
that for some  positive integers $C$ and $C'$, depending only on
$\mcA$, $\mcA'$, $H$, $H'$ and the vocabulary, we have
\begin{align*}
\frac{\big| \mbS_n(\mcA', H') \big|}{\big| \mbS_n(\mcA, H) \big|} \ \sim \ 
\frac{C' \binom{n}{p'}}{C \binom{n}{p}} 
\exp_2\bigg( &kr\big[(p-q) \ - \ (p'-q')\big] n^{r-1} \\
&+ \ 
\big[\beta(p',q',s') \ - \ \beta(p,q,s)\big] n^{r-2}
\ + \ \mcO\big(n^{r-3}\big) \bigg).
\end{align*}
So if $p-q < p'-q'$ or if $p-q = p'-q'$ and $\beta(p,q,s) > \beta(p',q',s')$,
then this quotient approaches 0 as $n \to \infty$.
\hfill $\square$
\\

\noindent
{\bf \em For the rest of this section, whenever we denote structures by $\mcA$ or $\mcA'$, 
sometimes with indices, we assume that they have no fixed point. 
Also, whenever we denote groups by $H$ or $H'$, sometimes with indices, we assume that
they have no fixed point.} 
Sometimes these assumptions are repeated and sometimes they are not necessary. 

For different subgroups $H$ and $H'$ of $\Aut(\mcA)$ the sets 
$\mbS_n(\mcA, H)$ and $\mbS_n(\mcA, H')$ may have nonempty intersections,
which complicates the analysis of an asymptotic estimate of the cardinality of a
union like $\bigcup_{i=1}^m \mbS_n(\mcA, H_i)$. 
However, it turns out that for subgroups $H$ and $H'$ of $\Aut(\mcA)$, either
$\mbS_n(\mcA, H) = \mbS_n(\mcA, H')$ or 
$\big| \mbS_n(\mcA, H) \cap \mbS_n(\mcA, H') \big|$ is negligibly small for large enough $n$.
The results~\ref{consequence of H and H' being equivalent over A} 
--~\ref{moving out the summation sign} make this statement precise.

\begin{defin}\label{defin of equivalence relation between groups}{\rm
Suppose that $\mcA \in \mbS$ and that $H$ and $H'$ are subgroups of $\Aut(\mcA)$.
We write $H \approx_\mcA H'$ if there is an automorphism $g \in \Aut(\mcA)$ such
that, for every $t = 1, \ldots, r-1$, and every orbit $O$ of $H$ on $A^t$,
\[g(O) = \{(g(a_1), \ldots, g(a_t)) : (a_1, \ldots, a_t) \in O \}\]
is an orbit of $H'$ on $A^t$. 
}\end{defin}

\noindent
Observe that $\approx_\mcA$ is an equivalence
relation on the set of subgroups of $\Aut(\mcA)$.

\begin{lem}\label{consequence of H and H' being equivalent over A}
Suppose that $\mcA \in \mbS$ and that $H$ and $H'$ are subgroups of $\Aut(\mcA)$.
If $H \approx_\mcA H'$ then $\mbS_n(\mcA, H) = \mbS_n(\mcA, H')$.
\end{lem}

\noindent
{\bf Proof.}
Suppose that $H \approx_\mcA H'$.
Recall from the discussion after the statement of 
Proposition~\ref{estimate of cardinality of S-n(A)}
that $\mbS_n(\mcA, H)$ is the disjoint union of all sets of the form
\[\mbS_n(\mcA_X, H) \ = \ \{ \mcM \in \mbS_n(\mcA, H) : \mcM \uhrc \Spt^*(\mcM) = \mcA_X \},\]
where $X \subseteq [n]$, $|X| = |A|$, $\mcA_X$ has universe $X$ and $\mcA_X \cong \mcA$;
and similarly for $H'$.
Therefore it suffices to prove that for all such $X \subseteq [n]$ and $\mcA_X$ we have
$\mbS_n(\mcA_X, H) = \mbS_n(\mcA_X, H')$.
By~(\ref{union of all S-n(A-X, Pi)}),
\[
\mbS_n(\mcA_X, H) \ = \ \bigcup_{\Pi_1, \ldots, \Pi_{r-1}} \mbS_n(\mcA_X, \Pi_1, \ldots, \Pi_{r-1}),
\]
where the union ranges over all sequences $\Pi_1, \ldots, \Pi_{r-1}$ 
of $(\mcA_X, H)$-partitions (see Definition~\ref{definition of sequence of (A,H)-partitions})
and 
$\mbS_n(\mcA_X, \Pi_1, \ldots, \Pi_{r-1})$
was defined in~(\ref{definition of S-n-A-X-Pi}). The same holds for $H'$.
Hence it suffices to show that if $\Pi_1, \ldots, \Pi_{r-1}$ is a sequence of $(\mcA_X, H)$-partitions,
then $\Pi_1, \ldots, \Pi_{r-1}$ is a sequence of $(\mcA_X, H')$-partitions.

So suppose that $\Pi_1, \ldots, \Pi_{r-1}$ is a sequence of $(\mcA_X, H)$-partitions and hence
there is an isomorphism $f : \mcA \to \mcA_X$ such that, for each $t = 1, \ldots, r-1$,
$\Pi_t$ is the set of orbits of $H_f = \{f\sigma f^{-1} : \sigma \in H\}$ on $X^t$.
As we assume that $H \approx_\mcA H'$, there is an automorphism $g \in \Aut(\mcA)$
such that, for every $t = 1, \ldots, r-1$ and every orbit $O$ of $H$ on $A^t$, 
$g(O)$ is an orbit of $H'$ on $A^t$. 
It follows that $f' = fg^{-1} : \mcA \to \mcA_X$ is an isomorphism and for
each $t$ and each orbit $O'$ of $H'$ on $A^t$, $g^{-1}(O')$ is an orbit of $H$ on $A^t$.
Consequently, for each $t$, $\Pi_t$ is the set of orbits of 
$H_{f'} = \{f' \sigma (f')^{-1} : \sigma \in H'\}$ on $X^t$, so
$\Pi_1, \ldots, \Pi_{r-1}$ is a sequence of $(\mcA_X, H')$-partitions.
\hfill $\square$

\begin{lem}\label{first consequence of H and H' not being equivalent over A}
Suppose that $\mcA \in \mbS$ and that $H$ and $H'$ are subgroups of $\Aut(\mcA)$
such that $H \not\approx_\mcA H'$.
Let $X \subseteq [n]$, $|X| = |A|$ and let $\mcA_X$ have universe $X$ and $\mcA_X \cong \mcA$.
If $\Pi_1, \ldots, \Pi_{r-1}$ is a sequence of $(\mcA_X, H)$-partitions
and $\Pi'_1, \ldots, \Pi'_{r-1}$ is a sequence of $(\mcA_X, H')$-partitions,
then $(\Pi_1, \ldots, \Pi_{r-1}) \neq (\Pi'_1, \ldots, \Pi'_{r-1})$.
\end{lem}

\noindent
{\bf Proof.}
Suppose that $H \not\approx_\mcA H'$, $\Pi_1, \ldots, \Pi_{r-1}$ is a sequence of $(\mcA_X, H)$-partitions
and $\Pi'_1, \ldots, \Pi'_{r-1}$ is a sequence of $(\mcA_X, H')$-partitions.
Towards a contradiction, assume that 
$(\Pi_1, \ldots, \Pi_{r-1}) = (\Pi'_1, \ldots, \Pi'_{r-1})$.
Then there are isomorphisms $f : \mcA \to \mcA_X$ and $f' : \mcA \to \mcA_X$
such that, for every $t = 1, \ldots, r-1$, 
$\Pi_t$ is the set of orbits of $H_f = \{f\sigma f^{-1} : \sigma \in H\}$ on $X^t$ and
$\Pi_t$ is also the set of orbits of $H_{f'} = \{f' \sigma (f')^{-1} : \sigma \in H'\}$ on $X^t$.
So $H_f$ and $H_{f'}$ have the same orbits on $X^t$, for each $t$. 
It follows that $g = (f')^{-1}f : \mcA \to \mcA$ is an automorphism such that 
for every $t = 1, \ldots, r-1$ and every orbit $O$ of $H$ on $A^t$, $g(O)$ is an orbit of $H'$ on $A^t$.
Hence $H \approx_\mcA H'$ which contradicts our assumption.
\hfill $\square$

\begin{lem}\label{second consequence of H and H' not being equivalent over A}
Suppose that $\mcA \in \mbS$ and that $H$ and $H'$ are subgroups of $\Aut(\mcA)$
such that $H \not\approx_\mcA H'$.
Suppose that $X \subseteq [n]$, $|X| = |A|$ and that $\mcA_X$ is a structure with 
universe $X$ such that $\mcA_X \cong \mcA$.
If $\Pi_1, \ldots, \Pi_{r-1}$ is a sequence of $(\mcA_X, H)$-partitions and
$\Pi'_1, \ldots, \Pi'_{r-1}$ is a sequence of $(\mcA_X, H')$-partitions, then
\[
\frac{\big| \mbS_n(\mcA_X, \Pi_1, \ldots, \Pi_{r-1}) \ \cap \ \mbS_n(\mcA_X, \Pi'_1, \ldots, \Pi'_{r-1}) \big|}
{\big| \mbS_n(\mcA_X, \Pi_1, \ldots, \Pi_{r-1}) \ \cup \ \mbS_n(\mcA_X, \Pi'_1, \ldots, \Pi'_{r-1}) \big|}
\ \leq \ \varepsilon(n),
\]
where $\varepsilon(n) \to 0$ as $n \to \infty$ and the function $\varepsilon : \mbbN \to \mbbR$ only depends
on $\mcA$, $H$, $H'$ and the vocabulary.
\end{lem}

\noindent
{\bf Proof.}
The assumptions and
Lemma~\ref{first consequence of H and H' not being equivalent over A}
imply that $(\Pi_1, \ldots, \Pi_{r-1}) \neq (\Pi'_1, \ldots, \Pi'_{r-1})$.
Lemma~\ref{an unusual property} is applicable to each one of the sequences
$\Pi_1, \ldots, \Pi_{r-1}$ and $\Pi'_1, \ldots, \Pi'_{r-1}$.
Now observe that the proof of 
Lemma~\ref{a negligible intersection} works out in exactly the same way even if 
$\Pi_1, \ldots, \Pi_{r-1}$ is a sequence of $(\mcA_X, H)$-partitions and
$\Pi'_1, \ldots, \Pi'_{r-1}$ is a sequence of $(\mcA_X, H')$-partitions;
the proof of Lemma~\ref{a negligible intersection} only uses the assumption that the sequences 
$\Pi_1, \ldots, \Pi_{r-1}$ and $\Pi'_1, \ldots, \Pi'_{r-1}$ are different.
Hence Lemma~\ref{second consequence of H and H' not being equivalent over A} is proved.
\hfill $\square$

\begin{rem}\label{auxilliary before corollary to second consequence of H and H' not being equivalent over A}{\rm
If $\mcA \in \mbS$ and $H \subseteq \Aut(\mcA)$ is a subgroup,
then, by Lemma~\ref{estimate of S-n(A-X,H)} and the argument between 
Proposition~\ref{estimate of cardinality of S-n(A)} and Lemma~\ref{estimate of S-n(A-X,H)}, 
\[\big| \mbS_n(\mcA, H) \big| \ \sim \ C \binom{n}{|A|}\big| \mbS_n(\mcA_X, H) \big|, \]
where $C$ is a constant that depends only on $\mcA$, $H$ and the vocabulary,
$X \subseteq [n]$, $\mcA_X$ is a structure with universe $X$ such that $\mcA_X \cong \mcA$
(and $\mbS_n(\mcA_X, H)$ is as defined in~(\ref{definition of S-n-A-X-H})).
}\end{rem}

\begin{cor}\label{corollary to second consequence of H and H' not being equivalent over A}
Suppose that $\mcA \in \mbS$ and that $H$ and $H'$ are subgroups of $\Aut(\mcA)$
such that $H \not\approx_\mcA H'$.
Then
\[
\frac{\big| \mbS_n(\mcA, H) \ \cap \ \mbS_n(\mcA, H') \big|}
{\big| \mbS_n(\mcA, H) \ \cup \ \mbS_n(\mcA, H') \big|}
\ \leq \ \varepsilon(n),
\]
where $\varepsilon(n) \to 0$ as $n \to \infty$ and the function $\varepsilon$ only depends
on $\mcA$, $H$, $H'$ and the vocabulary.
\end{cor}

\noindent
{\bf Proof.}
Suppose that $\mcA \in \mbS$ and that $H$ and $H'$ are subgroups of $\Aut(\mcA)$
such that $H \not\approx_\mcA H'$.
By Remark~\ref{auxilliary before corollary to second consequence of H and H' not being equivalent over A},
it suffices to prove that there is a function $\varepsilon(n)$, depending only on $\mcA$, $H$ and
the vocabulary, such that $\lim_{n\to\infty} \varepsilon(n) = 0$ and for every
$X \subseteq [n]$ and $\mcA_X$ as above,
\[
\frac{\big| \mbS_n(\mcA_X, H) \ \cap \ \mbS_n(\mcA_X, H') \big|}
{\big| \mbS_n(\mcA_X, H) \ \cup \ \mbS_n(\mcA_X, H') \big|}
\ \leq \ \varepsilon(n).
\]
Recall from~(\ref{union of all S-n(A-X, Pi)}) that 
\[\mbS_n(\mcA_X, H) \ = \ \bigcup_{\Pi_1, \ldots, \Pi_{r-1}} \mbS_n(\mcA_X, \Pi_1, \ldots, \Pi_{r-1})\]
where the union ranges over all sequences of $(\mcA_X, H)$-partitions.
Given $X$ and $\mcA_X$ there is a finite bound $\alpha$,
depending only on $\mcA$, $H$, $H'$ and the vocabulary,
such that there are at most $\alpha$
sequences $\Pi_1, \ldots, \Pi_{r-1}$ of $(\mcA_X, H)$-partitions 
and at most $\alpha$ sequences $\Pi'_1, \ldots, \Pi'_{r-1}$ of $(\mcA_X, H')$-partitions.
Therefore the bound we are looking for is a fixed multiple of the bound given
by Lemma~\ref{second consequence of H and H' not being equivalent over A}.
\hfill $\square$

\begin{lem}\label{moving out the summation sign}
Suppose that $\mcA \in \mbS$ and that $H_i$, $i = 1, \ldots, m$, are subgroups of $\Aut(\mcA)$
such that if $i \neq j$, then $H_i \not\approx_\mcA H_j$ and 
$\big| \mbS_n(\mcA, H_i) \big| \Big/ \big| \mbS_n(\mcA, H_j) \big|$ converges to
a positive rational number.
Then
\[ \Bigg| \bigcup_{i=1}^m \mbS_n(\mcA, H_i) \Bigg| \ \sim \ 
\sum_{i=1}^m \big| \mbS_n(\mcA, H_i) \big|. \]  
\end{lem}

\noindent
{\bf Proof.}
From Corollary~\ref{corollary to second consequence of H and H' not being equivalent over A}
it follows that if $i \neq i'$ then
\begin{equation}\label{intersection smaller that o(1) times sum}
\big| \mbS_n(\mcA, H_i) \ \cap \ \mbS_n(\mcA, H_{i'}) \big| \ \leq \ 
o(1) \Big( \big| \mbS_n(\mcA, H_i) \big| \ + \ \big| \mbS_n(\mcA, H_{i'}) \big| \Big),
\end{equation}
where the bound $o(1)$ depends only on $\mcA$, $H_1, \ldots, H_m$ and the vocabulary.
Now the assumption that $\big| \mbS_n(\mcA, H_i) \big| \Big/ \big| \mbS_n(\mcA, H_{i'}) \big|$ 
converges to a positive rational number 
and~(\ref{intersection smaller that o(1) times sum}) implies that
if $i \neq i'$, then
\[
\big| \mbS_n(\mcA, H_i) \ \cap \ \mbS_n(\mcA, H_{i'}) \big| \ \leq \ 
o(1) \big| \mbS_n(\mcA, H_i) \big|,
\]
for some bound $o(1)$ which depends only on $\mcA$, $H_1, \ldots, H_m$ and the vocabulary.
If we let $\mbU_n$ be the union of all intersections
\[
\mbS_n(\mcA, H_i) \ \cap \ \mbS_n(\mcA, H_{i'})
\]
where $\{i,i'\}$ range over all subsets of $[m]$ with cardinality 2,
then we get, for every $i$,
\[
\big| \mbS_n(\mcA, H_i) \big| \ - \ \big| \mbU_n \big| \ \geq \ \big(1 - o(1)\big)\big| \mbS_n(\mcA, H_i) \big|,
\]
where the bound $o(1)$ depends only on $\mcA$ and $H_1, \ldots, H_m$.
Now we get
\[
\Bigg| \bigcup_{i=1}^m \mbS_n(\mcA, H_i) \Bigg| \ \geq \ 
\sum_{i=1}^m \Big( \big| \mbS_n(\mcA, H_i) \big| \ - \ \big| \mbU_n \big| \Big) \ \geq \ 
\big(1 - o(1)\big) \sum_{i=1}^m \big| \mbS_n(\mcA, H_i) \big|.
\]
Since also 
\[ \Bigg| \bigcup_{i=1}^m \mbS_n(\mcA, H_i) \Bigg| \ \leq \ 
\sum_{i=1}^m \big| \mbS_n(\mcA, H_i) \big|
\]
the proof of the lemma is finished.
\hfill $\square$

\begin{prop}\label{general quotient proposition}
Let $\mcA_1, \ldots, \mcA_m, \mcA'_1, \ldots, \mcA'_{m'} \in \mbS$ 
be such that none of them has any fixed point.
Suppose that for every $i = 1, \ldots, m$ and $j = 1, \ldots, l_i$, 
$H_{i,j}$ is a subgroup of $\Aut(\mcA_i)$ without any fixed point and that
for every $i = 1, \ldots, m'$ and $j = 1, \ldots, l'_i$
$H'_{i,j}$ is a subgroup of $\Aut(\mcA'_i)$ without any fixed point.
Then the following limit exists in $\mbbQ \cup \{\infty\}$:
\begin{equation}\label{general limit of quotients of big unions}
\lim_{n\to\infty} \frac{\Big| \bigcup_{i=1}^{m'} \bigcup_{j=1}^{l'_i} \mbS_n(\mcA'_i, H'_{i,j}) \Big|}
{\Big| \bigcup_{i=1}^{m} \bigcup_{j=1}^{l_i} \mbS_n(\mcA_i, H_{i,j}) \Big|}.
\end{equation}
\end{prop}

\noindent
{\bf Proof.}
By just removing, if necessary, some $\mcA_i$ or $\mcA'_i$ and relabelling the rest of the structures 
and groups, we may assume that $\mcA_i \not\cong \mcA_j$ if $i\neq j$ and
$\mcA'_i \not\cong \mcA'_j$ if $i \neq j$.
Also, by Lemma~\ref{consequence of H and H' being equivalent over A},
we may assume that $H_{i,j} \not\approx_{\mcA_i} H_{i,j'}$ if $j \neq j'$ and that
$H'_{i,j} \not\approx_{\mcA'_i} H'_{i,j'}$ if $j \neq j'$.

By Proposition~\ref{comparisson of structures with different support-substructures}~(i),
for all $1 \leq i \leq m$ and all $1 \leq j,j' \leq l_i$,
\[
\big| \mbS_n(\mcA_i, H_{i,j'}) \big| \Big/ \big| \mbS_n(\mcA_i, H_{i,j}) \big|
\]
converges to a rational number or approaches infinity, as $n \to \infty$.
The same holds for all $1 \leq i \leq m'$, all $1 \leq j,j' \leq l'_i$
and $\big| \mbS_n(\mcA'_i, H'_{i,j'}) \big| \Big/ \big| \mbS_n(\mcA'_i, H'_{i,j}) \big|$.
Therefore it suffices to prove~(\ref{general limit of quotients of big unions})
under the assumption that for all $1 \leq i \leq m$ and all $1 \leq j,j' \leq l_i$,
$\big| \mbS_n(\mcA_i, H_{i,j'}) \big| \Big/ \big| \mbS_n(\mcA_i, H_{i,j}) \big|$
converges to a {\em positive rational number} and for all
$1 \leq i \leq m'$ and all $1 \leq j,j' \leq l'_i$,
$\big| \mbS_n(\mcA'_i, H'_{i,j'}) \big| \Big/ \big| \mbS_n(\mcA'_i, H'_{i,j}) \big|$
converges to a {\em positive rational number}.

From our assumptions we have
$\mbS_n(\mcA_i, H_{i,j}) \cap \mbS_n(\mcA_{i'}, H_{i',j'}) = \es$ if $i \neq i'$
(and the same for $\mcA'_i$, $\mcA'_j$, $H'_{i,j}$ and $H'_{i,j'}$).
By applying Lemma~\ref{moving out the summation sign} and the assumptions,
we now get
\begin{align*}
&\frac{\Big| \bigcup_{i=1}^{m'} \bigcup_{j=1}^{l'_i} \mbS_n(\mcA'_i, H'_{i,j}) \Big|}
{\Big| \bigcup_{i=1}^{m} \bigcup_{j=1}^{l_i} \mbS_n(\mcA_i, H_{i,j}) \Big|} \ = \ 
\frac{\sum_{i=1}^{m'} \sum_{j=1}^{l'_i} \big| \mbS_n(\mcA'_i, H'_{i,j}) \big|}
{\sum_{i=1}^{m} \sum_{j=1}^{l_i} \big| \mbS_n(\mcA_i, H_{i,j}) \big|} \\
= \ &\frac{\big| \mbS_n(\mcA'_1, H_{1,1}) \big|}
{\sum_{i=1}^{m} \sum_{j=1}^{l_i} \big| \mbS_n(\mcA_i) \big|} 
\ + \ \ldots \ + \ 
\frac{\big| \mbS_n(\mcA'_{m'}, H_{m' l_{m'}}) \big|}
{\sum_{i=1}^{m} \sum_{j=1}^{l_i} \big| \mbS_n(\mcA_i, H_{i,j}) \big|} \\
= \ &\Bigg( \sum_{i=1}^{m} \sum_{j=1}^{l_i}
\frac{\big| \mbS_n(\mcA_i, H_{i,j}) \big|}{\big| \mbS_n(A'_1, H'_{1,1}) \big|} \Bigg)^{-1}
\ + \ \ldots \ + \ 
\Bigg( \sum_{i=1}^{m} \sum_{j=1}^{l_i}
\frac{\big| \mbS_n(\mcA_i, H_{i,j}) \big|}{\big| \mbS_n(A'_{m'}, H'_{m', l_{m'}}) \big|} \Bigg)^{-1}.
\end{align*}
Note that Proposition~\ref{estimate of cardinality of S-n(A)} implies that, for all
$i$, $j$ and all sufficiently large $n$, 
$\big| \mbS_n(\mcA_i, H_{i,j}) \big| > 0$, and similarly for
and $\mcA'_i$ and $H'_{i,j}$, so we do not divide by zero in the above expression if $n$ is large enough.
By Proposition~\ref{comparisson of structures with different support-substructures}~(i),
for every choice of $i, i', j$ and $j'$,
$\big| \mbS_n(\mcA_i, H_{i,j}) \big| \Big/ \big| \mbS_n(\mcA'_{i'}, H'_{i',j'}) \big|$ 
converges to a rational number or approaches $\infty$.
This implies~(\ref{general limit of quotients of big unions}) so the proposition is proved.
\hfill $\square$

\begin{prop}\label{comparisson of two subgroups}
Let $G$ and $G'$ be finite groups.
Then the following limit exists in $\mbbQ \cup \{\infty\}$:
\[
\lim_{n\to\infty} \frac{\big| \big\{ \mcM \in \mbS_n : G' \leq \Aut(\mcM) \big\} \big|}
{\big| \big\{ \mcM \in \mbS_n : G \leq \Aut(\mcM) \big\} \big|}.
\]
\end{prop}

\noindent
{\bf Proof.}
Let $G$ and $G'$ be finite groups.
Lemma~\ref{different sets equals a union of S(A,H)} implies that there are finitely many
\begin{align*}
&\text{structures } \ \mcA_1, \ldots, \mcA_m, \mcA'_1, \ldots, \mcA'_{m'} \in \mbS,\\
&\text{subgroups $H_{i,1}, \ldots, H_{i,l_i} \subseteq \Aut(\mcA_i)$, for $i = 1, \ldots, m$, and} \\
&\text{subgroups $H'_{i,1}, \ldots, H'_{i,l'_i} \subseteq \Aut(\mcA'_i)$, for $i = 1, \ldots, m'$,}
\end{align*}
such that the following hold:
\begin{itemize}
\item[(i)] None of the permutation groups 
$\Aut(\mcA_i)$, $\Aut(\mcA'_i)$, $H_{i,j}$ or $H_{i,j}$ has any fixed point.

\item[(ii)] $\big| \{\mcM \in \mbS_n : G \leq \Aut(\mcM) \big\} \big| \sim 
\Big| \bigcup_{i=1}^{m} \bigcup_{j=1}^{l_i} \mbS_n(\mcA_i, H_{i,j}) \Big|$ as $n \to \infty$.

\item[(iii)] $\big| \{\mcM \in \mbS_n : G' \leq \Aut(\mcM) \big\} \big| \sim 
\Big| \bigcup_{i=1}^{m'} \bigcup_{j=1}^{l'_i} \mbS_n(\mcA'_i, H'_{i,j}) \Big|$ as $n \to \infty$.
\end{itemize}
Hence it suffices to prove that
\begin{equation*}
\lim_{n\to\infty} \frac{\Big| \bigcup_{i=1}^{m'} \bigcup_{j=1}^{l'_i} \mbS_n(\mcA'_i, H'_{i,j}) \Big|}
{\Big| \bigcup_{i=1}^{m} \bigcup_{j=1}^{l_i} \mbS_n(\mcA_i, H_{i,j}) \Big|}
\end{equation*}
exists in $\mbbQ \cup \{\infty\}$.
But this follows immediately from 
Proposition~\ref{general quotient proposition}.
\hfill $\square$
\\

\noindent
By the definition of $\mbS_n(\mcA, H)$ 
(Definition~\ref{definition of S-n(A)}), 
for every $\mcM \in \mbS_n(\mcA, H)$,
$\Aut(\mcM) \uhrc \Spt^*(\mcM)$ has a subgroup $H_f$ such that $H_f \cong_P H$.
The next lemma shows that for almost all $\mcM \in \mbS(\mcA, H)$ 
any such $H_f$ has the same orbits as $\Aut(\mcM) \uhrc \Spt^*(\mcM)$.

\begin{lem}\label{automorhism group almost surely has same orbits as H}
Suppose that $\mcA \in \mbS$ has no fixed point and that $H$ is a 
subgroup of $\Aut(\mcA)$ without any fixed point.
There is a function $\varepsilon(n)$, depending only on $\mcA$, $H$ and
the vocabulary, such that $\lim_{n\to\infty} \varepsilon(n) = 0$ and
the proportion of $\mcM \in \mbS_n(\mcA, H)$ with the following property is at most $\varepsilon(n)$:
\begin{itemize}
\item[$(*)$] For some isomorphism $f : \mcA \to \mcM \uhrc \Spt^*(\mcM)$ such
that $H_f = \{f\sigma f^{-1} : \sigma \in H\}$ is a subgroup of $\Aut(\mcM) \uhrc \Spt^*(\mcM)$,
there is $t \in \{1, \ldots, r-1\}$ such that 
the orbits of $\Aut(\mcM) \uhrc \Spt^*(\mcM)$ on $\Spt^*(\mcM)^t$ are {\em not} the same
as the orbits of $H_f$ on $\Spt^*(\mcM)^t$.
\end{itemize}
\end{lem}

\noindent
{\bf Proof.}
Let $\mcM \in \mbS_n(\mcA, H)$,  $X = \Spt^*(\mcM)$ and $\mcA_X = \mcM \uhrc X$, so $\mcA_X \cong \mcA$.
Moreover, let $f : \mcA \to \mcA_X$ be an isomorphism and assume that
$H_f = \{f \sigma f^{-1} : \sigma \in H\}$ is a subgroup of $\Aut(\mcM) \uhrc X$.
Suppose that for some $t \in \{1, \ldots, r-1\}$
the orbits of $\Aut(\mcM) \uhrc X$ on $X^t$ are not the same as the orbits of $H_f$ on $X^t$.
It follows that $\Aut(\mcM) \uhrc X$ has fewer orbits on $X^t$ than $H_f$.
Hence there is a subgroup $H'$ of $\Aut(\mcA)$ such that $H \subseteq H'$,
$H'$ has fewer orbits than $H$ on $A^t$ and $\mcM \in \mbS_n(\mcA, H')$.
It follows that $H' \not\approx_\mcA H$ and that
\[\mcM \in \mbS_n(\mcA, H') \cap \mbS_n(\mcA, H).\]
Now Corollary~\ref{corollary to second consequence of H and H' not being equivalent over A}
implies that 
\[
\big| \mbS_n(\mcA, H) \ \cap \ \mbS_n(\mcA, H') \big| \ \leq \ 
\varepsilon(n) \big| \mbS_n(\mcA, H) \ \cup \ \mbS_n(\mcA, H') \big|,
\]
where $\varepsilon(n) \to 0$ as $n \to \infty$ and $\varepsilon(n)$ only depends
on $\mcA$, $H$, $H'$ and the vocabulary.
Since $H$ is a subgroup of $H'$ we have
\[\mbS_n(\mcA, H') \ \subseteq \ \mbS_n(\mcA, H),\]
which implies that
\begin{equation}\label{H' is negligibly small compared with H} 
\big| \mbS_n(\mcA, H') \big| \ \leq \ \varepsilon(n) \big| \mbS_n(\mcA, H) \big|. 
\end{equation}
We have proved that if $\mcM \in \mbS_n(\mcA, H)$ and satisfies~$(*)$ then
$\mcM \in \mbS_n(\mcA, H')$ for some subgroup $H'$ of $\Aut(\mcA)$ such 
that~(\ref{H' is negligibly small compared with H}) holds.
As the number of subgroups $H'$ of $\Aut(\mcA)$ is finite and depends only on $\mcA$ the lemma follows.
\hfill $\square$

\begin{defin}\label{definition of full permutation group}{\rm
Suppose that $\mcA \in \mbS$ has no fixed point and that $H$ is a 
subgroup of $\Aut(\mcA)$ without any fixed point.
For $\mcM \in \mbS_n(\mcA, H)$ we say that $H$ is the {\em full automorphism group of $\mcM$}
if for every isomorphism $f : \mcA \to \mcM \uhrc \Spt^*(\mcM)$ such
that $H_f = \{f \sigma f^{-1} : \sigma \in H\}$ is a subgroup of $\Aut(\mcM) \uhrc \Spt^*(\mcM)$
we have $H_f = \Aut(\mcM) \uhrc \Spt^*(\mcM)$.
}\end{defin}

\begin{lem}\label{0-1 law for probability that the subgroup is the whole automorphism group}
Suppose that $\mcA \in \mbS$ has no fixed point and that $H$ is a 
subgroup of $\Aut(\mcA)$ without any fixed point.
The proportion of $\mcM \in \mbS_n(\mcA, H)$ such that $H$ is the full automorphism group of $\mcM$
converges to either 0 or 1 as $n \to \infty$.
\end{lem}

\noindent
{\bf Proof.}
By Lemma~\ref{automorhism group almost surely has same orbits as H},
it suffices to consider $\mcM \in \mbS_n(\mcA, H)$ with the following property:
\begin{itemize}
\item[] For every isomorphism $f : \mcA \to \mcM \uhrc \Spt^*(\mcM)$ such that
$H_f$ is a subgroup of $\Aut(\mcM) \uhrc \Spt^*(\mcM)$, 
$H_f$ and $\Aut(\mcM) \uhrc \Spt^*(\mcM)$ have the same orbits on $\Spt^*(\mcM)^t$
for all $t = 1, \ldots, r-1$.
\end{itemize}
For such $\mcM$ the question whether 
there is $g \in \Aut(\mcM) \uhrc \Spt^*(\mcM)$ such that $g \notin H_f$
depends only on the isomorphism type of $\mcA$, $H$ and the isomorphism $f : \mcA \to \mcM \uhrc \Spt^*(\mcM)$.
In fact, it depends only on the isomorphism type of $\mcA$ and $H$.
For if $f$ and $f'$ are isomorphisms from $\mcA$ to $\mcM \uhrc \Spt^*(\mcM)$,
$H_f = \Aut(\mcM) \uhrc \Spt^*(\mcM)$ and $g \in \Aut(\mcM) \uhrc \Spt^*(\mcM)$ does not
belong to $H_{f'}$, then, since $|H_f| = |H_{f'}|$ (because 
$f'f^{-1}$ is an isomorphism from $H_f$ to $H_{f'}$ as permutation groups), we get
$|H_{f'}| < |\Aut(\mcM)| = |H_f| = |H_{f'}|$, which is impossible.
\hfill $\square$

\begin{lem}\label{0-1 law for probability that the subgroup is the whole automorphism group, a sharpening}
Suppose that $\mcA \in \mbS$ has no fixed point and that $H \subseteq \Aut(\mcA)$ is a subgroup
without any fixed point. For every group $G \leq H$, the proportion of $\mcM \in \mbS_n(\mcA, H)$
such that $G \cong \Aut(\mcM)$ converges to either 0 or 1 as $n \to \infty$.
\end{lem}

\noindent
{\bf Proof.}
Suppose that $\mcA \in \mbS$ has no fixed point, that $H \subseteq \Aut(\mcA)$ is a subgroup
without any fixed point and $G \leq H$.
Since $\Aut(\mcM) \cong \Aut(\mcM) \uhrc \Spt^*(\mcM)$ for every $\mcM \in \mbS$,
Lemma~\ref{0-1 law for probability that the subgroup is the whole automorphism group} 
implies that the proportion of $\mcM \in \mbS_n(\mcA, H)$ such that $H \cong \Aut(\mcM)$
converges to either 0 or 1 as $n \to \infty$.
If $G \cong H$ it follows that the proportion of $\mcM \in \mbS_n(\mcA, H)$ such that $G \cong \Aut(\mcM)$
converges to either 0 or 1 as $n \to \infty$.
If $G$ is isomorphic to a proper subgroup of $H$ then, since $H \leq \Aut(\mcM)$ for
every $\mcM \in \mbS_n(\mcA, H)$, it follows that $G \not\cong \Aut(\mcM)$ for every $\mcM \in \mbS_n(\mcA, H)$.
\hfill $\square$

\begin{prop}\label{proportion of structures with subgroup G that have automorphism group G}
If $G$ is a finite group then there is a rational number $0 \leq a \leq 1$ such that
\[
\lim_{n\to\infty} \frac{\big| \big\{ \mcM \in \mbS_n : G \cong \Aut(\mcM) \big\} \big|}
{\big| \big\{ \mcM \in \mbS_n : G \leq \Aut(\mcM) \big\} \big|} \ = \ a.
\]
\end{prop}

\noindent
{\bf Proof.}
Let $G$ be a finite group.
By Lemma~\ref{different sets equals a union of S(A,H)}, there are finitely many
\begin{align*}
&\text{structures } \ \mcA_1, \ldots, \mcA_m \in \mbS \text{ and} \\
&\text{subgroups $H_{i,1}, \ldots, H_{i,l_i} \subseteq \Aut(\mcA_i)$, for $i = 1, \ldots, m$} 
\end{align*}
such that:
\begin{itemize}
\item[(i)] None of the permutation groups $\Aut(\mcA_i)$ or $H_{i,j}$ has a fixed point.

\item[(ii)] $G \leq H_{i,j}$ for all $i$ and $j$.

\item[(iii)] $\big| \{\mcM \in \mbS_n : G \leq \Aut(\mcM) \big\} \big| \sim 
\Big| \bigcup_{i=1}^{m} \bigcup_{j=1}^{l_i} \mbS_n(\mcA_i, H_{i,j}) \Big|$ as $n \to \infty$.
\end{itemize}
Lemma~\ref{0-1 law for probability that the subgroup is the whole automorphism group, a sharpening}
says that for every $\mcA_i$ and every $H_{i,j}$ the proportion of $\mcM \in \mbS_n(\mcA_i, H_{i,j})$
such that $G \cong \Aut(\mcM)$ converges to either 0 or 1.
Let $(\mcA'_i, H'_{i,j})$, $i = 1, \ldots, m'$, $j = 1, \ldots, l'_i$, enumerate
the pairs $(\mcA_i, H_{i,j})$ such that the proportion of $\mcM \in \mbS_n(\mcA_i, H_{i,j})$
for which $G \cong \Aut(\mcM)$ converges to 1.
Then
\[
\frac{\big| \big\{ \mcM \in \mbS_n : G \cong \Aut(\mcM) \big\} \big|}
{\big| \big\{ \mcM \in \mbS_n : G \leq \Aut(\mcM) \big\} \big|} \ \sim \ 
\frac{\Big| \bigcup_{i=1}^{m'} \bigcup_{j=1}^{l'_i} \mbS_n(\mcA'_i, H'_{i,j}) \Big|}
{\Big| \bigcup_{i=1}^{m} \bigcup_{j=1}^{l_i} \mbS_n(\mcA_i, H_{i,j}) \Big|} \ \to \ a
\ \ \text{ as } n \to \infty
\]
for some rational $0 \leq a \leq 1$, by Proposition~\ref{general quotient proposition}.
\hfill $\square$

\begin{prop}\label{comparisson of two groups as automorphism groups}
Let $G$ and $G'$ be finite groups.
Then the following limit exists in $\mbbQ \cup \{\infty\}$:
\[
\lim_{n\to\infty} \frac{\big| \big\{ \mcM \in \mbS_n : G' \cong \Aut(\mcM) \big\} \big|}
{\big| \big\{ \mcM \in \mbS_n : G \cong \Aut(\mcM) \big\} \big|}.
\]
\end{prop}

\noindent
{\bf Proof.}
By Lemma~\ref{different sets equals a union of S(A,H)}, there are finitely many
\begin{align*}
&\text{structures } \ \mcA_1, \ldots, \mcA_m, \mcA'_1, \ldots, \mcA'_{m'} \in \mbS,\\
&\text{subgroups $H_{i,1}, \ldots, H_{i,l_i} \subseteq \Aut(\mcA_i)$, for $i = 1, \ldots, m$, and} \\
&\text{subgroups $H'_{i,1}, \ldots, H'_{i,l'_i} \subseteq \Aut(\mcA'_i)$, for $i = 1, \ldots, m'$,}
\end{align*}
such that:
\begin{itemize}
\item[(i)] None of the permutation groups 
$\Aut(\mcA_i)$, $\Aut(\mcA'_i)$, $H_{i,j}$ or $H_{i,j}$ has any fixed point.

\item[(ii)] $G \leq H_{i,j}$ and $G' \leq H'_{i,j}$ for all $i, j$.

\item[(iii)] $\big| \{\mcM \in \mbS_n : G \leq \Aut(\mcM) \big\} \big| \sim 
\Big| \bigcup_{i=1}^{m} \bigcup_{j=1}^{l_i} \mbS_n(\mcA_i, H_{i,j}) \Big|$ as $n \to \infty$.

\item[(iv)] $\big| \{\mcM \in \mbS_n : G' \leq \Aut(\mcM) \big\} \big| \sim 
\Big| \bigcup_{i=1}^{m'} \bigcup_{j=1}^{l'_i} \mbS_n(\mcA'_i, H'_{i,j}) \Big|$ as $n \to \infty$.
\end{itemize}
As in the proof of 
Proposition~\ref{proportion of structures with subgroup G that have automorphism group G}
we now use 
Lemma~\ref{0-1 law for probability that the subgroup is the whole automorphism group, a sharpening}.
So let $(\mcA^*_i, H^*_{i,j})$, $i = 1, \ldots, m^*$, $j = 1, \ldots, l^*_i$, 
enumerate all pairs $(\mcA_i, H_{i,j})$ such that
the proportion of $\mcM \in \mbS_n(\mcA_i, H_{i,j})$ for which $G \cong \Aut(\mcM)$ converges to 1. 
Similarly, let $(\mcA^+_i, H^+_{i,j})$, $i = 1, \ldots, m^+$, $j = 1, \ldots, l^+_i$, 
enumerate all pairs $(\mcA'_i, H'_{i,j})$ such that
the proportion of $\mcM \in \mbS_n(\mcA'_i, H'_{i,j})$ for which $G' \cong \Aut(\mcM)$ converges to 1. 
Then
\[
\frac{\big| \big\{ \mcM \in \mbS_n : G' \cong \Aut(\mcM) \big\} \big|}
{\big| \big\{ \mcM \in \mbS_n : G \cong \Aut(\mcM) \big\} \big|} \ \sim \ 
\frac{\Big| \bigcup_{i=1}^{m^+} \bigcup_{j=1}^{l^+_i} \mbS_n(\mcA^+_i, H^+_{i,j}) \Big|}
{\Big| \bigcup_{i=1}^{m^*} \bigcup_{j=1}^{l^*_i} \mbS_n(\mcA^*_i, H^*_{i,j}) \Big|},
\]
where, by Proposition~\ref{general quotient proposition}, 
the right side converges to a rational number or tends to infinity as $n \to \infty$.
\hfill $\square$

\begin{rem}\label{remark on limit for other complexity measures than subgroups}{\rm
By the use of Lemma~\ref{different sets equals a union of S(A,H)}
and arguments similar to those already carried out one can prove that 
$a_n / b_n$ converges to a rational number as $n \to \infty$ if, for example,
$a_n = |\mbS_n(\spt \geq k_1)|$ and $b_n = |\mbS_n(\spt \geq k_2)|$.
}\end{rem}

\section{Logical limit laws}\label{Logical limit laws}

\noindent
The main results of this section are 
Theorems~\ref{zero-one law for A,H}
and~\ref{limit law for first-order logic}, where the later implies
Theorem~\ref{main theorem about limit laws for first-order logic}.
In Remark~\ref{remark on limit laws} we observe that we do not have
a zero-one law in Theorems~\ref{limit law for first-order logic} 
or~\ref{main theorem about limit laws for first-order logic}.

\begin{theor}\label{zero-one law for A,H}
Suppose that $\mcA \in \mbS$ has no fixed point and let $H$ be a subgroup of $\Aut(\mcA)$
without any fixed point. Then $\mbS(\mcA, H)$ has a zero-one law.
\end{theor}

\noindent
Before proving Theorem~\ref{zero-one law for A,H} we derive:

\begin{theor}\label{limit law for first-order logic}
(i) Let $\mcA_1, \ldots, \mcA_m \in \mbS$ be such that none of them has any fixed point.
Suppose that for every $i = 1, \ldots, m$ and $j = 1, \ldots, l_i$, $H_{i,j}$ is a subgroup
of $\Aut(\mcA_i)$ without any fixed point.
Then $\displaystyle \bigcup_{i=1}^m \bigcup_{j=1}^{l_i} \mbS(\mcA_i, H_{i,j})$
has a limit law. Moreover the limit always belongs to a finite set of rational numbers which is
determined only by the structures $\mcA_i$ and the permutation groups $H_{i,j}$.\\
(ii) The following sets have limit laws: for every finite group $G$, 
$\{\mcM \in \mbS : G \cong \Aut(\mcM)\}$ and
$\{\mcM \in \mbS : G \leq \Aut(\mcM)\}$, and for every integer $m \geq 2$, 
$\mbS(\spt^* = m)$, $\mbS(\spt \geq m)$ and $\mbS(\spt^* \geq m)$. 
Moreover, in each case there is a finite set $Q \subseteq \mbbQ$ such that
for every sentence $\varphi$, the proportion of structures in which $\varphi$ is true
converges to a number in $Q$.
\end{theor}

\noindent
{\bf Proof.} Part (i) is immediate from Theorem~\ref{zero-one law for A,H}
and Proposition~\ref{general quotient proposition}.
For part~(ii) let $\mbX$ be any one of the sets of structures considered.
By Lemmas~\ref{spt-star equals a union of S(A,H)}
and~\ref{different sets equals a union of S(A,H)} 
(and in one case the proof
of Proposition~\ref{proportion of structures with subgroup G that have automorphism group G}),
there are structures
$\mcA_1, \ldots, \mcA_l \in \mbS$ without any fixed point
and for every $i = 1, \ldots, l$ and $j = 1, \ldots, l_i$, a subgroup $H_{i,j}$
of $\Aut(\mcA_i)$ without any fixed point, such that if $\mbX_n = \mbX \cap \mbS_n$ then
\[
\big|\mbX_n\big| \ \sim \  \Bigg| \bigcup_{i=1}^l \bigcup_{j=1}^{l_i} \mbS_n(\mcA_i, H_{i,j}) \Bigg|.
\]
Now part~(ii) follows from part~(i).
\hfill $\square$

\subsection{Proof of Theorem~\ref{zero-one law for A,H}}\label{proof of zero-one law}

Suppose that $\mcA \in \mbS$ has no fixed point and let $H$ be a subgroup of $\Aut(\mcA)$
without any fixed point. We will define a theory $T_{\mcA, H}$ and show that it is consistent and complete
and that for every finite subset $\Delta \subseteq T_{\mcA, H}$, the proportion of $\mcM \in \mbS_n(\mcA, H)$
such that $\mcM \models \Delta$ approaches 1 as $n \to \infty$.
Then the compactness theorem implies that if $T_{\mcA, H} \models \varphi$ then 
the proportion of $\mcM \in \mbS_n(\mcA, H)$ in which $\varphi$ is true approaches 1 as $n \to \infty$;
otherwise that proportion approaches 0.
The argument follows a well known path, using so-called extension axioms.
What makes things more complicated here, compared with Fagin's original proof 
\cite{EF, Fag76} that for every $k \in \mbbN$ the proportion
of $\mcM \in \mbS_n$ satisfying the $k$-extension axiom approaches 1 as $n \to \infty$,
is that all members of $\mbS_n(\mcA, H)$ have nonempty support (of cardinality $|A|$).

To make the main ideas of the argument more transparent, to avoid heavy formulations and notation and to
expose more clearly how our argument differs from the ``standard argument''
(in \cite{EF, Fag76} for example), we will prove Theorem~\ref{zero-one law for A,H}
in the special case when the vocabulary consists of only one binary relation symbol $R$.
It is straightforward to generalise the argument to any finite relational vocabulary
with at least one relation symbol of arity at least 2, but it comes at the expense of longer
definitions and heavier notation and formulations in order to keep track of all data and
how it can be combined.
Moreover, the arguments can be modified to the case when we assume that some (possibly all)
relation symbols are always interpreted as irreflexive relations, or when we assume that some (possibly all)
relation symbols are always interpreted as irreflexive {\em and} symmetric relations.

For any structure $\mcM$ and formula $\varphi(x)$ we let
\[\varphi(\mcM) \ = \ \{ a \in M : \mcM \models \varphi(a) \}.\]
Let $p = |A|$, $A = \{a_1, \ldots, a_p\}$ and let $x_1, \ldots, x_p$ be distinct variables.
For $a_i, a_j \in A$, let $a_i \approx a_j$ mean that $a_i$ and $a_j$ belong to the same orbit of $H$.
Let $\chi_\mcA(x_1, \ldots, x_p)$ be a formula which describes the isomorphism type of $\mcA$.
More precisely, $\chi_\mcA(x_1, \ldots, x_p)$ is the conjunction of all formulas of the following form:
$x_i \neq x_j$ for $i \neq j$; $R(x_i, x_j)$ if $\mcA \models R(a_i, a_j)$;
and $\neg R(x_i, x_j)$ if $\mcA \models \neg R(a_i, a_j)$.

We will define formulas $\theta(x)$, $\xi(x,y)$ such that the proportion of $\mcM \in \mbS_n(\mcA, H)$
such that the following hold approaches 1 as $n \to \infty$:
\begin{itemize}
\item[(I)] For every $a \in M$, $\mcM \models \theta(a)$ if and only if $a \in \Spt^*(\mcM)$.
\item[(II)] $\mcM$ satisfies following sentence, denoted $\psi$:
\begin{align*}
\exists x_1, \ldots, x_m \Bigg(
&\chi_\mcA(x_1, \ldots, x_m) \ \wedge \ 
\forall y \bigg[ \theta(y) \ \longleftrightarrow \ \bigvee_{i=1}^m y = x_i \bigg] \ \wedge \\ 
&\bigwedge_{a_i \approx a_j} \xi(x_i, x_j) \ \wedge \ 
\bigwedge_{a_i \not\approx a_j} \neg\xi(x_i, x_j) \ \wedge \\ 
&\bigwedge_{i,j=1}^m \forall y 
\bigg[ \neg\theta(y) \ \wedge \ \xi(x_i, x_j) \ \longrightarrow \\ 
&\Big(\xi(y, x_i) \longleftrightarrow \xi(y, x_j)\Big) \ \wedge \ 
\Big(\xi(x_i, y) \longleftrightarrow \xi(x_j, y)\Big) \bigg]
\Bigg).
\end{align*}
\end{itemize}

\noindent
It is straightforward to define, for every $k \in \mbbN$, a sentence
$\varphi_k$ such that for every, possibly infinite, structure $\mcM$:
\begin{itemize}
\item[(III)] If $\mcM \models \varphi_k$ then the following hold:
\begin{itemize}
\item[(a)] $\big| \theta(\mcM) \big| = m$ and the relation defined by $\xi(x,y)$ restricted
to $\theta(\mcM)$ is an equivalence relation.

\item[(b)] For every choice of $i \in \{0,1\}$, $B \subseteq M \setminus \theta(\mcM)$ with $|B| = k$,
sets $E, E'$ of $\xi$-equivalence classes on $\theta(\mcM)$ 
and $Y, Y' \subseteq B$, there is $c \in M \setminus \theta(\mcM)$
such that
\begin{itemize}
\item[(i)] $\mcM \models R(c,c) \ \Longleftrightarrow \ i = 1$,
\item[(ii)] for every $a \in \theta(\mcM)$, \\
$\mcM \models R(c,a) \ \Longleftrightarrow $ $a$ belongs to some class in $E$ and \\ 
$\mcM \models R(a,c) \ \Longleftrightarrow $ $a$ belongs to some class in $E'$, and
\item[(iii)] for every $b \in B$, \\
$\mcM \models R(c,b) \ \Longleftrightarrow \ b \in Y$ and \\
$\mcM \models R(b,c) \ \Longleftrightarrow \ b \in Y'$.
\end{itemize}
\end{itemize}
\end{itemize}

\noindent
We call $\varphi_k$ the {\em $k$-extension axiom}.
Assuming that we have $\theta(x)$, $\xi(x,y)$, $\psi$ and $\varphi_k$, $k \in \mbbN$, as above, we let
\[ T_{\mcA, H} \ = \ \{\psi\} \ \cup \ \{\varphi_k : k \in \mbbN \}. \]
Note that every model of $T_{\mcA, H}$ is infinite.
We postpone the proof that $T_{\mcA, H}$ has a model to the end of the argument.
By \L os' and Vaught's categoricity theorem (see \cite{Roth}, Theorem~8.5.1, for instance), $T_{\mcA, H}$ is complete
if we can prove that it is countably categorical.

\begin{lem}\label{countable categoricity of T}
If $\mcM_1$ and $\mcM_2$ are countable models of $T_{\mcA, H}$ then $\mcM_1 \cong \mcM_2$.
\end{lem}

\noindent
{\bf Proof.}
This is a standard back-and-forth argument, so we only sketch it.
Suppose that $\mcM_1$ and $\mcM_2$ are countable models of $T_{\mcA, H}$.
Since both $\mcM_1$ and $\mcM_2$ satisfy $\psi$ it follows that $\theta(\mcM_1)$
and $\theta(\mcM_2)$ are finite and that there is
an isomorphism $f_0 : \mcM_1 \uhrc \theta(\mcM_1) \to \mcM_2 \uhrc \theta(\mcM_2)$
such that for all $a,b \in \theta(\mcM_1)$, 
$\mcM_1 \models \xi(a,b)$ if and only if $\mcM_2 \models \xi(f_0(a),f_0(b))$.
Therefore it suffices to prove the following statement:
\\

\noindent
{\bf Claim.} {\em Suppose that $B_1 \subseteq M_1 \setminus \theta(\mcM_1)$, 
$B_2 \subseteq M_2 \setminus \theta(\mcM_2)$ and that
\[f : \mcM_1 \uhrc (\theta(\mcM_1) \cup B_1) \to \mcM_2 \uhrc (\theta(\mcM_2) \cup B_2)\]
is an isomorphism.
If $c_1 \in M_1 \setminus (\theta(\mcM_1) \cup B_1)$ (or $c_2 \in M_2 \setminus (\theta(\mcM_2) \cup B_2)$)
then there is $c_2 \in M_2 \setminus (\theta(\mcM_2) \cup B_2)$
(or $c_1 \in M_1 \setminus (\theta(\mcM_1) \cup B_1)$) such that $f$ can be extended to
an isomorphism from 
$\mcM_1 \uhrc (\theta(\mcM_1) \cup B_1 \cup \{c_1\})$
to
$\mcM_2 \uhrc (\theta(\mcM_2) \cup B_2 \cup \{c_2\})$.}
\\

\noindent
Let $k = |B_1| = |B_2|$.
The claim follows in a straigthforward way since $\mcM$ and $\mcN$ are models of $\{\psi, \varphi_k\}$.
\hfill $\square$
\\

\noindent
It remains to show that there are $\theta(x)$ and $\xi(x,y)$ such that, for every $k$, the proportion
of $\mcM \in \mbS_n(\mcA, H)$ that satisfy~(I) and the sentences $\psi$ and $\varphi_k$ 
approaches 1 as $n \to \infty$.
Recall that, with the notation from Section~\ref{asymptotic estimates}, 
\[ \mbS_n(\mcA, H) \ = \ \bigcup_X \bigcup_{\mcA_X} \mbS_X(\mcA_X, H), \]
where the first union ranges over all subsets of $[n]$ with cardinality $m = |A|$, and for each such
subset $X$, the second union ranges over all structures $\mcA_X$ with universe $X$ that are isomorphic
to $\mcA$. As observed in that section, if $X \neq X'$ then $\mbS_n(\mcA_X, H)$ is disjoint from
$\mbS_n(\mcA_{X'}, H)$. Moreover, if $\mcA_X$ and $\mcA'_X$ are different structures with universe
$X$ then $\mbS_n(\mcA_X, H)$ is disjoint from $\mbS_n(\mcA'_X, H)$.
Recall our assumption in this proof that there is only one relation symbol $R$ and it has arity $r = 2$.
In Section~\ref{asymptotic estimates} we also saw (recall~(\ref{union of all S-n(A-X, Pi)}))
that for each $\mbS_n(\mcA_X, H)$,
\[ \mbS_n(\mcA_X, H) \ = \ \bigcup_{\Pi_1} \mbS_n(\mcA_X, \Pi_1), \]
where the union ranges over all $(\mcA_X, H)$-partitions $\Pi_1$ of $X$; 
see Definition~\ref{definition of sequence of (A,H)-partitions}. 
The number of $(\mcA_X, H)$-partitions of $X$ is the same for every sufficiently large
$n$, every $X \subseteq[n]$ with $|X| = m$ and every $\mcA_X \cong \mcA$.
Therefore it suffices to prove that there are $\theta(x)$, $\xi(x,y)$
and, for every $k$, $0 < \alpha < 1$ such
that for every $X \subseteq [n]$
with $|X| = m$, every $\mcA_X \cong \mcA$ with universe $X$ and every $(\mcA_X, H)$-partition $\Pi_1$,
the proportion of $\mcM \in \mbS_n(\mcA_X, \Pi_1)$ that satisfy~(I)
and the sentences $\psi$ and $\varphi_k$
is at least $1 - \alpha^n$.

{\bf \em For the rest of this section we fix $X \subseteq [n]$ with $|X| = |A| = m$
and $\mcA_X \cong \mcA$ with universe $X$.}
The results below refer to all large enough $n$ with respect to other parameters that occur.

\begin{defin}\label{definition of k-extension property}{\rm
We say that $\mcM \in \mbS_n(\mcA_X, \Pi_1)$ has the {\em $k$-extension property}
if~(III)~(b) holds when `$\theta(\mcM)$' is replaced by `$X$' and
`$\xi$-equivalence classes' with `parts of the partition $\Pi_1$ (of $X$)'.
}\end{defin}

\begin{lem}\label{upper bound on structures with with fixed support structure and partition}
For every $k \in \mbbN$
there is $0 < \alpha_k < 1$, depending only on $k$ and $\mcA$,  
such that the proportion of 
$\mcM \in \mbS_n(\mcA_X, \Pi_1)$ which does not have the $k$-extension property is at most $\alpha_k^n$.
\end{lem}

\noindent
{\bf Proof.}
Recall that $\mbS_n(\mcA_X, \Pi_1) \subseteq \mbT_n(\mcA_X, \Pi_1)$,
where $\mbT_n(\mcA_X, \Pi_1)$ was defined in~(\ref{definition of T-n-A-X-Pi}) in 
Section~\ref{asymptotic estimates}.
From Lemma~\ref{estimate for Pi-1,...,Pi-r-1} we know that
$\big| \mbS_n(\mcA_X, \Pi_1) \big| \Big/ \big| \mbT_n(\mcA_X, \Pi_1) \big| \to 1$ as $n \to \infty$,
so it suffices to prove the statement of the lemma for $\mbT_n(\mcA_X, \Pi_1)$
in the place of $\mbS_n(\mcA_X, \Pi_1)$. The reason for doing this is that $\mbT_n(\mcA_X, \Pi_1)$
is easier to work with because its members do not have the constraint that the support of the structure is
exactly $X$ (but from the arguments in Section~\ref{asymptotic estimates}
we know that for every $\mcM \in \mbT_n(\mcA_X, \Pi_1)$, $X \subseteq \Spt^*(\mcM)$).

Let $\mcM \in \mbT_n(\mcA_X, \Pi_1)$.
A subset $B \subseteq [n] \setminus X$ of cardinality $k$ can be chosen in no more than $n^k$ ways.
Once $B \subseteq [n] \setminus X$ with $|B| = k$ is fixed, the number of 
ways to choose $i \in \{0,1\}$, $E, E' \subseteq \Pi_1$
and $Y, Y' \subseteq B$ is bounded where the bound depends only on $k$ and $\mcA$. 
Therefore it suffices to show, for an arbitrary fixed $B \subseteq [n] \setminus X$ with $|B| = k$
and an arbitrary choice of $i \in \{0,1\}$, $E, E' \subseteq \Pi_1$ and $Y, Y' \subseteq B$, 
that the proportion of $\mcM \in \mbT_n(\mcA_X, \Pi_1)$ such that there is {\em no} $c \in M$
such that the conjunction 
of~(i)--(iii) of~(III) is satisfied
is at most $\alpha_k^n$ for some constant $0 < \alpha_k < 1$ that depends only on $k$ and $\mcA$.

For arbitrary $c \in [n] \setminus (X \cup B)$ we estimate the probability 
that at least one of~(i)--(iii) of~(III) fails.
We consider $\mbT_n(\mcA_X, \Pi_1)$ as a probability space by giving each member the same probability.
From the definition of $\mbT_n(\mcA_X, \Pi_1)$ we see that the probability that 
$\mcM \in \mbT_n(\mcA_X, \Pi_1)$ satisfies~(i)--(iii) of~(III) is
$2^{-\beta}$ for some $\beta > 0$ depending only on $|X|$ and $|B|$, and independently of what the case
is for other elements than $c$ in $[n] \setminus (X \cup B)$. 
The probability that, for every $c \in [n] \setminus (X \cup B)$,
the conjunction of~(i)--(iii) does not hold is therefore
\[\big(1 - 2^{-\beta}\big)^{n-|X \cup B|}.\]
As $B$ can be chosen in at most $n^k$ ways
it follows that the probability that the conjunction of~(i)--(iii) 
is not satisfied in $\mcM \in \mbT_n(\mcA_X, \Pi_1)$ is at most 
$\alpha_k^n$ for some $0 < \alpha_k < 1$ that depends only on $k$ and $\mcA$.
\hfill $\square$
\\

\noindent
Remember that $m = |A| = |X|$.
Let $\theta(x)$ denote the following formula:
\begin{align*}
\exists y_1, \ldots, y_{m-1} \Bigg(
&\bigwedge_{i=1}^{m-1} x \neq y_i \ \wedge \ \bigwedge_{i \neq j} y_i \neq y_j \ \wedge \\
&\forall z \Bigg[ \Big( z \neq x \ \wedge \ \bigwedge_{i=1}^{m-1} z \neq y_i \Big) \ \longrightarrow \ 
\Big( R(x,z) \ \longleftrightarrow \ R(y_1, z) \Big) \Bigg] \Bigg).
\end{align*}

\begin{lem}\label{definability of support in structures with fixed support}
Suppose that $\mcM \in \mbS_n(\mcA_X, \Pi_1)$ has the $2$-extension property.
Then for all $a \in M$, $a \in X = \Spt^*(\mcM)$ if and only if $\mcM \models \theta(a)$.
\end{lem}

\noindent
{\bf Proof.}
Suppose that $\mcM \models \theta(a)$.
Then there are distinct $b_1, \ldots, b_{m-1} \in M$ different from $a$ such that
for all $c$ different from $b_1, \ldots, b_{m-1}$ and from $a$,
\[ \mcM \models R(a,c) \ \longleftrightarrow \ R(b_1, c).\]
As $\mcM$ has the $2$-extension property this is only possible if $a, b_1 \in X = \Spt^*(\mcM)$
and $a$ and $b_1$ belong to the same part of $\Pi_1$.

Now suppose that $a \in \Spt^*(\mcM) = X$.
Let $b_1, \ldots, b_{m-1}$ be such that
\[X = \{a, b_1, \ldots, b_{m-1}\} \]
and $a$ and $b_1$ belong to the same part of $\Pi_1$.
By the definition of $\mbS_n(\mcA_X, \Pi_1)$, there is an automorphism of $\mcM$ which sends
$a$ to $b_1$ and fixes every element outside of $X$ and therefore we must have
\[ 
\mcM \models \forall z \Bigg[ \Big( z \neq a \ \wedge \ 
\bigwedge_{i=1}^{m-1} z \neq b_i \Big) \ \longrightarrow \
\Big( R(a,z) \ \longleftrightarrow \ R(b_1,z) \Big) \Bigg].
\]
\hfill $\square$
\\

\noindent
Let $\xi(x_1, x_2)$ be the formula
\[ \forall z \Big( \neg \theta(z) \ \longrightarrow \ 
\big[ R(z,x_1) \ \longleftrightarrow \ R(z,x_2) \big] \Big). \]

\begin{lem}\label{formula defining the orbits}
Suppose that $\mcM \in \mbS_n(\mcA_X, \Pi_1)$ has the $2$-extension property.
Then for all $a_1, a_2 \in X = \Spt^*(\mcM)$, $a_1$ and $a_2$ belong to the same part of $\Pi_1$
if and only if $\mcM \models \xi(a_1, a_2)$.
\end{lem}

\noindent
{\bf Proof.}
Suppose that $a_1, a_2 \in X = \Spt^*(\mcM)$ and $a_1$ and $a_2$ belong to the same part of $\Pi_1$.
By the definition of $\mbS_n(\mcA_X, \Pi_1)$, for every $c \in M \setminus X$ there is an automorphism
which sends $a_1$ to $a_2$ and fixes every element outside of $X$. 
From Lemma~\ref{definability of support in structures with fixed support} 
it follows that $\mcM \models \xi(a_1, a_2)$.

Now suppose that $a_1, a_2 \in X = \Spt^*(\mcM)$ and $\mcM \models \xi(a_1, a_2)$.
From Lemma~\ref{definability of support in structures with fixed support}
it follows that for all $c \in M \setminus X$ 
\[ \mcM \models R(c,a_1) \Longleftrightarrow R(c,a_2). \]
Since we assume that $\mcM$ has the $2$-extension property this is only possible if $a_1$ and
$a_2$ belong to the same part of $\Pi_1$.
\hfill $\square$
\\

\noindent
According to the arguments before Definition~\ref{definition of k-extension property}
and the compactness theorem,
the following corollary concludes the proof of Theorem~\ref{zero-one law for A,H}.

\begin{cor}\label{corollary giving that the k-extension axioms are almost always satisfied}
For every $k \in \mbbN$, there is $0 < \alpha < 1$, depending only on $k$ and $\mcA$, such
that the proportion of $\mcM \in \mbS_n(\mcA_X, \Pi_1)$ that satisfy~(I)
and the sentences $\psi$ and $\varphi_l$ for $l = 0, \ldots, k$
is at least $1 - \alpha^n$.
\end{cor}

\noindent
{\bf Proof.}
Let $k' = \max(2,k,m)$ (where $m = |A|$).
By Lemma~\ref{upper bound on structures with with fixed support structure and partition}
there is $0 < \alpha < 1$, depending only on $k'$ and $\mcA$ such that
the proportion of $\mcM \in \mbS_n(\mcA_X, \Pi_1)$ with the $l$-extension property for every $l \leq k'$
is at least $1 - \alpha^n$. From Lemmas~\ref{definability of support in structures with fixed support}
and~\ref{formula defining the orbits} it follows that all $\mcM \in \mbS_n(\mcA_X, \Pi_1)$
with the $l$-extension property for every $l \leq k'$ satisfy~(I)
and the sentences $\psi$ and $\varphi_l$ for $l = 0, \ldots, k$.
\hfill $\square$

\begin{rem}\label{remark on limit laws}{\rm
Let $\mbS'$ be any one of the sets of structures in part~(ii) of
Theorem~\ref{limit law for first-order logic}
and let $\mbS'_n = \mbS' \cap \mbS_n$. We assume that if a finite group $G$ is
involved in the definition of $\mbS'$ then $G$ is nontrival.
We will show that $\mbS'$ does not satisfy a zero-one law.
By Lemmas~\ref{spt-star equals a union of S(A,H)}
and~\ref{different sets equals a union of S(A,H)} (and in one case the proof of 
Proposition~\ref{proportion of structures with subgroup G that have automorphism group G}),
there are mutually nonisomorphic $\mcA_1, \ldots, \mcA_l \in \mbS$ without any fixed point and,
for $i = 1, \ldots, l$ and $j = 1, \ldots, l_i$,
subgroups $H_{i,j} \subseteq \Aut(\mcA_i)$ without any fixed point such that 
\[ 
\big| \mbS'_n \big| \ \sim \ 
\Bigg| \bigcup_{i=1}^l \bigcup_{j=1}^{l_i} \mbS_n(\mcA_i, H_{i,j}) \Bigg|.
\]
If $\mbS'$ is $\{\mcM \in \mbS : G \leq \Aut(\mcM)\}$ or $\{\mcM \in \mbS : G \cong \Aut(\mcM)\}$,
then we may also assume that $G \leq H_{i,j}$ or $G \cong H_{i,j}$, respectively, for all $i$ and $j$.

Now observe the following: Suppose that $\mcA \in \mbS$ has no fixed point and that
$H$ is a subgroup of $\Aut(\mcA)$ without any fixed point.
Let $\mcA'$ and $\mcA''$ have the same universe $A$ as $\mcA$ and assume that for every
relation symbol $R$, $R^{\mcA'} = \es$ and $R^{\mcA''} = A^i$ if $R$ is $i$-ary.
Then $H$ is a subgroup of $\Aut(\mcA')$ and of $\Aut(\mcA'')$ and, from 
Proposition~\ref{estimate of cardinality of S-n(A)}, it follows that
\[
\frac{\big| \mbS_n(\mcA', H) \big|}{\big| \mbS_n(\mcA, H) \big|} \quad \text{ and } \quad 
\frac{\big| \mbS_n(\mcA'', H) \big|}{\big| \mbS_n(\mcA, H) \big|}
\]
converge to the same $c \in \mbbQ$ as $n \to \infty$.
From the assumption that $\mbS'$ is one of the sets of structures in part~(ii) of
Theorem~\ref{limit law for first-order logic} (and $G$ is assumed to be nontrivial) it follows that
there must be $i, i', j, j'$ such that $\mcA_i \not\cong \mcA_{i'}$ and both
$\big| \mbS_n(\mcA_i, H_{i,j}) \big| \big/ \big| \mbS'_n \big|$ and
$\big| \mbS_n(\mcA_{i'}, H_{i',j'} \big| \big/ \big| \mbS'_n \big|$
converge to positive numbers $c$ and $c'$ as $n \to \infty$.
With the help of the formula $\theta$ from the proof of 
Theorem~\ref{zero-one law for A,H} one can easily construct a sentence
$\varphi$ which, in almost all $\mcM \in \mbS'$,
expresses that ``$\mcM \uhrc \Spt^*(\mcM) \cong \mcA_i$''. 
Then the proportion of $\mcM \in \mbS'_n$ in which $\varphi$ is true converges to
some number $0 < d < 1$.
}\end{rem}

\section{Unlabelled structures}\label{Unlabelled structures}

\noindent
The main result of this final section is Theorem~\ref{unlabelled limit laws}, which implies
Theorem~\ref{main theorem about unlabelled structures},
which says that Theorems~\ref{main theorem about comparissons between different groups}
and~\ref{main theorem about limit laws for first-order logic} hold also for unlabelled structures.

\begin{defin}\label{definition of isomorphism counting notation}{\rm
(i) For every $\mcM \in \mbS$, let
$[\mcM] \ = \ \{ \mcN \in \mbS : \mcN \cong \mcM \}$.\\
(ii) For every $\mbX \subseteq \mbS$, let 
$[\mbX] = \{[\mcM] : \mcM \in \mbX\}$.\\
(iii) We say that a set $\mbX \subseteq \mbS$ is {\em closed under isomorphism}
if $\mcM \in \mbX$, $\mcN \in \mbS$ and $\mcN \cong \mcM$ implies that $\mcN \in \mbX$.
}\end{defin}

\noindent
The next lemma is a generalisation of Lemma~4.3.10 in \cite{EF}.

\begin{lem}\label{connection between labelled and unlabelled structures}
If $\mbX_n \subseteq \mbS_n$ is closed under isomorphism then
\[ \liso \mbX_n \riso n! = 
\sum_{\pi \in Sym_n} \big| \mbS_n(\pi) \cap \mbX_n \big|.\]
\end{lem}

\noindent
{\bf Proof.}
For every $\mcM \in \mbS_n$ and $\pi \in Sym_n$, let $\pi(\mcM)$ denote the 
unique structure $\mcN \in \mbS_n$ such that $\pi$ is an isomorphism from $\mcM$ onto $\mcN$.
Fix an arbitrary $\mcM \in \mbX_n$ and let $H = \Aut(\mcM)$.
Then $H$ is a subgroup of $Sym_n$ and we consider the left cosets of $H$ in $Sym_n$.
Note that for every $\mcN \in \mbX_n$ we have $\mcN \cong \mcM$ if and only if there
is $\pi \in Sym_n$ such that $\pi(\mcM) = \mcN$.
For all $\pi, \sigma \in Sym_n$ we have 
\[
\pi H = \sigma H \ \Longleftrightarrow \ H = \pi^{-1}\sigma H 
\ \Longleftrightarrow \ \pi^{-1}\sigma \in H = \Aut(\mcM) 
\ \Longleftrightarrow \ \pi(\mcM) = \sigma(\mcM)
\]
As we assume that $\mbX_n$ is closed under isomorphism it follows that 
\begin{equation*}
\big| \{\mcN \in \mbX_n : \mcN \cong \mcM\} \big| 
\ = \ 
\text{the number of cosets (the index) of $\Aut(\mcM)$ in $Sym_n$.}
\end{equation*}
Hence
\[
\big| \Aut(\mcM) \big| \cdot \big| \{ \mcN \in \mbX_n : \mcN \cong \mcM \} \big| 
\ = \ \big| Sym_n \big| \ = \ n!,
\]
and, as $\big| \Aut(\mcN) \big| = \big| \Aut(\mcM) \big|$ if $\mcN \cong \mcM$, we get
\begin{align*}
\sum_{\substack{\mcN \in \mbX_n \\ \mcN \cong \mcM}} \big| \Aut(\mcN) \big| 
\ = \ 
\sum_{\substack{\mcN \in \mbX_n \\ \mcN \cong \mcM}} \big| \Aut(\mcM) \big| 
\ = \ 
\big| \{ \mcN \in \mbX_n : \mcN \cong \mcM \} \big| \cdot \big| \Aut(\mcM) \big| \ = \ n!.
\end{align*}
If $\mcM_1, \ldots, \mcM_m$ is a sequence containing exactly one representative from
every isomorphism class that is represented in $\mbX_n$, then $m = \liso \mbX_n \riso$ and
\[ 
\sum_{\mcM \in \mbX_n} \big| \Aut(\mcM) \big| \ = \ 
\sum_{i=1}^m \sum_{\substack{\mcN \in \mbX_n \\ \mcN \cong \mcM_i}} \big| \Aut(\mcN) \big| \ = \ 
\sum_{i=1}^m n! \ = \ \liso \mbX_n \riso \cdot n!.
\]
We also have 
\[ \sum_{\mcM \in \mbX_n} \big| \Aut(\mcM) \big| \ = \ 
\big| \{(\mcM, \pi) : \mcM \in \mbX_n \text{ and } \pi \in \Aut(\mcM) \} \big| 
\ = \ \sum_{\pi \in Sym_n} \big| \mbS_n(\pi) \cap \mbX_n \big|,
\] 
which concludes the proof of the lemma.
\hfill $\square$

\begin{lem}\label{relationship between counting labelled and unlabelled structures}
If $\mbY \subseteq \mbS$ is closed under isomorphism and $p \geq 2$ is fixed, then
\[ \big| \mbS_n(\spt^* \leq p) \ \cap \ \mbY \big| \ \sim \ 
n! \liso \mbS_n(\spt^* \leq p) \ \cap \ \mbY \riso \quad \text{ as } n \to \infty. \]
\end{lem}

\noindent
{\bf Proof.}
For every permutation $\pi$ of $[n]$ and $\mcM \in \mbS_n$, let
$\pi(\mcM)$ denote the unique structure $\mcM' \in \mbS_n$ such that $\pi$ is an isomorphism from
$\mcM$ to $\mcM'$.
If $\mcM \in \mbS_n(\spt^* \leq p)$, $\pi$ is a permutation of $[n]$ and $\pi(\mcM) = \mcM$,
then $\Spt(\pi) \subseteq \Spt^*(\mcM)$. Hence there are at most $p!$ permutations
$\pi$ of $[n]$ such that $\pi(\mcM) = \mcM$.
Since we assume that $\mbY$ is closed under isomorphism we get
\[ \big| \mbS_n(\spt^* \leq p) \ \cap \ \mbY \big| \ \geq \ 
(n! - p!) \liso \mbS_n(\spt^* \leq p) \ \cap \ \mbY \riso. \]
It is also clear that
\[ \big| \mbS_n(\spt^* \leq p) \ \cap \ \mbY \big| \ \leq \ 
n! \liso \mbS_n(\spt^* \leq p) \ \cap \ \mbY \riso. \]
Since $(n! - p!) \sim n!$ as $n \to \infty$, it follows that
\[ \big| \mbS_n(\spt^* \leq p) \ \cap \ \mbY \big| \ \sim \ 
n! \liso \mbS_n(\spt^* \leq p) \ \cap \ \mbY \riso. \]
\hfill $\square$
\\

\begin{prop}\label{automorphisms with large support are unusual for unlabelled structures}
Suppose that $m, t \in \mbbN$, $f_1, \ldots, f_s \in Sym_n$,
$\spt(f_1, \ldots, f_s) = m$ and $t > 2r(m! - 1)m/m! + 1$, where $r \geq 2$ is the maximal
arity of the relation symbols.
Then there is $\lambda > 0$ such that for all sufficiently large $n$,
\[
\frac{\liso \mbS_n(\spt \geq t) \riso}{\liso \mbS_n(f_1, \ldots, f_s) \riso} \ \leq \ 
2^{-\lambda n^{r-1}}.
\]
\end{prop}

\noindent
{\bf Proof.}
Suppose that $m \in \mbbN$, $f_1, \ldots, f_s \in Sym_n$ and
$\spt(f_1, \ldots, f_s) = m$.
Let 
\[ 
\widehat{\mbS}_n(f_1, \ldots, f_s) \ = \ 
\big\{\mcM \in \mbS_n : \mcM \cong \mcN \text{ for some } \mcN \in \mbS_n(f_1, \ldots, f_s) \big\}
\]
and observe that $\big[ \widehat{\mbS}_n(f_1, \ldots, f_s) \big] = \big[ \mbS_n(f_1, \ldots, f_s) \big]$.
By Propositions~\ref{automorphisms with large support are unusual}
and~\ref{support of structure is bounded in terms of support of automorphisms}, 
there are constants $p, \alpha > 0$
such that for all sufficiently large $n$,
\[
\frac{\big| \mbS_n(\spt^* > p) \big|}{\big| \mbS_n(f_1, \ldots, f_s) \big|} \ \leq \ 
2^{-\alpha n^{r-1} \pm \mcO\big(n^{r-2}\big)}.
\]
Since $\big| \mbS_n(f_1, \ldots, f_s) \big| \leq \big| \widehat{\mbS}_n(f_1, \ldots, f_s) \big|$, we get
\[
\frac{\big| \mbS_n(\spt^* > p) \big|}{\big| \widehat{\mbS}_n(f_1, \ldots, f_s) \big|} \ \leq \ 
2^{-\alpha n^{r-1} \pm \mcO\big(n^{r-2}\big)},
\]
which implies that
\begin{equation}\label{almost all members of wide-hat-S-n of the functions}
\big| \widehat{\mbS}_n(f_1, \ldots, f_s) \ \cap \ \mbS_n(\spt^* \leq p) \big| \ \sim \ 
\big| \widehat{\mbS}_n(f_1, \ldots, f_s) \big|.
\end{equation}
Lemma~\ref{relationship between counting labelled and unlabelled structures} with
$\mbY = \widehat{\mbS}_n(f_1, \ldots, f_s)$ gives
\[
\big| \widehat{\mbS}_n(f_1, \ldots, f_s) \ \cap \ \mbS_n(\spt^* \leq p) \big| \ \sim \ 
n! \liso \widehat{\mbS}_n(f_1, \ldots, f_s) \ \cap \ \mbS_n(\spt^* \leq p) \riso.
\]
This and~(\ref{almost all members of wide-hat-S-n of the functions}) gives
\begin{equation}\label{reduction of unlabelled structures to labelled with the given automorphisms}
n! \liso \widehat{\mbS}_n(f_1, \ldots, f_s) \ \cap \ \mbS_n(\spt^* \leq p) \riso \ \sim \ 
\big| \widehat{\mbS}_n(f_1, \ldots, f_s) \big|.
\end{equation}
Suppose that $t > 2r(m! - 1)m/m! + 1$.
By Lemma~\ref{connection between labelled and unlabelled structures}
with $\mbX_n = \mbS_n(\spt \geq t)$ we get
\begin{equation}\label{the transfer to labelled structures}
\liso \mbS_n(\spt \geq t) \riso \cdot n! \ = \ 
\sum_{\pi \in Sym_n} \big| \mbS_n(\pi) \ \cap \ \mbS_n(\spt \geq t) \big|.
\end{equation}
For every $\pi \in Sym_n$, $\big| \mbS_n(\pi) \ \cap \ \mbS_n(\spt \geq t) \big| \ \leq \ 
\big| \mbS_n(\spt \geq t) \big|$ and there are not more than $(t-1)!n^{t-1}$
permutations $\pi \in Sym_n$ such that $\spt(\pi) < t$.
Therefore,
\begin{equation}\label{upper bound on permutations with small support}
\sum_{\substack{\pi \in Sym_n \\ \spt(\pi) < t}} \big| \mbS_n(\pi) \ \cap \ \mbS_n(\spt \geq t) \big|
\ \leq \ 
(t-1)!n^{t-1} \big| \mbS_n(\spt \geq t) \big|.
\end{equation}
If $\pi \in Sym_n$ and $\spt(\pi) \geq t$ then
$\mbS_n(\pi) \ \cap \ \mbS_n(\spt \geq t) = \mbS_n(\pi)$, so we get
\begin{equation}\label{first simplification of a sum}
\sum_{\substack{\pi \in Sym_n \\ \spt(\pi) \geq t}} 
\big| \mbS_n(\pi) \ \cap \ \mbS_n(\spt \geq t) \big| \ = \ 
\sum_{\substack{\pi \in Sym_n \\ \spt(\pi) \geq t}} \big| \mbS_n(\pi) \big|.
\end{equation}
Now we get
\begin{align}\label{getting rid of an intersection}
&\sum_{\pi \in Sym_n} \big| \mbS_n(\pi) \ \cap \ \mbS_n(\spt \geq t) \big| \\ 
\leq \ 
&(t-1)! n^{t-1} \big| \mbS_n(\spt \geq t) \big| \ + \ 
\sum_{\substack{\pi \in Sym_n \\ \spt(\pi) \geq t}} \big| \mbS_n(\pi) \big| 
\qquad \qquad \text{ by (\ref{upper bound on permutations with small support})
and (\ref{first simplification of a sum})} \nonumber \\
\leq \ 
&(t-1)! n^{t-1} \sum_{\substack{\pi \in Sym_n \\ \spt(\pi) \geq t}} \big| \mbS_n(\pi) \big| \ + \ 
\sum_{\substack{\pi \in Sym_n \\ \spt(\pi) \geq t}} \big| \mbS_n(\pi) \big|  
\qquad \quad \text{ by (\ref{upper bound for spt less or equal to m})} \nonumber \\ 
\leq \ 
&2(t-1)! n^{t-1} \sum_{\substack{\pi \in Sym_n \\ \spt(\pi) \geq t}} \big| \mbS_n(\pi) \big|.
\nonumber
\end{align}
Moreover, as $\liso \mbS_n(f_1, \ldots, f_s) \riso \ = \ \liso \widehat{\mbS}_n(f_1, \ldots, f_s) \riso$ we have
\begin{align}
&\frac{\liso \mbS_n(\spt \geq t) \riso}{\liso \mbS_n(f_1, \ldots, f_s) \riso} \ = \ 
\frac{\liso \mbS_n(\spt \geq t) \riso}{\liso \widehat{\mbS}_n(f_1, \ldots, f_s) \riso} 
\label{reduction to expression in the proof for the labelled case}\\ 
\leq \ 
&\frac{\liso \mbS_n(\spt \geq t) \riso}
{\liso \widehat{\mbS}_n(f_1, \ldots, f_s) \ \cap \ \mbS_n(\spt^* \leq p) \riso}  \nonumber \\ 
= \ 
&\frac{n! \cdot \liso \mbS_n(\spt \geq t) \riso}
{n! \cdot \liso \widehat{\mbS}_n(f_1, \ldots, f_s) \ \cap \ \mbS_n(\spt^* \leq p) \riso}  \nonumber \\ 
\sim \ 
&\frac{\sum_{\pi \in Sym_n} \big| \mbS_n(\pi) \ \cap \ \mbS_n(\spt \geq t) \big|}
{\big| \widehat{\mbS}_n(f_1, \ldots, f_s) \big|} 
\qquad \qquad \text{ by (\ref{the transfer to labelled structures}) and 
(\ref{reduction of unlabelled structures to labelled with the given automorphisms})}  \nonumber \\ 
\leq \ 
&\frac{2(t-1)! n^{t-1} \sum_{\substack{\pi \in Sym_n \\ \spt(\pi) \geq t}} \big| \mbS_n(\pi) \big|}
{\big| \mbS_n(f_1, \ldots, f_s) \big|}
\qquad \qquad \text{ by (\ref{getting rid of an intersection})}  \nonumber \\
\leq \ 
&2(t-1) n^{t-1} 
\sum_{\substack{\pi \in Sym_n \\ \spt(\pi) \geq t}}
\exp_2\Bigg( \sum_{i=1}^r k_i \orb^i(\pi) \ - \ \sum_{i=1}^r k_i \orb^i(f_1, \ldots, f_s)\Bigg)  \nonumber \\
&\text{ by (\ref{equality for automorphisms}) and (\ref{upper bound for spt less or equal to m})}.
\nonumber
\end{align}
Since $t > 2r(m! - 1)m/m! + 1$ it follows from the final estimates of the proof of 
Proposition~\ref{automorphisms with large support are unusual}
that there is $\beta > 0$ such that 
\[ 
\sum_{\substack{\pi \in Sym_n \\ \spt(\pi) \geq t}}
\exp_2\Bigg( \sum_{i=1}^r k_i \orb^i(\pi) \ - \ \sum_{i=1}^r k_i \orb^i(f_1, \ldots, f_s)\Bigg) \ \leq \
2^{-\beta n^{r-1} \pm \mcO(n^{r-2})}.
\]
This together with~(\ref{reduction to expression in the proof for the labelled case})
implies that there is $\lambda > 0$ such that
\[
\frac{\liso \mbS_n(\spt \geq t) \riso}{\liso \mbS_n(f_1, \ldots, f_s) \riso} \ \leq \ 
2^{-\lambda n^{r-1}}
\]
for all large enough $n$.
\hfill $\square$

\begin{cor}\label{corollary to automorphisms with large support are unusual in unlabelled case}
Let $m, t \in \mbbN$.\\
(i) If $t > 2r(m! - 1)m/m! + 1$ then
\[
\lim_{n\to\infty} \ \frac{\liso \mbS_n(\spt \geq t) \riso}{\liso \mbS_n(\spt \geq m) \riso} \ = \ 
\lim_{n\to\infty} \ \frac{\liso \mbS_n(\spt \geq t) \riso}{\liso \mbS_n(\spt^* \geq m) \riso} \ = \ 0.
\]
(ii) There is $T > m$ such that 
\begin{align*}
&\lim_{n\to\infty} \frac{\liso \mbS_n(\spt \geq m) \cap \mbS_n(\spt^* \leq T) \riso}
{\liso \mbS_n(\spt \geq m) \riso} \ = \\
&\lim_{n\to\infty} \frac{\liso \mbS_n(\spt^* \geq m) \cap \mbS_n(\spt^* \leq T) \riso}
{\liso \mbS_n(\spt^* \geq m) \riso} \ = \ 1.
\end{align*}
\end{cor}

\noindent
{\bf Proof.}
Part~(i) follows immediately from 
Proposition~\ref{automorphisms with large support are unusual for unlabelled structures}, 
because if $f \in Sym_n$ and $\spt(f) = m$,
then $\mbS_n(f) \subseteq \mbS_n(\spt \geq m)
\subseteq \mbS_n(\spt^* \geq m)$.

Part~(ii) is proved like 
Corollary~\ref{almost surely bounded support of automorphisms implies bounded support of structure},
but with part~(i) instead of 
Corollary~\ref{corollary to automorphisms with large support are unusual}.
\hfill $\square$

\begin{cor}\label{having G as a subgroup almost allways implies bounded support of automorphisms in unlabelled case}
For every finite group $G$ there is $T \in \mbbN$ such that 
\[
\lim_{n\to\infty} 
\frac{\liso \{ \mcM \in \mbS_n : G \leq \Aut(\mcM) \text{ and } \spt^*(\mcM) \leq T \} \riso}
{\liso \{\mcM \in \mbS_n : G \leq \Aut(\mcM) \} \riso} \ = \ 1. 
\]
\end{cor}

\noindent
{\bf Proof.}
Let $G$ be isomorphic to a permutation group without fixed points on $[m]$ for some $m \in \mbbN^+$.
Let $t = 2r(m! - 1)m/m! + 1$.
In the same way as we proved
Corollary~\ref{having G as a subgroup almost allways implies bounded support of automorphisms},
but using Proposition~\ref{automorphisms with large support are unusual for unlabelled structures}
instead of Proposition~\ref{automorphisms with large support are unusual}, we get
\begin{equation*}
\lim_{n \to \infty} 
\frac{\liso \{\mcM \in \mbS_n : G \leq \Aut(\mcM) \text{ and } \spt(\mcM) \leq t) \} \riso}
{\liso \{ \mcM \in \mbS_n : G \leq \Aut(\mcM) \} \riso} \ = \ 1.
\end{equation*}
By Proposition~\ref{support of structure is bounded in terms of support of automorphisms}
the sought after $T$ exists.
\hfill $\square$

\begin{theor}\label{unlabelled limit laws}
For each result in the previous sections which, for some sequence $\mbS'_n \subseteq \mbS_n$,
$n \in \mbbN^+$, and
set $\mbX \subseteq \mbS$ that is closed under
isomorphism, can be stated in the form 
\[ \lim_{n\to\infty} \frac{\big| \mbS'_n \cap \mbX \big|}{\big| \mbS'_n \big|} \ = \ c \]
where $0 \leq c \leq 1$, we also have
\[ \lim_{n\to\infty} \frac{\liso \mbS'_n \cap \mbX \riso}{\liso \mbS'_n \riso} \ = \ c \]
for the same constant $c$.
\end{theor}

\begin{rem}\label{remark about unballed limit laws}{\rm
The statement in Theorem~\ref{unlabelled limit laws} that we get exactly the same limit
$c$ in both the labelled and unlabelled case may seem counter intuitive, because we consider
structures with a nontrivial automorphism. Roughly speaking, the reason why we indeed
get exactly the same limit in the labelled and the unlabelled case is that 
for each $\mbS' = \bigcup_{n \in \mbbN^+ }\mbS'_n$ considered, there is $p$ such that
$\big| \mbS'_n \big| \sim \big| \mbS'_n \cap \mbS_n(\spt^* \leq p) \big|$ and therefore
Lemma~\ref{relationship between counting labelled and unlabelled structures} can
be applied in the proof of Theorem~\ref{unlabelled limit laws}.
}\end{rem}

\begin{exam}{\rm
Here are three examples of applications of Theorem~\ref{unlabelled limit laws}.\\
(i) Let $t \geq 2$, let $\varphi$ be a sentence and let 
\[ \mbX_\varphi \ = \ \big\{ \mcM \in \mbS : \mcM \models \varphi \big\}. \]
By Theorem~\ref{limit law for first-order logic}, 
$\big| \mbS_n(\spt \geq t) \ \cap \ \mbX_\varphi \big| \Big/ \big| \mbS_n(\spt \geq t) \big|$
converges to some $0 \leq c \leq 1$ as $n \to \infty$. 
Now Theorem~\ref{unlabelled limit laws} implies that
\[ \lim_{n\to\infty} \frac{\liso \mbS_n(\spt \geq t) \ \cap \ \mbX_\varphi \riso}
{\liso \mbS_n(\spt \geq t) \riso} \ = \ c.\]
(ii) Let $G$ be a finite group, $\varphi$ a sentence and $\mbX_\varphi$ as above.
By Theorem~\ref{limit law for first-order logic}, 
\[ \frac{\big| \big\{ \mcM \in \mbS_n : G \leq \Aut(\mcM) \big\} \ \cap \ \mbX_\varphi \big|}
{\big| \big\{ \mcM \in \mbS_n : G \leq \Aut(\mcM) \big\} \big|} \]
converges to some $0 \leq c \leq 1$. Theorem~\ref{unlabelled limit laws} implies that
\[ \lim_{n\to\infty} 
\frac{\liso \big\{ \mcM \in \mbS_n : G \leq \Aut(\mcM) \big\} \ \cap \ \mbX_\varphi \riso}
{\liso \big\{ \mcM \in \mbS_n : G \leq \Aut(\mcM) \big\} \riso} \ = \ c.\]
}\end{exam}

\noindent
{\bf Proof of Theorem~\ref{unlabelled limit laws}.}
Suppose that $\mbS'_n \subseteq \mbS_n$,
for $n \in \mbbN^+$, that
$\mbX \subseteq \mbS$ is closed under
isomorphism and that we have proved (in previous sections) that
\begin{equation}\label{labelled convergence of X to c}
\lim_{n\to\infty} \frac{\big| \mbS'_n \cap \mbX \big|}{\big| \mbS'_n \big|} \ = \ c 
\end{equation}
for some $0 \leq c \leq 1$. 
In all of these cases it is clear that $\mbS' = \bigcup_{n=1}^{\infty} \mbS'_n$
is closed under isomorphism.
It also follows, either by definition or by 
results that have been proved,
that there is an integer $p$ such that
\begin{align}\label{almost all structures in S' have bounded support}
&\lim_{n\to\infty}
\frac{ \big| \mbS'_n \ \cap \ \mbS_n(\spt^* \leq p) \big| }{ \big| \mbS'_n \big| } \ = \ 1 
\quad \text{ and} \\
&\lim_{n\to\infty}
\frac{ \liso \mbS'_n \ \cap \ \mbS_n(\spt^* \leq p) \riso }{ \liso \mbS'_n \riso } \ = \ 1.
\label{almost all iso-classes of structures in S' have bounded support}
\end{align}
It follows from~(\ref{labelled convergence of X to c}) 
and~(\ref{almost all structures in S' have bounded support}) that 
\begin{equation}\label{labelled convergence of S' S-n and X}
\lim_{n\to\infty}
\frac{\big| \mbS_n(\spt^* \leq p) \ \cap \ \mbS' \ \cap \ \mbX \big|}
{\big| \mbS_n(\spt^* \leq p) \ \cap \ \mbS' \big|} \ = \ c. 
\end{equation}

\noindent
Lemma~\ref{relationship between counting labelled and unlabelled structures}
with $\mbY = \mbS'$ gives
\[ \big| \mbS_n(\spt^* \leq p) \ \cap \ \mbS' \big| \ \sim \ 
n! \liso \mbS_n(\spt^* \leq p) \ \cap \ \mbS' \riso, \]
and with $\mbY = \mbS' \cap \mbX$ it gives
\[ \big| \mbS_n(\spt^* \leq p) \ \cap \ \mbS' \ \cap \ \mbX \big| \ \sim \ 
n! \liso \mbS_n(\spt^* \leq p) \ \cap \ \mbS' \ \cap \ \mbX \riso. \]
This together with~(\ref{labelled convergence of S' S-n and X}) gives 
\begin{align*}
\frac{\liso \mbS_n(\spt^* \leq p) \ \cap \ \mbS' \ \cap \ \mbX \riso}
{\liso \mbS_n(\spt^* \leq p) \ \cap \ \mbS' \riso} 
\ &= \ 
\frac{n! \liso \mbS_n(\spt^* \leq p) \ \cap \ \mbS' \ \cap \ \mbX \riso}
{n! \liso \mbS_n(\spt^* \leq p) \ \cap \ \mbS' \riso} \\
&\sim \
\frac{\big| \mbS_n(\spt^* \leq p) \ \cap \ \mbS' \ \cap \ \mbX \big|}
{\big| \mbS_n(\spt^* \leq p) \ \cap \ \mbS' \big|} 
\ \to \ c \quad \text{as } n \to \infty.
\end{align*}
Combining this with~(\ref{almost all iso-classes of structures in S' have bounded support}) gives
\begin{align}\label{unlabelled convergence of S' S-n and X over S'}
&\frac{\liso \mbS_n(\spt^* \leq p) \ \cap \ \mbS'_n \ \cap \ \mbX \riso}{\liso \mbS'_n \riso} \\
= \ 
&\frac{\liso \mbS_n(\spt^* \leq p) \ \cap \ \mbS'_n \ \cap \ \mbX \riso}
{\liso \mbS_n(\spt^* \leq p) \ \cap \ \mbS'_n \riso} 
\ \cdot \ 
\frac{\liso \mbS_n(\spt^* \leq p) \ \cap \ \mbS'_n \riso}{\liso \mbS'_n \riso} 
\ \to \ c \quad \text{ as } n \to \infty. \nonumber
\end{align}
Finally we have
\begin{align*}
\frac{\liso \mbS'_n \ \cap \ \mbX \riso}{\liso \mbS'_n \riso} 
\ = \ 
\frac{\liso \mbS_n(\spt^* \leq p) \ \cap \ \mbS'_n \ \cap \ \mbX \riso}{\liso \mbS'_n \riso}
\ + \ 
\frac{\liso \mbS_n(\spt^* \geq p+1) \ \cap \ \mbS'_n \ \cap \ \mbX \riso}{\liso \mbS'_n \riso}
\end{align*}
which tends to $c$ as $n \to \infty$, 
because of~(\ref{unlabelled convergence of S' S-n and X over S'})
and~(\ref{almost all iso-classes of structures in S' have bounded support}).
This concludes the proof of Theorem~\ref{unlabelled limit laws}.
\hfill $\square$

\end{document}